\documentclass[11pt]{amsart}
\usepackage{fullpage}
\usepackage{amssymb,amsmath}
\usepackage{amsthm, etoolbox}
\usepackage[x11names]{xcolor}
\usepackage{hyperref}
\usepackage{graphicx} 
\usepackage[utf8]{inputenc}
\usepackage[T1]{fontenc}
\usepackage{amsfonts}
\usepackage{amsthm}
\usepackage{comment}
\usepackage{tikz-cd}
\usetikzlibrary{cd}
\usepackage[export]{adjustbox}
\usepackage[toc,page]{appendix}
\usepackage{enumerate}   
\usepackage{caption}
\usepackage[skip=1ex, belowskip=2ex]{subcaption}
\usepackage{bm}
\usepackage{enumitem}
\usepackage{stmaryrd}

\usetikzlibrary{decorations.pathmorphing, decorations.markings}
\usepackage[
backend=biber,
style=ext-alphabetic,    
maxnames=4,
minnames=3,
url=true,
doi=true,
maxbibnames=99,
maxcitenames=99,
backref=true
]{biblatex}

\DeclareFieldFormat{postnote}{#1}

\AtEveryBibitem{
	\iffieldundef{doi}{}{
		\clearfield{url}
	}
}

\DeclareCiteCommand{\citeplain}
{\unspace\usebibmacro{prenote}}
{\printtext[bibhyperref]{%
		\iffieldundef{labelalpha}
		{\printfield{labelnumber}}
		{\printfield{labelalpha}}}}
{\multicitedelim}
{\usebibmacro{postnote}\unspace}

\allowdisplaybreaks

\makeatletter
\newcommand{\myitem}[1]{%
	\item[#1]\protected@edef\@currentlabel{#1}%
}
\makeatother

\AtEveryBibitem{%
	\iffieldundef{doi}{}{%
		\clearfield{url}%
	}%
	\iffieldundef{eprint}{}{%
		\clearfield{url}%
	}%
}

\newcommand{\pa}[1]{\left(#1\right)}
\newcommand{\brc}[1]{\left\{#1\right\}}
\newcommand{\brk}[1]{\!\left[#1\right]}
\newcommand{\ticI}{\tilde {\mathcal I}}
\newcommand{\tif}{\tilde f}
\newcommand{\tix}{\tilde x}
\newcommand{\tbf}{\textbf}
\newcommand{\cD}{\mathcal D}
\newcommand{\cR}{\mathcal R}
\newcommand{\gr}{\text{gr}}
\newcommand{\codim}{\text{codim}}
\newcommand{\cA}{\mathcal A}

\newcommand{\mf}{\mathfrak}
\newcommand{\maxinv}{\text{maxinv}}

\newcommand{\mx}{\mathbf x}
\newcommand{\my}{\mathbf y}
\newcommand{\inv}{\text{inv}}
\newcommand{\cF}{\mathcal F}
\newcommand{\cG}{\mathcal G}
\newcommand{\cJ}{\mathcal J}
\newcommand{\cI}{\mathcal I}
\newcommand{\cO}{\mathcal O}

\newcommand{\maxord}{\text{maxord}}
\newcommand{\ord}{\text{ord}}

\newcommand{\Z}{\mathbb Z}
\newcommand{\ra}{\rightarrow}

\newcommand{\Spec}{\text{Spec }}

\newcommand{\fr}{\mathfrak}
\newcommand{\call}{\mathcal}
\newcommand{\bb}{\mathbb}

\newcommand{\mz}{\mathbf z}
\newcommand{\mbf}{\mathbf}

\newif\ifqedbarused
\newcommand{\qedbar}{%
	\unskip\nobreak\hfill % Push content to the right
	\hspace{1em} % Add space before the bar (adjust ''1em'' as needed)
	\rule{3em}{0.4pt} % Horizontal bar
}
\newcommand{\qedbarhere}{%
	\qedbar%
	\global\qedbarusedtrue% Suppress automatic bar
}

\newcommand{\thicksquigarrow}{%
	\begin{tikzpicture}[baseline=-0.5ex]
		\draw[line width=1.5pt, decoration={
			snake,
			amplitude=2pt,
			segment length=6pt,
			post length=0.5pt, % Adjust post length to control the flat end
			pre length=0.5pt   % Adjust pre length to control the flat end
		}, decorate] (0,0) -- (0.8,0); % Squiggly part (shorter)
		
		\draw[line width=1.5pt, -{>[line width=2pt]}] (0.8,0) -- (1.0,0); % Flat end with arrowhead
	\end{tikzpicture}%
}

\usepackage{todonotes}

\usepackage{xcolor,hyperref}
\definecolor{darkred}{rgb}{.8,0,0}
\definecolor{tocolor}{rgb}{.1,.1,.1}
\definecolor{urlcolor}{rgb}{.3,.2,.6}
\definecolor{linkcolor}{rgb}{.1,.4,.5}
\definecolor{citecolor}{rgb}{.5,.2,.3}
\definecolor{gray}{rgb}{.8,.8,.8}

\hypersetup{
	colorlinks=true,
	urlcolor=urlcolor,
	linkcolor=linkcolor,
	citecolor=citecolor,
	pdfauthor = {Maxim Brais},
}

\usepackage{aliascnt}

\newcommand{\thdef}[2]{
	\newaliascnt{#1}{theorem}  
	\newtheorem{#1}[#1]{#2}
	\aliascntresetthe{#1}  
	\newtheorem*{#1*}{#2}
	\expandafter\newcommand\expandafter{\csname #1autorefname\endcsname}{#2}
}

\newtheorem{theorem}{Theorem}[subsection]
\newtheorem*{theorem*}{Theorem}
\thdef{conjecture}{Conjecture}
\thdef{problem}{Problem}
\thdef{lemma}{Lemma}
\thdef{corollary}{Corollary}
\thdef{proposition}{Proposition}

\makeatletter
\newtheorem*{rep@theorem}{\rep@title}
\newcommand{\newreptheorem}[2]{%
	\newenvironment{rep#1}[1]{%
		\def\rep@title{#2 \ref{##1}}%
		\begin{rep@theorem}}%
		{\end{rep@theorem}}}
\makeatother

\newreptheorem{theorem}{Theorem}

\theoremstyle{definition}
\thdef{definition}{Definition}
\thdef{notation}{Notation}
\thdef{question}{Question}
\thdef{remark}{Remark}
\thdef{example}{Example}
\thdef{construction}{Construction}
\thdef{caution}{Caution}

\AtEndEnvironment{remark}{%
	\ifqedbarused\else\qedbar\fi%
	\global\qedbarusedfalse% Reset toggle
}
\AtEndEnvironment{definition}{%
	\ifqedbarused\else\qedbar\fi%
	\global\qedbarusedfalse%
}
\AtEndEnvironment{example}{%
	\ifqedbarused\else\qedbar\fi%
	\global\qedbarusedfalse%
}
\AtEndEnvironment{notation}{%
	\ifqedbarused\else\qedbar\fi%
	\global\qedbarusedfalse%
}
\AtEndEnvironment{caution}{%
	\ifqedbarused\else\qedbar\fi%
	\global\qedbarusedfalse%
}
\AtEndEnvironment{construction}{%
	\ifqedbarused\else\qedbar\fi%
	\global\qedbarusedfalse%
}
\numberwithin{equation}{section}

\newtheorem{method}{Method}

\title{Streamlining resolution of singularities with weighted blow-ups}
\author{Maxim Jean-Louis Brais}\thanks{Department of Mathematics, Universität Bonn, \href{mailto:s37mbrai@uni-bonn.de}{\texttt{s37mbrai@uni-bonn.de}}.\\ This research was supported by the ``2024 Undergraduate Research Award'' of the Natural Sciences and Engineering Research Council of Canada, and by the ``Supplément aux bourses du 1er cycle en milieu académique du CRSNG'' of the Fonds de recherche du Québec}
\date{\today}

\begin{document}
\delimitershortfall=-1pt

\begin{abstract}
In 2019, Abramovich--Temkin--W\l{}odarczyk and McQuillan used weighted blow-ups to obtain very fast and functorial algorithms for resolution of singularities in characteristic zero. Recently, Abramovich--Quek--Schober simplified the construction of the centre of blow-up introduced by Abramovich--Temkin--W\l{}odarczyk in the case of plane curves by using the Newton graph of the defining function. Their work follows the line of Schober's previous polyhedral analysis of the Bierstone--Milman invariant. In this paper, we extend their graphical approach to varieties of arbitrary (co)dimension in characteristic zero. This yields a factorial reduction in complexity in comparison with Abramovich--Temkin--W\l{}odarczyk, as previously achieved by W\l{}odarczyk. Our approach builds on the formalism of weighted blow-ups via filtrations of ideals developed and used by Loizides--Meinrenken, Quek--Rydh and W\l{}odarczyk, and on its interplay with systems of parameters. All constructions and proofs\,---\,including that of resolution of varieties by Deligne--Mumford stacks\,---\,are self-contained.
\end{abstract}
\maketitle
\setcounter{tocdepth}{1}
\tableofcontents
% !TEX root = main.tex
\section*{Introduction}\label{sec: intro}
In 2019, \cite{abramovich2024functorial}, \cite{mcquillan2019fastfunctorialeasyresolution} and \cite{Marzo2019} developed algorithms for resolution of singularities in characteristic zero based on weighted blow-ups, thereby giving faster proofs of Hironaka's classical result \cite[Main Theorem I]{Hironaka1964:I}.
These new algorithms differ significantly from earlier ``classical'' desingularisation algorithms based on classical (non-weighted) blow-ups.

Classical resolution algorithms must rely on what is called ``history'' or ``memory''; that is, any given step of the procedure is informed by the \emph{previous} blow-ups. The equation $x^2=y^2z^2$ in $\bb A^3$, studied in \cite[Example 1.13]{Frühbis-Krüger2023}, provides an insightful example of how this need arises. The singular locus is given by $V(x,zy)$, i.e. the $y$- and $z$-axes of the $x=0$ plane meeting at the origin. We want to consider blow-ups at smooth centres, and so there are essentially three candidates: the origin, the $y$-axis and the $z$-axis. Since there is an involution interchanging the two axes, choosing to blow up any of them would require choosing at random, which would prevent canonicity of the resolution process. Therefore, there is no choice but to blow up the origin. In one of the standard blow-up charts, we obtain the same singularity, but one of the lines is now contained in the exceptional divisor. Hence, we are able to \emph{canonically} choose this line as the centre for the next blow-up, which visibly improves the singularities by removing this singular line. The following observation can be drawn from this example: classical algorithms cannot operate by blowing up the ``worst'' singularity, unless we adapt the meaning of the ``worst'' singularity to be dependent on the resolution process itself. Weighted blow-ups can be used to correct this.

After Hironaka proved resolution of singularities in characteristic zero (\cite{Hironaka1964:I,hir64p2}), successive efforts have been made to improve his techniques. \cite{villamayor1989constructiveness} and \cite{Bierstone1997} achieved canonicity of the resolution (i.e. that the resolution process is invariant under automorphisms). Moreover, \cite{Schwartz1992} and \cite{wlodarzcyk2005simple} achieved functoriality for smooth morphisms. The proof of resolution itself was also made simpler by explicitly constructing different algorithms. These algorithms used different types of \tbf{invariants}. Let $Y$ be a smooth variety. Given a singular subvariety $X\subset Y$ to be resolved, an invariant is (very roughly) a function on the closed points of $X$ taking values in a well-ordered set whose maximal level set determines the centre of blow-up. The goal is for the maximal (value of the) invariant to decrease on the proper transform of $X$ under the blow-up. As a consequence of the above discussion, invariants of classical algorithms must take history into account. One of these invariants, defined in \cite{Bierstone1997} and henceforth referred to as the \emph{Bierstone--Milman invariant}, has values of the form
\begin{equation}
\inv_X^{(\text{BM})}(p)=(a_1,s_1,\frac{a_2}{a_1},s_2,\frac{a_3}{a_2},\dots, s_{k-1},\frac{a_{k}}{a_{k-1}})\label{eq: BM inv}
\end{equation}
at a given point $p\in X$. Here, $a_1\leq\cdots\leq a_k$ are positive rational numbers and $s_1,\dots,s_{k-1}$ are non-negative integers that reflect the history of the resolution process. In particular, before the resolution process starts (often referred to as ``the year zero''), we have $s_1=\cdots=s_{k-1}=0$.

The discovery of \cite{abramovich2024functorial,mcquillan2019fastfunctorialeasyresolution}, building on results from \cite{Panazzolo2006} and \cite{McquillanPanazzolo}, is that by allowing (stack-theoretic) weighted blow-ups one can successfully use an invariant that has no history component. In particular, by removing the history component of the Bierstone--Milman invariant \eqref{eq: BM inv}, they consider the invariant
\begin{equation}
	\inv_X^{(\text{ATW})}(p)=\inv_X(p):=(a_1,\dots,a_k).\notag
\end{equation}
If $(a_1,\dots,a_k)$ is the maximal invariant on $X$, and $Z$ is the locus where this maximum value is attained, then the maximal invariant decreases on the proper transform of $X$ under the weighted blow-up of $Y$ at $Z$ with weights $d/a_1,\dots,d/a_k$ assigned to the normal directions in a manner that shall be made precise later. Here, $d$ is the least integer such that $d/a_i$ is integral for all $i\leq k$. The centre of blow-up\,---\,namely, $Z$ with the weights assigned appropriately to the normal directions\,---\,is called the \tbf{associated centre}. Since $\inv_X$ only takes into account the geometry of $X$ and ignores the history of the resolution, we may reasonably regard $Z$ as the ``worst'' singular locus. We shall see later that the level sets are locally closed; $\inv_X$ may thus be regarded as stratifying $X$ from the smooth locus to the worst singular locus. Note also that since weighted blow-ups provide more flexibility to adapt to the local geometry of the singularities, far fewer blow-ups are typically required to achieve resolution.

 In \cite{abramovich2024functorial}, the associated centre is defined recursively using a certain type of coefficient ideals. As constructed in their paper, the coefficient ideal of an ideal is obtained by raising each of its derivatives to a specified power. The exponents in this construction involve factorials intended to calibrate the ``weight'' of each derivative. As a consequence, they grow out of control quite rapidly; even the resolution of mild singularities becomes uncomputable once a relatively low threshold in dimension or degree is exceeded. Recently, in \cite{abramovich2024resolvingplanecurvesusing}, the construction of the associated centre was simplified in the case of plane curves in order to extend these resolution techniques to perfect fields of positive characteristic under some assumptions. In particular, the construction of the associated centre is entirely formulated in terms of the Newton polygon of the defining function, in line with the approach of \cite{Schober2013, schober2014polyhedralapproachinvariantbierstone} for the Bierstone--Milman invariant, and thereby avoids coefficient ideals and factorial growth in complexity. 
 
 In this paper, we present two methods for constructing the associated centre. The first works directly in parameters and is best suited for computations in a local ring. The second uses derivatives of ideals, which is optimal for globalising the construction. Our approach to the construction of the associated centre generalises that of \cite{abramovich2024resolvingplanecurvesusing} to singular varieties of any (co)dimension \emph{over a field of characteristic zero}. This results in a factorial reduction in complexity compared to \cite{abramovich2024functorial}. We note that presentations of the Bierstone--Milman invariant that do not exhibit factorial growth were already obtained by \cite{Bierstone1997} and \cite{Schober2013, schober2014polyhedralapproachinvariantbierstone}. \cite{wlodarczyk2023functorialresolutiontorusactions,Włodarczyk2023} also achieved an efficient presentation of \cite{abramovich2024functorial}, with which our perspective shares some conceptual similarities. 
\subsection{Our formalism: weightings, smooth centres, and marked centres}
Some formalism is needed to properly define weighted blow-ups. Let $Y$ be a smooth variety and let $Z$ be a smooth subvariety. Recall that the non-weighted blow-up of $Y$ at $Z$ is a modification of $Y$ which, in colloquial terms, replaces $Z$ by the projectivisation of its normal bundle in $Y$. In particular, we replace $Z$ by a projective bundle. The notion of weighted blow-up generalises that of a non-weighted blow-up by using a weighted projective bundle instead of a projective bundle. Therefore, at the centre of blow-up $Z$, we need to keep track of which ``normal direction'' has which weight. We can do this by locally using a sequence of parameters that cut out $Z$, and assigning a ``weight'' or ``degree'' to each parameter. 

Recall that the automorphism group of a  weighted projective space $\bb P(w_1,\dots,w_k)$ is given by the graded automorphisms of $\Bbbk[x_1,\dots,x_k]$ modulo scaling, where this ring has grading given by $\deg(x_i)=w_i$. In other words, a change of coordinates on a weighted projective space is a change of coordinates on the corresponding affine space that preserves the degrees. This phenomenon carries over directly to weighted blow-ups: a change of parameters cutting out $Z$ which preserves the degree of the ``lower terms'' results in the same blow-up (see \autoref{lem:equal-weightings} below). In fact, we will need a parameter-invariant way to package this information. Here, there are multiple equivalent candidates, and we choose the notion of a \tbf{weighting}, following its study in \cite{loizides2021differentialgeometryweightings}, which builds on ideas dating back to \cite{melrose1996differential}. \cite{quek2021weighted} and \cite{wlodarczyk2023functorialresolutiontorusactions,Włodarczyk2023} also investigate this notion specifically in the context of weighted blow-ups. Let $Y$ be a smooth variety. A weighting on $Y$ is a descending filtration 
\[
\cO_Y=\cF_0\supset\cF_1\supset\cF_2\supset\cdots
\]
 of ideals on $Y$ satisfying the following condition: there exist weights $w_1,\dots,w_k\in \bb Z_{> 0}$, where $k\leq \dim Y=n$, such that locally near each point of $Y$, there is a system of parameters $\mx=(x_1,\dots,x_n)$ for which $\cF_j$ is the ideal generated by monomials of total weight at least $j$, that is, 
\[
\cF_j=\pa{\mx^\beta:\sum_{i=1}^k\beta_iw_i\geq j},
\]
where given a multi-index $\beta=(\beta_1,\dots,\beta_n)$, $\mx^\beta$ denotes the monomial $x_1^{\beta_1}\cdots x_n^{\beta_n}$. We say that the system of parameters $\mx$ is \tbf{compatible} with $\cF_\bullet$. We shall refer to the closed subvariety $Z=V(\cF_1)$, which we assume to be non-empty, as the \tbf{underlying subvariety} of $\cF_\bullet$. In the context of weighted blow-ups, we always assume that $\gcd(w_1,\dots,w_k)=1$, and we call such a weighting a \tbf{smooth centre}.
 
Given a weighting $\cF_\bullet$, we attach a valuation $v_{\cF_\bullet}$ to each open set.  Expanding locally a function $f=\sum_{\beta}c_\beta\mx^\beta$ in compatible parameters, where the coefficients $c_\beta$ are elements of the ground field $\Bbbk$, we have 
\[
v_{\cF_\bullet}(f)=\min_{\beta:\;c_\beta\neq 0}\brc{\sum_{i=1}^k\beta_iw_i},
\]
so that the valuation may be regarded as a weighted order of vanishing along $Z$. The weighting can be recovered from the valuation via $\cF_j=\{f:v_{\cF_\bullet}(f)\geq j\}$.

A  \tbf{marked centre} is a pair $\cJ=(\cF_\bullet,d)$ where $\cF_\bullet$ is a smooth centre and $d$ is a positive integer called the \tbf{marking}. The marking will be used to measure the degree of tangency of the smooth centre to the singularity we are trying to resolve. If $\cF_\bullet$ has weights $w_1,\dots,w_k$ and locally compatible parameters $\mx$, we write $\cJ$ locally as $(x_1^{a_1},\dots,x_k^{a_k})$, where each $a_i=\frac{d}{w_i}$. Conversely, since we assume the weights to be coprime, the expression $(x_1^{a_1},\dots,x_k^{a_k})$ where $a_1,\dots,a_k\in\bb Q_{>0}$, corresponds to a unique marked centre. We say that $\inv(\cJ)=(a_1,\dots,a_k)$ is the \tbf{invariant} of the marked centre $\cJ$, and we \emph{always assume this sequence is non-decreasing}. $\cJ$ also has a valuation $v_{\cJ}$ defined by $v_{\cJ}=\frac{1}{d}v_{\cF_\bullet}$. This gives 
\[
v_{\cJ}(f)=\min_{\beta:\;c_\beta\neq 0}\brc{\sum_{i=1}^k\frac{\beta_i}{a_i}}.
\]

\subsection{The associated centre as the optimal quasi-homogeneous approximation}
Let $X\subset Y$ be the variety whose singularities are to be resolved. For simplicity, in this introduction, we assume that $X=V(f)$ is a hypersurface given by a principal ideal $(f)\subset \cO_Y$. In the body of the paper, we shall treat the general case. Recall that a function $f$ is quasi-homogeneous with respect to a system of parameters $\mx=(x_1,\dots,x_n)$ and non-negative integral weights $w_1,\dots,w_k$ if it is a linear combination of monomials of the same weighted degree $d$. Equivalently, in terms of weightings, each monomial $\mx^\beta$ in the expansion has $v_{\cF_\bullet}(\mx^\beta)=d$, where $\cF_\bullet$ is the weighting with compatible parameters $\mx$ and weights $w_1,\dots,w_k$. Working with the marked centre $\cJ=(\cF_\bullet,d)$, this is equivalent to each such monomial satisfying $v_{\cJ}(\mx^\beta)=1$. In general, it might be impossible to find a system of parameters in which every monomial of $f$ has valuation $1$. However, we may regard a marked centre $\cJ$ for which $v_{\cJ}(f)=1$ as a \emph{quasi-homogeneous approximation} of $f$: there are some quasi-homogeneous terms with valuation $1$, and higher terms with higher valuation.

Consider things locally at a closed point $p\in X$ in the singular locus. Thus, marked centres are taken at $p$, meaning that the filtrations are on the local ring $\cO_{Y,p}$.
One theorem at the cornerstone of \cite{abramovich2024functorial}\,---\,which appeared in a different formulation already in \cite{Bierstone1997} as explained in \autoref{rem: graph int thm} below\,---\,is that there is a unique marked centre $\cJ=(x_1^{a_1},\dots,x_k^{a_k})$ at $p$ such that $v_\cJ(f)=1$ and the sequence $(a_1,\dots,a_k)$ is maximal in the lexicographic order; see \autoref{thm: unique maximum} below.\footnote{Since we are dealing with sequences of rationals with different lengths, we add the rule that a truncated sequence is \emph{larger}.} $\cJ$ is the \tbf{associated marked centre} to $X\subset Y$ at $p$ and its underlying smooth centre is the \tbf{associated centre} to $X\subset Y$ at $p$. Our invariant function then has value
\[
\inv_{X\subset Y}(p)=\inv(p):=\inv(\cJ)=(a_1,\dots,a_k)
\]
at $p$. We may think of $\cJ$ as the ``optimal'' quasi-homogeneous approximation of $f$ at $p$.

\subsection{Blowing up the optimal quasi-homogeneous approximation improves the singularities}
Let us try to understand why the optimal quasi-homogeneous approximation is relevant to resolution of singularities.  To this end, we first ought to understand some key geometric properties of weighted blow-ups. These were first explained in 2020 by \cite{loizides2021differentialgeometryweightings} in the differential-geometric context, where weighted blow-ups are presented as a quotient of their \emph{weighted deformation space}\,---\,the differential-geometric version of the \emph{degeneration to the weighted normal bundle}. In 2021, W\l{}odarczyk introduced in his Oberwolfach lecture the concept of \emph{cobordant blow-ups} (now in the algebro-geometric context) which is based on the notion of \emph{birational cobordism} and emphasises birational aspects. This notion was further developed in 2022 in an earlier version of \cite{wlodarczyk2023functorialresolutiontorusactions}, wherein a proof of resolution of singularities based on $\bb G_m$-actions\,---\,which we adapt in this paper\,---\,was presented for the first time. Simultaneously, \cite{quek2021weighted} also developed weighted blow-ups in the algebro-geometric context,
using the degeneration to the weighted normal cone/bundle.\footnote{We thank Jaros\l{}aw W\l{}odarczyk for clarifying parts of this history.} All these constructions are defined using a suitable \emph{extended Rees algebra}, which we cover in \autoref{sec: blow-ups}. In this paper, we choose to adopt the latter degenerational viewpoint. Indeed, we want to emphasise that the proof that the invariant decreases on the blow-up is to be \emph{degenerated} to the weighted normal bundle, where there is a $\bb G_m$-action simplifying things, something which we now explain.

Consider the following situation: let $X$ be the variety over a field $\Bbbk$ of characteristic zero to be resolved and let $Z$ be the worst singular locus\,---\,that is, the locus of maximal invariant\,---\,which is itself smooth and assumed to be connected for simplicity. Suppose further that $X$ is embedded in a smooth variety $Y$ that carries a $\bb G_m$-action such that $X$ is an invariant subvariety. Finally, suppose that $Z=Y^{\bb G_m}$ is the fixed locus, and that for all $p\in Y$, its limit $\lim_{t\to 0}t\cdot p\in Z$ exists. Note that by the Bia\l{}ynicki-Birula stratification \cite[Theorem 4.3]{bialynicki1973some}, $Y$ is an affine space fibre bundle over $Z$. By Sumihiro's theorem \cite[Corollary 2]{Sumihiro1974EquivariantCompletion}, $Y$ is covered by affine open invariant subsets; we can therefore use Luna's étale slice theorem \cite[III. Lemme]{MSMF_1973__33__81_0} to get a local étale equivariant morphism to the tangent space of $Y$ at any point $p\in Z$. Hence, the action has a description in local parameters with certain non-negative ``canonical'' weights, that is, it induces a weighting $\cF_\bullet$ on $Y$ whose underlying subvariety is $Z$. The blow-up $\tilde Y$ of $Y$ with respect to $\cF_\bullet$ admits a nice geometric description: it is the total space of the tautological line bundle over the stack-theoretic quotient $\brk{\pa{Y\setminus Z}/\bb G_m}=:E$, i.e. the line bundle associated to the $\bb G_m$-torsor\footnote{Here, the exceptional divisor is defined as a \emph{stack-theoretic} quotient to avoid introducing quotient singularities in the process of resolution.} 
\[
Y\setminus Z\twoheadrightarrow E.
\]
The blow-up $\tilde Y$ is birational over $Y$ with exceptional divisor $E$, which is a fibre bundle over $Z$ whose fibre is a stack-theoretic version of weighted projective space. The proper transform $\tilde X$ of $X$ under the blow-down is the tautological line bundle over $\brk{\pa{X\setminus Z}/\bb G_m}$. It is straightforward to see why the singularities have improved: in the process of forming $\tilde X$, we have essentially \emph{excised} the worst locus $Z$. Of course, we do not always have such a $\bb G_m$-action at our disposal in resolution of singularities. However, our quasi-homogeneous approximation is sufficiently good so that we can degenerate $X$ to a certain (weighted) normal cone where such an action exists, while preserving enough structure of the local singularities.

Recall that the \tbf{degeneration to the normal cone} to $X$ at $Z$ is a flat family $BX\to \bb A^1$ whose fibre $BX_s$ at $s\neq 0$ is canonically isomorphic to $X$ and whose special fibre $BX_0$ is the normal cone to $X$ at $Z$; see \cite[Chapter 5]{Fulton1998}.\footnote{Note that \cite{Fulton1998} defines the degeneration to the normal cone as the blow-up of $Y\times \bb P^1$ at $X\times \{\infty\}$. Thus, defining the blowup using the degeneration to the normal cone may appear to be circular reasoning. However,  \emph{op.\ cit.} already noted that the degeneration to the normal cone could equally be defined using the extended Rees algebra\,---\,as done in the present paper\,---\,which avoids the circularity.}  $BX$ carries a $\bb G_m$-action: $BX_{s\neq 0}$ is equivariantly isomorphic to $X\times \bb G_m$, and $\bb G_m$ acts by dilation of the normal cone $BX_0$. The subvariety $Z \times \bb G_m\hookrightarrow BX_{s\neq 0}$ extends to an invariant subvariety $\mbf V:=Z\times \bb A^1\hookrightarrow BX$. The blow-up of $X$ at $Z$ is then the quotient
\[
\tilde X:=\pa {BX\setminus\mbf V}/\bb G_m.
\]

A fundamental insight of \cite{loizides2021differentialgeometryweightings,quek2021weighted}\,---\,and of \cite{wlodarczyk2023functorialresolutiontorusactions} in the context of cobordant blow-ups\,---\,is that this construction works equally well in the presence of an embedding $X\subset Y$ and a weighting $\cF_\bullet$ on $Y$ whose underlying subvariety is $Z$; the main difference is that the special fibre $BX_0$ is now a weighted version of the normal cone. The equations for $BX_0\subset BY_0$ are given by the leading quasi-homogeneous terms of the equations for $X\subset Y$ with respect to the weighting $\cF_\bullet$. The blow-up of $Y$ at the weighting $\cF_\bullet$ is the stack-theoretic quotient
\[
\tilde Y:= \brk{\pa{BY\setminus \mbf V}/\bb G_m},
\]
 and the proper transform of $X$ under the blow-up is the quotient
\[
\tilde X=\brk{\pa{BX\setminus \mbf V}/\bb G_m}.
\]
We will show in \autoref{sec: global associated centre} that the definition of the associated centre globalises, and that the underlying subvariety of the (global) associated centre is the worst singular locus $Z$. If we take our weighting $\cF_\bullet$ to be the associated centre, showing that the singularities of $\tilde X$ have improved will amount to proving that our invariant (and therefore our optimal quasi-homogeneous approximation) is canonical enough in two ways:
\begin{enumerate}
	\item if $\pi:Y'\to Y$ is smooth, then $\inv_{\pi^{-1}(X)\subset Y'}=\pi^*\inv_{X\subset Y}$ (see \autoref{prop: local functoriality}),\label{eq: inv func smooth}
	\item $\maxinv_{BX\subset BY}=\maxinv_{X\subset Y}$ and the locus where this maximal invariant is attained on $BX$ is $\mbf V$ (see \ref{eq: C}, \autoref{prop: reduction conditions} and \autoref{prop: invariant to normal cone}).\label{eq: can deg normal}
\end{enumerate}
Note that  \eqref{eq: inv func smooth} and \eqref{eq: can deg normal} are part of the \emph{resolution principle} of \cite[3.3.33]{wlodarczyk2023functorialresolutiontorusactions} and \cite[4.4.18]{Włodarczyk2023}; see also \cite{abramovich2025logarithmicresolutionsingularitiescharacteristic,abramovich2025resolutionsingularitiesdynamicalmathematician}.

Using the smooth and surjective morphism
\[
\pi: BY\setminus \mbf V\to \brk{
\pa{
BY\setminus \mbf V
}
/\bb G_m
}
=\tilde Y,
\]
we find that
\[
\maxinv_{\tilde X\subset \tilde Y}\stackrel{\eqref{eq: inv func smooth}}{=}\maxinv_{BX\setminus \mbf V\subset BY\setminus \mbf V}<\maxinv_{BX\subset BY}\stackrel{\eqref{eq: can deg normal}}{=}\maxinv_{X\subset Y},
\]
i.e. the singularities have improved on $\tilde X$. This discussion can thus be summarised by the following slogan:
\begin{quotation}
The singularities on the weighted normal cone $BX_0$ always improve upon blowing up $Z$ with canonical weights; hence, by degenerating to the normal cone, the same is true for any variety $X$.
\end{quotation}
\subsection{Method 1: constructing the associated centre in parameters}\label{subsec: method 1}
We now attempt to construct the associated marked centre using local parameters. We fix a closed point $p\in X=V(f)$ in the singular locus and consider marked centres at $p$. We want to find the unique marked centre $\cJ=(x_1^{a_1},\dots,x_k^{a_k})$ satisfying $v_{\cJ}(f)=1$ and whose invariant $(a_1,\dots,a_k)$ is lexicographically maximal. 

Our strategy is to successively determine each minimal weight $1/a_j$ along with a compatible parameter $x_j$. That is, we will recursively construct marked centres $\cJ^{(j)}=(x_1^{a_1},\dots,x_j^{a_j})$ for $j\leq k$ which approximate the associated marked centre $\cJ=(x_1^{a_1},\dots,x_k^{a_k})$. It is important to recall that we assume $a_1\leq \cdots\leq a_k$. For $j< k$, the marked centre $\cJ^{(j)}$ has $v_{\cJ^{(j)}}(f)< 1$, meaning that we do not yet have a quasi-homogeneous approximation. We may complete $x_1,\dots,x_j$ into a full system of parameters $x_1,\dots,x_j,y_{j+1},\dots,y_n$. We think of $x_1,\dots,x_j$ as the ``good parameters'', and of $y_{j+1},\dots,y_n$ as the ``bad parameters'', i.e. those which still need to be refined. For each $b\geq a_j$, define the marked centre $\cJ^{(j)}[b]:=(x_1^{a_1},\dots,x_j^{a_j},y_{j+1}^{b},\dots,y_n^{b})$. By construction, the marked centre $\cJ^{(j)}[a_j]$ is a quasi-homogeneous approximation, i.e. it satisfies $v_{\cJ^{(j)}[a_j]}(f)= 1$. However, it is not an ``optimal'' quasi-homogeneous approximation, inasmuch as we may be able to find a smaller weight $1/a_{j+1}$ which can be assigned to the bad variables $y_{j+1},\dots,y_n$ so that $\cJ^{(j)}[a_{j+1}]$ is also a quasi-homogeneous approximation. In fact, the next entry in the invariant is the maximal number $a_{j+1}$ such that $\cJ^{(j)}[a_{j+1}]$ satisfies $v_{\cJ^{(j)}[a_{j+1}]}(f)= 1$. Now, the ``bad parameters'' $y_{j+1},\dots,y_n$ may not be the optimal candidates for playing this minimising game in the next steps. Therefore, we need to find the next ``good parameter'' $x_{j+1}$. We achieve this by choosing $x_{j+1}$ to be an element of $\cO_Y$ whose zero locus is maximally tangent to $V(f)$ with respect to the current ``good parameters'' and their weights.

Let us illustrate this process with monomials. Suppose that we have found our $j$th approximation $\cJ^{(j)}=(x_1^{a_1},\dots,x_j^{a_j})$ and that $y_{j+1},\dots, y_n$ are our ``bad parameters'' which still need to be refined. We first expand
\begin{equation*}
	f=\sum_{\beta}c_\beta x_1^{\beta_1}\cdots x_j^{\beta_j}y_{j+1}^{\beta_{j+1}}\cdots y_n^{\beta_n},
\end{equation*}
where $c_\beta\in\Bbbk$.
 We now seek the greatest number $a_{j+1}$ such that for each monomial $\beta$ appearing in the expansion, 
\begin{equation}v_{\cJ^{(j)}[a_{j+1}]}\pa{
x_1^{\beta_1}\cdots x_j^{\beta_j}y_{j+1}^{\beta_{j+1}}\cdots y_n^{\beta_n}
}
\geq 1.\label{eq: mon val a_{j+1}}
\end{equation}
Let $\beta$ be a monomial appearing in the expansion. Assuming $\sum_{i=1}^j\frac{\beta_i}{a_i}<1$, finding the greatest rational number $b$ such that 
\[
v_{\cJ^{(j)}[b]}
\pa{
x_1^{\beta_1}\cdots x_j^{\beta_j}y_{j+1}^{\beta_{j+1}}\cdots y_n^{\beta_n}
}
=\sum_{i=1}^j\frac{\beta_i}{a_i}+\frac{\sum_{i=j+1}^n\beta_i}{b}\geq 1\]
reduces to solving the equation 
\[
\sum_{i=1}^j\frac{\beta_i}{a_i}+\frac{\sum_{i=j+1}^n\beta_i}{b}=1,
\]
so that 
\[
b=\frac{\sum_{i=j+1}^n\beta_i}{1-\sum_{i=1}^j\frac{\beta_i}{a_i}}.
\]
Therefore, the minimal number $a_{j+1}$ such that \emph{all} monomials appearing in the expansion satisfy \eqref{eq: mon val a_{j+1}} is
\begin{equation}
a_{j+1}:=\min_\beta\brc{
\frac{\sum_{i=j+1}^n\beta_i}{1-\sum_{i=1}^j\frac{\beta_i}{a_i}}:c_\beta\neq 0\text{ and } \sum_{i=1}^j\frac{\beta_i}{a_i}<1
}
,
\label{eq: first a_{j+1}}
\end{equation}
assuming this minimum exists (see \autoref{lem: existence minimum}). Next, we must find the new ``good parameter'' $x_{j+1}$ to assign the weight $1/a_{j+1}$. As mentioned, this is done by imposing certain weighted tangency conditions. We proceed as follows: take $\beta$ to be a minimiser of \eqref{eq: first a_{j+1}}, and choose an index $l\geq j+1$ with $\beta_l>0$. We then define 
\[
x_{j+1}:=\partial_{x_1}^{\beta_1}\cdots \partial_{x_j}^{\beta_j}\partial_{y_{j+1}}^{\beta_{j+1}}\cdots\partial_{y_l}^{\beta_l-1}\cdots\partial_{y_n}^{\beta_n}f.
\]
That this gives a correct choice of parameter will be proved by a careful analysis of how differential operators interact with valuations (see \autoref{prop: ind semi-ass}).

The game we play starts at step $0$ with $\cJ^{(0)}=()$ (i.e. there are no ``good parameters'' yet), and ends at step $k$ whenever $\cJ=\cJ^{(k)}=(x_1^{a_1},\dots,x_k^{a_k})$ satisfies $v_{\cJ}(f)=1$. \autoref{k-associated implies associated} ensures that this procedure terminates with $\cJ$ as the associated centre.

\subsection{Method 2: constructing the associated centre with derivatives of ideals}\label{subsec: method 2} \hyperref[subsec: method 1]{Method 1} sketched above constructs the associated centre at a point. However, we want a global procedure, which is able to find the locus of maximal invariant $(a_1,\dots,a_k)$, and compute the associated centre there.
Although it can be done, \hyperref[subsec: method 1]{Method 1} is not optimal for such a globalisation, as it is parameter-dependent. We can formulate things more invariantly using differential operators as follows.

Given an ideal $\cI$, its $m$th derivative $\cD^{\leq m}\cI$ is defined by applying to $\cI$ all differential operators of order at most $m$. We recall this in \autoref{subsec: diff op}. These successive derivatives stratify $Y$ according to the order of vanishing of $\cI$, and this simplifies the globalisation process.

Suppose we wish to compute the associated centre to $X=V(\cI)\subset Y$. Moreover, suppose that $x_1,\dots,x_j$ have been chosen  as the  first $j$ parameters of the associated marked centre and that the first $j$ entries in the invariant $a_1,\dots,a_j$ are already known.

We recursively define certain ideals $\cD[\beta]=\cD[\beta;x_1,\dots,x_j;\cI]$ for any multi-index $\beta$ of length at most $j+1$. If $\beta=(\beta_1)$ is a multi-index of length $1$, we let 
\[
\cD\brk{\beta}:=\cD^{\leq\beta_1}\cI.
\]
Otherwise, we let 
 \[
\cD\brk{
\beta
}
:=\cD^{\leq\beta_{l}}\pa{
\cD\brk{
\pa{
\beta_1,\dots,\beta_{l-1}
}
}
|_{V(x_{l-1})}
}.
\]
The next entry in the invariant is 
\begin{equation}
a_{j+1}:=\min_{\beta=(\beta_1,\dots,\beta_j)}\left\{\frac{\beta_{j+1}}{1-\sum_{i=1}^j\frac{\beta_i}{a_i}}: \sum_{i=1}^j\frac{\beta_i}{a_i}<1\text{ and }\cD[\beta]=\cO_{Y,p}\right\}.\label{eq: method 2}
\end{equation}
If $\beta$ is a minimiser in \eqref{eq: method 2}, then any (lift to $\cO_Y$ of an)
element of $\cD[(\beta_1,\dots,\beta_j,\beta_{j+1}-1)]$ that vanishes to order $1$ at $p$ can be taken as the next parameter for the associated centre at $p$. Such an element is called a \tbf{maximal contact element} at $p$.
\subsection{Overview of the paper} In \autoref{sec: weightings and marked centres}, we introduce the notions of weighting and marked centre. We derive many simple yet subtle facts about them which shall be used throughout the paper. In \autoref{sec: blow-ups}, we quickly cover weighted blow-ups and show that if the invariant and associated centre satisfy a certain set of properties, then resolution of singularities follows with ease. That these properties are satisfied will be shown throughout the remainder of the paper. In \autoref{sec:associated centre}, we construct the associated centre using \hyperref[subsec: method 1]{Method 1}, and show simultaneously that its invariant is greatest amongst marked centres $\cJ$ with $v_{\cJ}(\cI)= 1$. In \autoref{sec: global associated centre}, we prove the existence of a \emph{global} associated centre. We then use numerical inequalities to derive \hyperref[subsec: method 2]{Method 2} directly from \hyperref[subsec: method 1]{Method 1}, and sketch an algorithm globalising \hyperref[subsec: method 2]{Method 2}. Finally, in \autoref{sec: eff}, we compare our approach with that of \cite{abramovich2024functorial}.

\subsection{Further applications} The methods of \cite{abramovich2024functorial} have been successfully generalised to the logarithmic setting and to foliations; see \cite{quek2022logarithmic,wlodarczyk2023functorialresolutiontorusactions,Włodarczyk2023,Temkin2023,Abramovich_2024,abramovich2025logarithmicresolutionsingularitiescharacteristic,abramovich2025principalizationlogarithmicallyfoliatedorbifolds}. We expect the methods developed in this paper to find application there as well.

\subsection*{Acknowledgements} The author wishes to thank their former supervisor, Brent Pym, who introduced them to resolution of singularities via weighted blow-ups two years ago and has been offering much appreciated help and guidance throughout the many stages of this project.  The author thanks Jacques Hurtubise for co-supervising their honours thesis, which contained several of the ideas that ultimately led to this work. Finally, the author also thanks Dan Abramovich, Bernd Schober and Jarosław Włodarczyk for kindly reading earlier versions of this paper and providing thoughtful comments that helped to situate it within the extensive literature on resolution of singularities.
\subsection*{Conventions}
By ``\tbf{variety}'', we mean a scheme that is separated, reduced, and of finite type over a field $\Bbbk$ of characteristic zero, which we now fix for the remainder of the paper. Throughout, $Y$ will denote a smooth variety of dimension $n$. By a (local) ``\tbf{function}'' on a variety, we always mean a (local) section of its structure sheaf. Unless otherwise specified, we work in the Zariski topology. By ``\tbf{ideal}'', we always mean a coherent ideal sheaf. By a ``\tbf{point}'' on a variety, we mean a closed point, and so $p\in Y$ denotes a closed point. We may also use the notation $|Y|$ to denote the set of closed points when specifically needed.

%  	
% !TEX root = main.tex
\section{Weightings and marked centres}\label{sec: weightings and marked centres}

Following \cite{loizides2021differentialgeometryweightings}, we adapt the notion of weighting for use in algebraic geometry. This structure is needed to define weighted blow-ups in a parameter-free manner and is equivalent to many other constructions used to define centres of weighted blow-ups; see \autoref{rem: terminology}.
\subsection{Definitions and first properties}
\begin{definition}[Parameters]\label{parameters}
	Let $U\subset Y$ be an open subvariety. A collection of elements $x_1,\dots, x_k\in\cO(U)$ is a \tbf{sequence of parameters} (or, more simply, \tbf{parameters}) if the exterior product $dx_1\wedge\cdots\wedge dx_k$ is a non-vanishing section of $\bigwedge^k\Omega_{U/\Bbbk}$. This is equivalent to the map $(x_1,\dots,x_k):U\ra \bb A^k$ being smooth. If $k=n$, we say that $x_1,\dots, x_n$ form a \tbf{system of parameters}. This is equivalent to the map $(x_1,\dots,x_n):U\ra \bb A^n$ being étale. If $x_1,\dots, x_k$ are parameters vanishing at $p\in Y$, we say that they are \tbf{parameters at $p$}, and similarly, if $k=n$, we say that they form a \tbf{system of parameters at $p$}.
\end{definition}
\begin{remark}
	Note that in contrast to coordinates (the analogous notion in differential and complex geometry), parameters need not be one-to-one. For example, consider the holomorphic function $z\mapsto z^m$ for $m\geq 2$. Its differential is $mz^{m-1}dz$, and so it is a parameter on $\bb C^*$, where it is $m$-to-one. Nevertheless, in the analytic topology, parameters are locally bijective, and therefore locally are coordinates. 
\end{remark}	
\begin{notation}
	Throughout, unless specified otherwise, bold Latin letters shall denote systems of parameters on $Y$ (hence shall \emph{always be of length $n$}), and the subindexed unbolded letter shall indicate the parameters themselves. For instance, we have $\mx=(x_1,\dots, x_n)$. 
	Similarly, Greek letters shall denote multi-indices with non-negative integral entries, and the subindexed letters shall denote the entries. Their length is at most $n$ and depends on the context; when specification is needed, we write  $l(\beta)=j$ to indicate that $\beta=(\beta_1,\dots,\beta_j)$ has length $j$. If $j=n$, we write $\mx^\beta$ for the monomial $x_1^{\beta_1}\cdots x_n^{\beta_n}$.
\end{notation}

\begin{definition}[Weighting]\label{def: weighting}
	A \tbf{weighting} on $Y$ is a descending filtration of ideals on $\cO_Y$
	\[
	\call O_Y=\call F_0\supset \call F_1\supset \call F_2\supset \cdots
	\]
	such that there exist weights $w_1,\dots w_k\in\Z_{>0}$ for $0\leq k\leq n$ and a non-empty closed smooth subvariety $Z\subset Y$ of pure codimension $k$ satisfying the following condition:
	\begin{quote}
		For all $p\in Y$, there exists a neighbourhood $U$ of $p$ and a sequence of parameters $x_1,\dots, x_k$ on $U$ that cut out $Z\cap U$, such that $\call F_j|_U$ is generated by those monomials $x_1^{\beta_1}\cdots x_k^{\beta_k}$ satisfying $\sum_{i=1}^k \beta_iw_i\geq j$.
	\end{quote}
	We say $Z$ is the \tbf{underlying subvariety} of $\cF_\bullet$. We consistently denote the underlying subvariety by $Z$, and may write $Z_{\cF_\bullet}$ when we find it convenient. Notice that $Z=V(\cF_1)$. 
	
	We say that the sequence $x_1,\dots,x_k$ or any system of parameters $\mx$ extending it is \tbf{compatible} with the weighting $\cF_\bullet$, or more succinctly \tbf{$\cF_\bullet$-compatible}. Note that we assume an order on the sequences, i.e. that the $i$th parameter should correspond to the weight $w_i$. As we shall see, a weighting has multiple compatible systems of parameters.
	
	By convention, we always assume that $w_1\geq w_2\geq\cdots\geq w_k$. We call $(w_1,\dots,w_k)$ the \tbf{weight sequence} of $\cF_\bullet$. We say a weighting is a \tbf{smooth centre} if the weights satisfy $\gcd(w_1,\dots, w_k)=1$.
	
	For a point $p\in Y$, we call a filtration on the local ring $\cO_{Y,p}$ a \tbf{weighting} (\emph{resp}. \tbf{smooth centre}) \tbf{at $p$} if it is the stalk at $p$ of a weighting (\emph{resp.} smooth centre) $\cF_\bullet$ on an open $U\subset Y$ with $p\in Z_{\cF_\bullet}$.
\end{definition}
\begin{remark}\label{conormal filtration} The weights are independent of the presentation of the weighting, and so the weight sequence is well-defined. Indeed, consider the filtration 
	\[
	\pa{
		\call N_{\cF}^\vee
	}
	_\bullet:=\frac{\cF_\bullet+\cF_1^2}{\cF_1^2}
	\]
	on the conormal sheaf of $Z$. Notice that $(\call N_\cF^\vee)_j=0$ for all $j> w_1$ since $\cF_j$ has no linear monomials. Furthermore, we also see that $(\call N_\cF^\vee)_j\supsetneq(\call N_\cF^\vee)_{j+1}$ if and only if $\cF_j$ has a linear monomial which is not in $\cF_{j+1}$. This happens exactly when $j=w_i$ for some $i\leq k$. Thus, the weights $w_1,\dots, w_k$ are exactly the degrees $j$ for which the associated graded piece $\text{gr}_j(\call N_\cF ^\vee)$ is non-zero. The rank of $\gr_j(\call N_{\cF}^\vee)$ is the number of times that the weight $j$ appears in the weight sequence. \end{remark}

\begin{remark}
	One can define weighted blow-ups in greater generality by dropping the requirement that $(x_1,\dots,x_k)$ be parameters, as done in \cite{quek2021weighted}. This is analogous to blowing up a scheme at a non-smooth subscheme, and accordingly exhibits more complicated geometry. We do not need this level of generality to develop our theory of desingularisation; whence our choice to work with \emph{smooth} centres, as defined above.
\end{remark}

\begin{example}\label{ex: zero weighting}
	We define the \tbf{zero weighting} to be the weighting given by the filtration 
	\[\cO_Y\supset 0\supset 0\supset\cdots.\]
	This is indeed a weighting since it is locally given by the empty sequence of parameters and the empty weight sequence. Its underlying subvariety is $Y$ itself. It will be useful in many of our inductive proofs.
\end{example}
\begin{example}\label{ex: gen in deg 1} Let $\mathfrak m_p$ be the maximal ideal at $p\in Y$. Consider the filtration given by 
	\[\cO_Y=\fr m_p^0\supset\fr m_p^1\supset\fr m_p^2\supset\cdots.\]
	Any system of parameters at $p$ is compatible with this weighting, and the weights are all $1$.
\end{example}
\begin{example}\label{ex: weighting associated to cusp}
	Let $Y=\bb A^2$. The filtration
	\[
	\cO_Y\supset\pa{
		x,y
	}
	\supset\pa{
		x,y
	}
	\supset \pa{
		x,y^2
	}
	\supset\pa{
		x^2,xy,y^2
	}
	\supset \pa{
		x^2,xy,y^3
	}
	\supset \pa{
		x^2,xy^2,y^3
	}
	\supset\cdots
	\]
	is a weighting with weights $(3,2)$ and compatible parameters $(x,y)$. One may check that the systems of parameters $(x+y^p,x^q+y)$ for $p\geq 2$ and $q\geq 1$ are also compatible; thus, this weighting admits many compatible systems of parameters.
\end{example}
\begin{proposition}\label{basic fact weighting} Let $\cF_\bullet$ be a weighting with weights $w_1,\dots,w_k$. Then, for all $j\geq 0$, $\cF_j/\cF_{j+1}$ is a locally free $\cO_Z$-module of rank equal to the number of multi-indices $\beta$ of length $k$ such that $ \sum\beta_i w_i=j$.
\end{proposition}
\begin{proof}
	Note that the sheaf $\cF_j/\cF_{j+1}$ is supported on $Z$. On a neighbourhood $U$ of $p\in Z$ on which $\cF_\bullet$ has compatible parameters $x_1,\dots, x_k$, the monomials $x_1^{\beta_1}\cdots x_k^{\beta_k}$ with $\sum\beta_iw_i=j $ modulo $\cF_{j+1}$ provide a free $\cO_{Z\cap U}$-basis of $\cF_j/\cF_{j+1}$.
\end{proof}
Weightings conveniently pull back along smooth morphisms.
\begin{proposition}\label{prop: pulling back weightings}
	Let $f:Y'\ra Y$ be a smooth map and $\cF_\bullet$ be a weighting on $Y$ such that the image of $Y'$ intersects $Z_{\cF_\bullet}$. Then, the pullback $\cF_\bullet\cdot\cO_{Y'}=f^*\cF_\bullet$ (defined filtered-piece-wise) is a weighting on $Y'$.
\end{proposition}
\begin{proof}
	$\cF_\bullet$-compatible parameters on $Y$ pull back to $(\cF_\bullet\cdot\cO_{Y'})$-compatible parameters on $Y'$.
\end{proof}
\subsection{Valuation and conditions for compatible parameters}\label{sub: valuation and condition for compatible parameters}

To a weighting $\call F_\bullet$, we associate a valuation on ideals and functions on $Y$. For an open $U\subset Y$ and an ideal $\cI\subset\cO_U$, we let 
\[
v_{\cF_\bullet}(\cI)=\max\brc{
	j:\cI\subset \cF_j
}
=\min\brc{
	j:\cI\neq 0\mod \cF_{j+1}
},
\]
where we set by convention this quantity to be infinite if the maximum/minimum does not exist. We should think of this valuation as a weighted order of vanishing of $\cI$ with respect to the weighting. If $\cI$ is the germ of an ideal at $p$, we define $v_{\cF_\bullet}(\cI)$ to be the limit of the valuation of $\cI$ evaluated on neighbourhoods of $p$. For a function $f\in\cO_Y$, we define its valuation to be the valuation of the ideal $(f)$.

One can equivalently define $v_{\cF_\bullet}$ as a section of the sheaf of valuations on (the Zariski-Riemann space of) $Y$; see \cite[Section 2]{abramovich2024functorial} for details.

\begin{remark}\label{rem: val recovers filt} Note that
	$\cF_j=\{f:v_{\cF_\bullet}(f)\geq j\}$ holds on any open, so we may recover the weighting from the valuation.
\end{remark}
The valuation may also be understood in parameters. A choice of compatible parameters $x_1,\dots,x_k$ on an open $U$ gives a trivialisation of the normal bundle of $Z\cap U=V(x_1,\dots,x_k)$ in $U$. Assuming $U$ is affine, we may choose a formal tubular neighbourhood: we have a (non-uniquely determined) isomorphism
\[
\widehat{\cO_U}:=\varprojlim_i\;\cO_U/\cA^i\simeq \cO_{Z\cap U}[\![x_1,\dots,x_k]\!]
\]
where $\cA=(x_1,\dots,x_k)$; see \cite[Corollaire 0.19.5.4]{grothendieck1964elements}. In particular, after embedding $\cO_U$ in $\widehat{\cO_U}\simeq \cO_{Z\cap U}[\![x_1,\dots,x_k]\!]$, we may expand any function $f\in\cO_U$ as 
\begin{equation}
	f=\sum_{l(\beta)=k} g_{\beta} x_1^{\beta_1}\cdots x_k^{\beta_k},\label{eq: expanding f}	
\end{equation}
where the coefficients $g_{\beta}$ are functions on $Z\cap U$ pulled back to $U$ by extending them constantly in the normal directions.

\begin{proposition}\label{valuation and expansions1} Let $(w_1,\dots,w_k)$ be the weights of $\cF_\bullet$. For $f\in\cO_U$ expanded as in \eqref{eq: expanding f}, we have
	\[
	v_{\cF_\bullet}(f)=\min_{\beta=(\beta_1,\dots,\beta_k):\;g_\beta\neq 0}\brc{
		\sum_{i=1}^k \beta_iw_i
	}
	.
	\]
\end{proposition}

\begin{proof}
	Let $c$ be the above minimum. Taking $f\mod \cF_c$ amounts to removing all the monomials $x_1^{\beta_1}\cdots x_k^{\beta_k}$ from $f$ which have $\sum_{i=1}^k\beta_i w_i\geq c$. Minimality ensures that $f\equiv 0\mod\cF_c$. However, since $f$ has at least one non-zero monomial with valuation $c$, we have that $f\not\equiv 0\mod \cF_{c+1}$. 
\end{proof}
Note that $v_{\cF_\bullet}$ is indeed a valuation in the traditional algebraic sense.
\begin{corollary}\label{cor: valuation is valuation}
	For any functions $f,g$ on $U$, $v_{\cF_\bullet}$ satisfies the following:
	\begin{itemize}
		\item $v_{\cF_{\bullet}}(f)=\infty$ if and only if $f= 0$;
		\item $v_{\cF_\bullet}(fg)=v_{\cF_\bullet}(f)+v_{\cF_\bullet}(g)$ and 
		\item $v_{\cF_\bullet}(f+g)\geq\min\{v_{\cF_\bullet}(f),v_{\cF_\bullet}(g)\}$, with equality if $v_{\cF_\bullet}(f)\neq v_{\cF_\bullet}(g).$
	\end{itemize}
\end{corollary}
\begin{proof}
	Each property follows directly by expanding locally $f$ and $g$ in compatible parameters.
\end{proof}
As a further corollary of \autoref{valuation and expansions1}, we can characterise the valuation using any system of parameters $\mx$ at $p\in U$ completing $x_1,\dots,x_k$.
\begin{corollary}\label{valuation and expansions}
	Expand $f=\sum_{\beta} c_\beta\mx^\beta$ where each $c_\beta\in\Bbbk$. Then, at $p$, we have 
	\[
	v_{\cF_\bullet}(f)=\min_{\beta=(\beta_1,\dots,\beta_n):\;c_\beta\neq 0}\brc{
		\sum_{i=1}^k\beta_iw_i
	}
	.
	\]
\end{corollary}
The valuation is determined by evaluating at a single point in $Z_{\cF_\bullet}$.
\begin{proposition}\label{prop: val from pt to global}
	Let $f\in\cO_Y$ and suppose that $Z$ is connected. Let $p\in Z$. If $v_{\cF_\bullet}(f)\geq j$ at $p$, then $v_{\cF_\bullet}(f)\geq j$ globally on $Y$.
\end{proposition}
\begin{proof} Let $U$ be a connected open set where we have a formal tubular neighbourhood and expand $f$ as in \eqref{eq: expanding f}. The valuation at $p$ is determined by whether $g_\beta=0$ near $p$. Since $Z\cap U$ is integral, $g_\beta= 0$ near $p$ if and only if $g_\beta=0$, so that we obtain $v_{\cF_\bullet}(f)\geq j$ at all points $q\in Z\cap U$. Covering $Z$ by other charts like $U$ and repeating the same argument (connectedness ensures that they must intersect with $Z\cap U$), we obtain $v_{\cF_\bullet}(f)\geq j$ at every $q\in Z$. Away from $Z$, the inequality is trivial.
\end{proof}

The next lemma provides a way to compare weightings using the valuation.
\begin{lemma}\label{lem: conditions for inclusion}
	Let $\cF_\bullet$ and $\cG_\bullet$ be two weightings on $Y$ and let $\cG_\bullet$ have compatible parameters $(x_1,\dots, x_k)$ on an open $U\subset Y$ with weights $(w_1,\dots,w_k)$. Then, $\cG_\bullet\subset\cF_\bullet$ on $U$ if and only if $v_{\cF_\bullet}(x_i)\geq w_i$ for $1\leq i\leq k$.
\end{lemma}
\begin{proof}
	It suffices to observe that $\cG_\bullet\subset \cF_\bullet$ on $U$ if and only if for all $f\in\cO_U$, $v_{\cF_\bullet}(f)\geq v_{\cG_\bullet}(f)$.
\end{proof}

Unlike the weights, the compatible parameters $(x_1,\dots,x_k)$ are not unique, as seen in \autoref{ex: weighting associated to cusp}. 
Nevertheless, we can give necessary and sufficient conditions for a sequence of parameters to be compatible with a weighting in terms of the valuation. We first need this lemma:

\begin{lemma}\label{lem: inclusion implies locally equal}
	Let $\call G_\bullet$ and $\call F_\bullet$ be two weightings on $Y$ with connected underlying subvarieties and the same weights $w_1,\dots,w_k$. If $\call G_\bullet \subset \call F_\bullet$, then $\call G_\bullet=\call F_\bullet$. 
	\begin{proof}
		Since $\call G_1\subset\call F_1$, and both ideals define a connected smooth variety with the same dimension, it follows that $\call G_1=\call F_1$. Now suppose $\call G_j=\call F_j$ for $j\geq 1$. Since $\call G_{j+1}\subset \call F_{j+1}$, we have a surjection $\call G_j/\call G_{j+1}=\call F_j/\call G_{j+1}\twoheadrightarrow\call F_j/\call F_{j+1} $. Since the weights of $\call G_\bullet$ and $\call F_\bullet$ agree, \autoref{basic fact weighting} implies both $\call G_j/\call G_{j+1}$ and $\call F_j/\call F_{j+1}$ are finite rank locally free sheaves on $Z$ of the same rank.  Therefore, the above surjection must be an isomorphism, which implies $\call G_{j+1}=\call F_{j+1}$.
	\end{proof}
\end{lemma}

The necessary and sufficient conditions are the following.

\begin{proposition}
	\label{lem:equal-weightings}
	Let $\cF_\bullet$ be a weighting with connected underlying subvariety and weight sequence $(w_1,\dots, w_k)$. Then, a sequence of parameters $(x_1,\dots,x_k)$ for which $V(x_1,\dots,x_k)$ is connected is compatible with $\cF_\bullet$ if and only if $v_{\cF_\bullet}(x_i)=w_i$ for all $i\leq k$.
\end{proposition}
\autoref{lem:equal-weightings} generalises \cite[Lemma~5.2.10]{abramovich2024functorial} and is equivalent to \cite[Corollary~2.3.11]{temkin2025dreamresolutionprincipalizationi}.
\begin{proof}
	That the valuations must be as specified is immediate from the definition of a weighting.  For the converse, suppose $v_{\cF_\bullet}(x_i)=w_i$ and
	let $\cG_\bullet$ be the weighting given by parameters $(x_1,\dots,x_k)$ and weights $(w_1,\dots, w_k)$. By \autoref{lem: conditions for inclusion}, we have that $\cG_\bullet\subset \cF_\bullet$, and so, \autoref{lem: inclusion implies locally equal} implies that $\cG_\bullet=\cF_\bullet$.\end{proof}
Note in particular that when considering weightings at a point, the connectedness hypothesis on $Z$ and $V(x_1,\dots,x_k)$ may be removed as we are considering germs of smooth varieties.

\autoref{lem:equal-weightings} ought to be thought of as a description of the trivialisations and transitions of weightings; indeed, the structure group of a weighting is given by those parameter-transformations that preserve the \emph{valuation}, i.e. the weighted degree of the weighted leading terms. In particular, weightings may be seen as generalising $\bb G_m$-actions of positive weights, for which the structure group consists of those parameter-transformations that preserve the \emph{degree} (i.e. preserve homogeneity).\footnote{It should be noted here that a weighting with degree-preserving parameter-transformations as transitions need not be induced from a $\bb G_m$-action; indeed, the parameters need not generate the algebra of functions.} This observation is useful to retain for \autoref{sec: blow-ups}.

We will need a technical result regarding the existence of certain compatible parameters.

\begin{proposition}\label{prop: nested param}
	Let $\cF_\bullet$ be a weighting with connected underlying subvariety and with weights $w_1,\dots, w_k$, and let $j\leq k$. Suppose that there is a system of parameters $\mx$ such that $v_{\cF_\bullet}(x_i)=w_i$ for all $i\leq j$, which is \emph{not} necessarily $\cF_\bullet$-compatible. Then, there exists an $\cF_\bullet$-compatible system of parameters $\my$ such that $y_i=x_i$ for $i\leq j$, and such that for $i\geq j+1$, $y_i$ is a function of $x_{j+1},\dots,x_n$ only in its $\mx$-expansion.
\end{proposition}
\begin{proof}
	By \autoref{lem:equal-weightings}, we may find an $\cF_\bullet$-compatible system of parameters $\mz$ with $z_i=x_i$ for $i\leq j$. For $i\geq j+1$, we expand $z_i=\sum_{l(\beta)=n}c_\beta^{(i)}\mx^\beta$. We define the set
	\[
	M:=\brc{\beta=(\beta_1,\dots,\beta_n):\beta_1=\cdots=\beta_j=0\text{ or } v_{\cF_\bullet}(\mx^\beta)\leq w_{j+1}},
	\]
	and observe it is finite. Consider 
	\begin{equation}
	y_i:=\sum_{\beta\in M}c_\beta^{(i)}\mx^\beta.\label{eq: exp of yi}
	\end{equation}
	This is a regular algebraic function (as $M$ is a finite set). Moreover, we note that 
	\[
	v_{\cF_\bullet}(z_i-y_i)=v_{\cF_\bullet}\pa{\sum_{\beta\notin M}c_\beta^{(i)}\mx^\beta}>w_{j+1}\geq w_i=v_{\cF_\bullet}(z_i)
	\]
	as for each $\beta\notin M$, $v_{\cF_\bullet}(\mx^\beta)>w_{j+1}$ and $v_{\cF_\bullet}$ is a valuation. Therefore, using again that $v_{\cF_\bullet}$ is a valuation, we conclude that
	\[
	v_{\cF_\bullet}(y_i)=v_{\cF_\bullet}(z_i-(z_i-y_i))=w_i.
	\]
	For $i\leq j$, we define $y_i:=x_i$. Using that $\mz$ is a system of parameters, we note that $\my$ is indeed also a system of parameters: for $i\geq j+1$, $y_i|_{V(x_1,\dots,x_j)}$ and $z_i|_{V(x_1,\dots,x_j)}$ have the same linear terms in their $\mx$-expansion. Thus, $\my$ is a compatible system of parameters by \autoref{lem:equal-weightings}, and it satisfies the required conditions of the proposition.
\end{proof}
\subsection{Differential operators}\label{subsec: diff op}
\begin{notation}
	Recall that given a parameter $x_i$, the definition of the partial derivative $\partial_{x_i}$ depends on a whole system of parameters $\mathbf x$. Usually, it is clear from the context in which system one works, but many of our proofs will need to account for this subtlety. When necessary, we will write $\partial_{\mathbf x_i}$ to denote the $x_i$ derivative in the $\mathbf x$-system and $\partial_{\mathbf x}^{\beta}$ to denote $\partial_{\mathbf x_1}^{\beta_1}\cdots\partial_{\mathbf x_n}^{\beta_n}$ where $\beta$ is a multi-index of length $n$. We also write $|\beta|$ for $\sum_{i=1}^n\beta_i$.
\end{notation}
Given a weighting $\cF_\bullet$ and an integer $m\geq 0$, we obtain a filtration on the sheaf $\cD_{Y/\Bbbk}^{\leq m}=\call D^{\leq m}$ of differential operators of degree $\leq m$ on $Y$.
Using \autoref{rem: val recovers filt}, we can define this filtration by specifying the valuation on differential operators.
For $\Phi\in \cD^{\leq m}$, we define locally
\[
v_{\cF_\bullet}(\Phi)=\min_{f\in \cO_Y}\brc{
	v_{\cF_\bullet}\pa{
		\Phi\pa{
			f
		}
	}
	-v_{\cF_\bullet}\pa{
		f
	}
}
.
\]
In other words, $v_{\cF_\bullet}(\Phi)$ is the minimal amount (possibly negative) by which applying $\Phi$ moves a function’s position in the filtration. The filtration $\cF_\bullet^{\leq m}$ on $\cD_Y^{\leq m}$ is then defined locally by
\[
\cF_j^{\leq m}=\brc{\Phi\in \cD_Y^{\leq m}:v_{\cF_\bullet}(\Phi)\geq j}.
\]
Let us shrink $Y$ and assume that it has a compatible system of parameters $\mx$, and let the weights be $w_1,\dots,w_k$. Observe that the (non-commutative) $\cO_Y$-algebra of differential operators $\cD=\cD_{Y/\Bbbk}$ is generated over $\cO_Y$ by $\partial_{x_i}$, for $i\leq n$. The valuation $v_{\cF_\bullet}$ on $\cD$ extends that on $\cO_Y$, and to understand it, it suffices to specify what it does on \emph{these} generators: 
\[v_{\cF_\bullet}(\partial_{x_i})=\begin{cases}
	-w_i & i\leq k\\ 0& i>k.
\end{cases}\]
The following proposition, which is an analogue of  \autoref{valuation and expansions}, follows directly. As a corollary, the filtration $\cF_\bullet^{\leq m}$ starts in degree $-mw_1$.

\begin{proposition}\label{prop: val and exp for diff op}
	Let $\Phi\in\cD^{\leq m}$, and expand $\Phi=\sum_{|\beta|\leq m}g_\beta\partial_\mx^\beta$.
	Then, 
	\[
	v_{\cF_\bullet}(\Phi)=\min_{|\beta|\leq m}\brc{
		v_{\cF_\bullet}(g_\beta)-\sum_{i=1}^n\beta_iw_i
	}
	.
	\]
\end{proposition}

Let $\cI\subset\cO_Y$ be an ideal. For $\Phi\in \cD$, define the ideal $\Phi(\cI)$ to be the ideal generated by differentiated sections of $\cI$, i.e.
\[
\Phi(\cI)=\pa{
	\Phi(f): f\in\cI
}
\]
We can reduce statements on the valuation of $\cI$ to statements on the valuation of certain differential operators.

\begin{proposition}[Monomial-differential-operator duality]\label{prop: duality} Let $\mx$ be $\cF_\bullet$-compatible parameters at $p\in Y$.  Then, the following quantities are equal, where the valuation and ideals are taken at the stalk at $p$.
	\begin{enumerate}
		\item $v_{\cF_\bullet}(\cI)$;
		\item $\min\{-v_{\cF_\bullet}(\Phi):\Phi\in\cD, v_{\cF_\bullet}(\Phi(\cI))=0\}$;
		\item $\min\{-v_{\cF_\bullet}(\Phi):\Phi\in\cD, \Phi(\cI)=(1)\}$;
		
		\item $\min\{-v_{\cF_\bullet}(\partial_\mx^\beta): v_{\cF_\bullet}(\partial_\mx^\beta(\cI))=0\}$;
		\item $\min\{-v_{\cF_\bullet}(\partial_\mx^\beta): \partial_\mx^\beta\cI=(1)\}$.
	\end{enumerate}
\end{proposition}
\begin{proof}
	Clearly, we have $(2)\leq (3),(4)$ and $(3),(4)\leq (5)$. So it suffices to show $(1)\leq (2)$ and $(5)\leq (1)$.
	
	$(1)\leq (2)$: Suppose $v_{\cF_\bullet}(\Phi(\cI))=0$. We may find $f\in\cI$ such that $v_{\cF_\bullet}(\Phi(f))=0$. We then have that $v_{\cF_\bullet}(\Phi)=\min_g\{v_{\cF_\bullet}(\Phi(g))-v_{\cF_\bullet}(g)\}\leq -v_{\cF_\bullet}(f)\leq -v_{\cF_\bullet}(\cI)$. 
	
	$(5)=(1)$: $\partial_\mx^\beta\cI=(1)$ if and only if there is $f\in\cI$ with the monomial $\mx^\beta$ appearing in its expansion with respect to $\mx$, and so \autoref{valuation and expansions}, implies the equality. 
\end{proof}
\begin{remark}
	Working with monomials is more intuitive, whereas working with differential operators yields cleaner proofs, since it avoids repeatedly specifying expansions. This is one reason to work with differential operators rather than monomials, but it is nevertheless instructive to retain both perspectives.
\end{remark}	
One may ask what can be said of $v_{\cF_\bullet}(\partial_\mx^{\beta})$ when $\mx$ is not a compatible system of parameters. If the first $j$ parameters have the right valuations, then the following technical lemma will be useful.

\begin{lemma}\label{lem: fdmtl chain rule lmma}
	Suppose $Z_{\cF_\bullet}$ is connected and let $\mx$ be a system of parameters with $v_{\cF_\bullet}(x_i)=w_i$ for all $i\leq j$, where $j\leq n$. Then, we have 
	\[
	v_{\cF_\bullet}(\partial_\mx^{\beta})\geq -\sum_{i=1}^j\beta_iw_i-\sum_{i=j+1}^n\beta_iw_{j+1}.
	\]
\end{lemma}	
\begin{proof}
	Choose an $\cF_\bullet$-compatible system of parameters $\my$ satisfying the conditions of \autoref{prop: nested param}. Using the chain rule, we have that $\partial_{\mx_i}=\partial_{\my_i}$ for $i\leq j$ and $\partial_{\mx_i}=\sum_{l=j+1}^n(\partial_{\mx_i}y_l)\partial_{\my_l}$ for $i>j$. Therefore, we have $v_{\cF_\bullet}(\partial_{\mx_i})=-w_i$ for $i\leq j$ and $v_{\cF_\bullet}(\partial_{\mx_i})\geq-w_{j+1}$ for $i>j$, implying the statement.
\end{proof}
Let $\cI\subset\cO_Y$ be an ideal. We denote by $\cD^{\leq m}\cI$ the ideal sheaf generated by the image of the evaluation map
\[
\mathrm{ev}:\cD^{\leq m}\times\cI\ra\cO_Y.
\]
Note that by the chain rule, locally where $\mx$ is a system of parameters, $\cD^{\leq m}\cI$ is generated by the ideals $\partial_\mx^\beta\cI$ for $|\beta|\leq m$. It follows that $\cD^{\le m}\cI$ vanishes precisely at those points where $\cI$ has order of vanishing greater than $m$.
\begin{definition}[Maximal contact]
	Let $a=\ord_p\;\cI$ be the order of vanishing of $\cI$ at $p\in V(\cI)$, i.e.
	\[
	a=\max\brc{i:\fr m^i\supset \cI}
	\]
	where $\fr m$ is the maximal ideal at $p$. Let $x\in\cD^{\leq a-1}\cI$. If $H=V(x)$ is smooth at $p$, we say that it is a \tbf{maximal contact hypersurface} at $p$, and that $x$ is a \tbf{maximal contact element} of $\cI$ at $p$.
\end{definition}
A useful heuristic is to regard maximal contact hypersurfaces as smooth hypersurfaces maximally tangent to $V(\cI)$ at $p$. They have been used in resolution of singularities since \cite{Giraud74}. In characteristic zero, they always exist locally.\footnote{This is where these techniques for resolution of singularities break in positive characteristic.} If $H$ is a maximal contact hypersurface at $p$, then, because $H$ is smooth near $p$, it remains a maximal contact hypersurface at every point in some neighbourhood of $p$ within $V(\cD^{\leq a-1}\cI)$. Using quasi-compactness, we can therefore find \emph{finitely} many hypersurfaces $H_1,\dots, H_m$ such that for each $p$ with $\ord_p\;\cI=a$, there is an $H_i$, for $i\leq m$, that is a maximal contact hypersurface at $p$. 
\subsection{Marked centres}
\begin{definition}[Marked centre]\label{marked centre}
	A \tbf{marked centre} on $Y$ is a pair $\call J=(\call F_\bullet,d)$, where $\call F_\bullet$ is a smooth centre on $Y$ and $d\in\bb Z_{>0}$. The integer $d$ is referred to as the \tbf{marking} of $\cF_\bullet$. If $w_1,\dots, w_k$ are the weights of $\call F_\bullet$, we call $\inv(\cJ)=(a_1,\dots a_k):=(d/w_1,\dots, d/w_k)$ the \tbf{invariant} of $\call J$. We say that a sequence (or system) of parameters is \tbf{compatible} with $\cJ$ when it is compatible with $\cF_\bullet$. 	 	\end{definition}
Notice that by our conventions, the invariant of a marked centre is always a non-decreasing sequence. 	
To a marked centre $\cJ=(\cF_\bullet,d)$, we again define a valuation on local functions and ideals by setting $v_{\cJ} := \frac{1}{d}\,v_{\cF_\bullet}$. Note that the value group of $v_\cJ$ is $\frac{1}{d}\,\Z$. If $\mx$ is a compatible system of parameters at $p$, then by expanding $f=\sum_{\beta}c_\beta\,\mx^\beta$ and applying \autoref{valuation and expansions}, we obtain
\[
v_{\cJ}(f)=\min_{\beta:\;c_\beta\neq 0}\brc{\sum_{i=1}^k\frac{\beta_i}{a_i}}
\]
at $p$.
\begin{definition}\label{def: pre-invariant} A \tbf{pre-invariant} is a non-empty non-decreasing sequence of positive rational numbers of length at most $n=\dim Y$.
\end{definition}
\begin{notation}
	For a pre-invariant $(a_1,\dots,a_k)$ and parameters $x_1,\dots,x_k$ on $U$, we will write $(x_1^{a_1},\dots x_k^{a_k})$ to denote the marked centre $\call J=(\call F_\bullet,d)$, where $\call F_\bullet$ is the smooth centre on $U$ given by parameters $x_1,\dots,x_k$ and weights $w_i=d/a_i$, and where $d$ is the unique positive integer such that $\text{gcd}(w_1,\dots, w_k)=1$. Since we require the underlying weighting of a marked centre to be a smooth centre (i.e. that the weights have no common divisor), every marked centre whose invariant has positive length locally has this form. By convention, we say that the \hyperref[ex: zero weighting]{zero weighting} with marking $1$ is the \tbf{zero marked centre}, which we denote by $()$.\end{notation}
\begin{example}
	Consider parameters $(x,y)$ on $\bb A^2$ and the marked centre $(x^2,y^3)$. Note that $d=6$ is the unique positive integer making $w_1=d/2=3$ and $w_2=d/3=2$ coprime. Therefore, the underlying smooth centre of this marked centre is the one defined in \autoref{ex: weighting associated to cusp}, with marking $6$.
\end{example}
\begin{remark}[Terminology in the literature]\label{rem: terminology}
	Our notion of \emph{marked centre} corresponds to the notion of \emph{center} in \cite{abramovich2024functorial} and to the notion of \emph{weighted centre} in \cite{temkin2025dreamresolutionprincipalizationi,temkin2025dreamresolutionprincipalizationii}. The latter is a $\mathbb{Q}$-ideal (or an \emph{idealistic exponent}, following \cite{Hir77,villamayor1989constructiveness,Hironaka2003,Hironaka2005,Schober2013,SCHOBER2021353}, or a \emph{presentation}, following \cite{Bierstone1997}) that locally has a presentation in parameters. Put simply, a $\mathbb{Q}$-ideal is a formal rational power of an ideal, \emph{modulo} identifying ideals with the same integral closure. The bridge between idealistic exponents and filtration of ideals appears in \cite{Kawanoue2007} and \cite[3.5.4--3.5.6]{quek2021weighted}. We have chosen the name ``marked centre'' because of its resemblance to the notion of a \emph{marked ideal} (mainly regarding how the marking works)\,---\,another notion equivalent to $\mathbb{Q}$-ideals and idealistic exponents\,---\,which has been used in resolution of singularities by \cite{wlodarzcyk2005simple,kollar2007resolutionsingularitiesseattle,BierstoneMilman2008}. Moreover, our notion of \emph{smooth centre} agrees with the notion of \emph{reduced center} of \cite{abramovich2024functorial} and with the notion of \emph{Rees algebra} (with a local presentation in parameters) in \cite{quek2021weighted}. Our choice of terminology is motivated by the desire for a ``centre'' to be that which is blown up and by the fact that weighted blow-ups do not depend on the marking.
\end{remark}

%
% !TEX root = main.tex
\section{Weighted blow-ups and resolution}\label{sec: blow-ups}
We will use smooth centres to define weighted blow-ups. As mentioned in \autoref{rem: terminology}, they offer a framework equivalent to the ones previously used in the literature. Our treatment of weighted blow-ups shall be terse, and we refer to \cite{loizides2021differentialgeometryweightings,wlodarczyk2023functorialresolutiontorusactions,quek2021weighted} for a more exhaustive exposition. After introducing weighted blow-ups, we will describe the key properties of the associated centre and the invariant that are required to achieve resolution, drawing in part on \cite{wlodarczyk2023functorialresolutiontorusactions,Włodarczyk2023},  \cite{abramovich2025logarithmicresolutionsingularitiescharacteristic} and \cite{abramovich2025resolutionsingularitiesdynamicalmathematician}.
\subsection{Weighted blow-ups}
Recall that the blow-up of a scheme $Y$ at an ideal $\cF_1$ is defined in \cite[Chapter 7]{Hartshorne1977} as 
\[
\text{Proj}_Y\pa{
	\cO_Y\brk{
		\cF_1 T
	}
}
:=\pa{
	\text{Spec}_Y\cO_Y\brk{
		\cF_1 T
	}
	\setminus V\pa{
		\cF_1 T
	}
}
/\bb G_m.
\]
The graded algebra $\cO_Y[\cF_1 T]$ is called the \tbf{Rees algebra} of the blow-up. Since we are removing the zero section $V(\cF_1 T)$ (i.e. the zero locus of the irrelevant ideal), we may as well consider the \tbf{extended Rees algebra} $\cO_Y[\cF_1 T, T^{-1}]$, which is a $\Z$-graded algebra. Indeed, the two projective spectra coincide: 
\[
\text{Proj}_Y\pa{
	\cO_Y\brk{
		\cF_1 T, T^{-1}
	}
}
:=\pa{
	\text{Spec}_Y\cO_Y\brk{
		\cF_1 T,T^{-1}
	}
	\setminus V\pa{
		\pa{
			\cF_1 T
		}
	}
}
/\bb G_m=\text{Proj}_Y\pa{
	\cO_Y\brk{
		\cF_1 T
	}
}
,
\]
where $(\cF_1T)\subset\cO_Y[\cF_1T,T^{-1}]$ is the ideal \emph{generated} by $\cF_1T$.  The extended Rees algebra was first introduced by Rees himself in \cite{Rees1961}, where it is referred to as the \emph{$\cF_1$-transform of $\cO_Y$.} Working with the extended Rees algebra has the following advantage: $\text{Spec}_Y\cO_Y[\cF_1 T,T^{-1}]$ is the degeneration to the normal bundle of $V(\cF_1)\subset Y$ (or degeneration to the normal cone if the embedding $V(\cF_1)\hookrightarrow Y$ is not regular), as defined in \cite[Chapter 5]{Fulton1998}. This provides insight into the geometry of blow-ups and also eases the process of deriving chart equations.

The extended Rees algebra can equally be constructed in the presence of weights, as developed by \cite{loizides2021differentialgeometryweightings,wlodarczyk2023functorialresolutiontorusactions,quek2021weighted}. It may therefore be used to define weighted blow-ups, as we now recall. Let $Y$ be a smooth variety and let $\cF_\bullet$ be a smooth centre on $Y$. The \tbf{extended Rees algebra of $\cF_\bullet$} is the $\Z$-graded $\cO_Y$-algebra
\begin{flalign*}
	&&&&&& \cR = \cR_{\cF_\bullet} := \bigoplus_{j\in\Z} \cF_j T^j && \text{where $\cF_j: = \mathcal{O}_Y$ for $j < 0$.}
\end{flalign*}
In terms of generators, this algebra may be written as
\[
\cR=\cO_Y\brk{
	T^{-1},\cF_jT^j:j\geq 0
}
=\cO_Y\brk{
	T^{-1},\cF_jT^j: 0\leq j\leq w_1
}
\subset \cO_Y\brk{
	T^{
		\pm 1
	}
}
\]
where $w_1$ is the highest weight.
Note that in the case where $\cF_\bullet$ has weights $(1,\dots,1)$, i.e. when $\cF_j=\cF_1^j$ for all $j\geq 0$, this agrees with the extended Rees algebra in the non-weighted setting. Let us now localise on $Y$ to where $\cF_\bullet$ has compatible parameters $x_1,\dots,x_k$ with weights $w_1,\dots,w_k$ and let $\tilde x_i:=x_iT^{w_i}$. Notice that the irrelevant ideal $\cR_+\subset \cR$ (i.e. the ideal \emph{generated} by positive degree terms) is equal to $(\tilde x_1,\dots,\tilde x_k)$. We refer to $\mbf V:=V(\cR_+)\subset \text{Spec}_Y\cR$ as the \tbf{vertex}.

Letting $s=T^{-1}$, we have a deformation map $\text{Spec}_Y\cR\ra\Spec\Bbbk[s]=\bb A^1$. We call 
\[
B=B_{\cF_\bullet}:=\text{Spec}_Y\cR
\]
the \tbf{degeneration to the weighted normal bundle} of $Z_{\cF_\bullet}=Z$ with respect to $\cF_\bullet$. Note that the graded structure of $\cR$ induces a $\bb G_m$-action on $B$. In \cite{wlodarczyk2023functorialresolutiontorusactions,Włodarczyk2023}, this $\bb G_m$-space is also called the \emph{full cobordant blow-up}; see \autoref{rem:cob-b-u}.

To justify our terminology, let us look at the fibres. At $s\neq 0$, we have $\cR[s^{-1}]=\cR[T]=\cO_Y[T^{\pm 1}]$ and so $B_{s\neq 0}=Y\times \bb G_m$. At the specialisation $s=0$, something noteworthy occurs:
\begin{equation}
	\cR/(s)=\frac{\bigoplus_{j\in\Z}\cF_j T^j}{\pa{T^{-1}}}=\bigoplus_{j\geq 0}\cF_j/\cF_{j+1}=\gr_{\cF_\bullet}\cO_Y.
	\label{eq: gr at special}
\end{equation}
We call
\[
B_0=\text{Spec}_Y\pa{
	\gr_{\cF_\bullet}\cO_Y
}
\]
the \tbf{weighted normal bundle} of $Z\subset Y$ with respect to $\cF_\bullet$, as it generalises the classical (non-weighted) normal bundle. Note that along the identification \eqref{eq: gr at special}, $\tix_i\mod (s)$ gets identified with $x_i\mod \cF_{w_i+1}$. $B_0$ is a fibre bundle over $Z$ with fibre $\bb A^k$, as it is locally trivialised by
\begin{equation}
	\gr_{\cF_\bullet}\cO_Y\simeq\cO_Z[\tix_1,\dots,\tix_k].\label{eq: equiv trans}
\end{equation}
 Since $(s)$ is a homogeneous ideal, $B_0\subset B$ is an invariant subvariety for the $\bb G_m$-action. The grading on $\gr_{\cF_\bullet}\cO_Y$ is given, along the trivialisation \eqref{eq: equiv trans}, by $\deg(\tilde x_i)=w_i$ for all $i\leq k$, and we therefore see that the $\bb G_m$-action on $B_0$ is one of weights $w_1,\dots,w_k$ on each fibre. The action being globally defined on $B_0$, the transitions must be equivariant. In other words, the structure group sits inside the group of $\bb G_m$-equivariant automorphisms of affine space, where $\bb G_m$ acts with weights $w_1,\dots,w_k$. Consequently, in general, $B_0$ \emph{fails} to be a vector bundle: such automorphisms need not be linear unless $w_1=\cdots=w_k$. This kind of fibre bundle has many names in the literature: it is called a \emph{twisted weighted vector bundle} in \cite[Definition 2.1.3]{quek2021weighted} (not to be confused with the notion of \emph{twisted sheaf}), a \emph{weighted affine bundle} in \cite[Definition 3.2]{Arena_2025} and a \emph{$\bb G_m$-fibration} in \cite{bialynicki1973some}. In the differential-geometric setting, it is also called a \emph{graded bundle}; see \cite{Grabowski_2012,loizides2021differentialgeometryweightings}. 

To understand the $\bb G_m$-action elsewhere on $B$, one uses the fact that $x_1,\dots,x_k$ form a sequence of parameters to get an identification
\begin{equation}
	\cR=\cO_Y\brk{
		s,\tix_1,\dots,\tix_k
	}
	/\pa{
		x_i-s^{w_i}\tix_i:i\leq k
	}
	;
	\label{eq: Rext}
\end{equation}
see \cite[Proposition 5.2.2]{quek2021weighted}. The functions $s,\tix_1,\dots,\tix_k$ form a sequence of parameters on $B$ and the $\bb G_m$-action may be described as follows:
\[
\lambda\cdot\pa{
	s,\tix_1,\dots,\tix_k
}
=\pa{
	\lambda^{-1}s,\lambda^{w_1}\tix_1,\dots,\lambda^{w_k}\tix_k
}
.
\]
Note that \eqref{eq: Rext} shows directly that $\mbf V\simeq Z\times \Spec\Bbbk [s]$, and that $\bb G_m$ only scales $s$ there.
On $B_{s\neq 0}= Y\times\Spec\Bbbk[s^{\pm 1}]$, we have parameters $(s,x_1,\dots,x_k)$. The action of $\lambda$ on $x_i=s^{w_i}\tix_i$ is trivial, and so the $\bb G_m$-action on $B_{s\neq 0}$ also contracts towards $s=0$. At $s=0$, $\bb G_m$ acts on each fibre with weights $w_1,\dots,w_k$, as shown above.

Suppose now that $X\subset Y$ is a subscheme containing $Z$. We call the closure $BX\subset B$ of $X\times \bb G_m\subset Y\times \bb G_m=B_{s\neq 0}$ inside of $B$ the \tbf{degeneration to the weighted normal cone} to $X$ at $\cF_\bullet$. We have that
\[
BX=\text{Spec}_X\pa{
	\cR|_X
}
,
\]
where by the restriction to $X$ we mean restriction of the generators as ideals (rather than pullback), that is,
\[
\cR|_X=\cO_X\brk{
	T^{-1},\pa{
		\cF_j|_X
	}
	T^j:0\leq j\leq w_1
}
\subset \cO_X\brk{
	T^{
		\pm 1
	}
}.
\]
By definition, away from the specialisation, $BX_{s\neq 0}=X\times \bb G_m$. The specialisation
\[
BX_0=\text{Spec}_X\pa{
	\gr_{\cF_\bullet|_X}\cO_X
}
\]
is the \tbf{weighted normal cone to $X$} at $\cF_\bullet$.
We can use a commutative diagram to illustrate the situation.
\begin{equation*}
	\begin{tikzcd}
		Z\times\bb G_m\ar[dd,hook]\ar[r, equals]&\mbf V_{s\neq 0}\ar[dd,hook]\ar[r, hook]&\mbf V\ar[dd,hook]\ar[r,equals]&Z\times \bb A^1\ar[ddddddl,bend left=40]&\mbf V_0\ar[dd, hook, bend right]\ar[l, hook']\ar[r, equals]& Z\ar[dddd, hook, bend right=50, "\substack{\text{\tiny zero}\\\text{\tiny section}}" description]
		\\
		\\
		X\times \bb G_m\ar[dd,hook]\ar[r,equals, crossing over] &BX_{s\neq 0}\ar[dd, hook]\ar[r, hook, crossing over]&BX\ar[dd,hook]&&BX_0\ar[dd,hook]\ar[ll,hook', crossing over]\ar[uu, bend right]
		\\
		\\
		Y\times\bb G_m\ar[ddr]\ar[r,equals]& B_{s\neq 0}\ar[dd]\ar[r, hook]&B\ar[dd,"
		\substack{
			\text{\tiny degeneration}\\
			\text{\tiny map}
		}
		" description]&&B_0\ar[dd, "\substack{\text{\tiny special}\\\text{\tiny fibre}}" description]\ar[ll,hook', crossing over]\ar[r,equals]&\text{Spec}_Y\pa{\gr_{\cF_\bullet}\cO_Y}\ar[uuuu, bend right=50, "\substack{\text{\tiny twisted}\\\text{\tiny weighted}\\\text{\tiny vector}\\\text{\tiny bundle}}" description]
		\\
		\\
		&\bb G_m\ar[r, hook]&\bb A^1&&\{0\}\ar[ll, hook']
	\end{tikzcd}
\end{equation*}
$\bb G_m$ acts on each of these spaces, by pulling the spaces on each row towards the fibre at $0$. All the maps are equivariant.

\begin{notation}
	For a function $f\in \cO_Y$, we let $\tilde f:=f\cdot T^{v_{\cF_\bullet}(f)}$.
\end{notation}
Let $\cI$ be the ideal of $X\subset Y$. The ideal defining $BX\subset B$ is $\ticI:=\left (\tif:f\in\cI\right )$. Over the special fibre, we have
\[
\tilde\cI\mod (s)\simeq\gr_{\cF_\bullet\cap\cI}\cI\subset \gr_{\cF_\bullet}\cO_Y
,
\]
which we call the \tbf{initialisation} of $\cI$ with respect to $\cF_\bullet$.
Consider a compatible system of parameters $\mx$ which completes $x_1,\dots,x_k$. We may consider the system of parameters $\tilde\mx$ on $B_0$. If
\[
f=\sum_{\beta}c_\beta\mx^\beta\in\cO_Y,
\]
then\footnote{Here, by $\tilde f$, we tacitly mean $\tilde f\mod (s)$.}
\begin{equation*}
	\tilde f=\sum_{\substack{\beta\\v_{\cF_\bullet}(\mx^\beta)=v_{\cF_\bullet}(f)}}c_\beta\tilde\mx^\beta\in\cO_{B_0}=\gr_{\cF_\bullet}\cO_Y
	,
\end{equation*}
justifying the word ``initialisation'': we are taking the initial form, i.e. the (weighted) leading terms of $f$. The next lemma shall be used later. Its proof is immediate by expanding in $\mx$ and $\tilde\mx$.
\begin{lemma}\label{lem: diff and hom}
	Let $\alpha,\beta$ be multi-indices of length $n$ satisfying $\alpha_i\leq \beta_i$ for each $i\leq n$, $v_{\cF_\bullet}(\mx^\beta)=v_{\cF_\bullet}(f)$ and $c_\beta\neq 0$. Then, $\widetilde{\partial_\mx^\alpha f}=\partial_{\tilde\mx}^\alpha\tilde f$ (in $\cO_{B_0}=\gr_{\cF_\bullet}\cO_Y$).
\end{lemma}
Let us now take a step back. Recall that weighted projective space $\bb P(w_1,\dots,w_k)$ is not smooth in general since the $\bb G_m$-action with weights $w_1,\dots,w_k$ on $\bb A^k\setminus\{\mbf 0\}$ has a $\mu_{w_i}$-stabiliser in the $i$th chart, where $\mu_{w_i}$ is the group of $w_i$th roots of unity. Likewise, the (classical) weighted blow-up $\text{Proj}_Y\cR=(B\setminus\mbf V)/\bb G_m$ may be singular if the weights are not all equal to $1$. However, if we intend to use such blow-ups to resolve singularities, they would need to be smooth. This is where stacks enter the scene. We define the (stack-theoretic) \tbf{weighted blow-up} of $Y$ at the smooth centre $\cF_\bullet$ to be the stack quotient
\[
Bl_{\cF_\bullet}Y:=\call Proj_Y\cR:=\brk{
	\pa{
		B\setminus \mbf V
	}
	/\bb G_m
}
.
\]
Smoothness of $\bb G_m$ ensures that this stack is smooth. The stack-theoretic projective spectrum $\call Proj$ was first introduced in \cite{AbramovichHassett2010}, and we refer to \cite[10.2.7]{olsson2016algebraic} for a complete treatment. By $\bb G_m$-invariance, the weighted normal bundle $B_0=\text{Spec}_Y(\gr_{\cF_\bullet}\cO_Y)$ descends to the \tbf{exceptional divisor} of the blow-up
\[ E:=\call Proj_Y\pa{
	\gr_{\cF_\bullet}\cO_Y
}
=\brk{
	\pa{
		B_0\setminus Z
	}
	/\bb G_m
}
.
\]
This is a fibre bundle whose fibre is the \tbf{weighted projective stack}
\[\call P(w_1,\dots,w_k):=\left[\left(\bb A^k\setminus\{\mbf 0\}\right) /\bb G_m\right ],\]
where $\bb G_m$ acts on $\bb A^k$ with weights $w_1,\dots,w_k$. The coarse moduli space of $\call P(w_1,\dots, w_k)$ is the corresponding weighted projective space $\bb P(w_1,\dots,w_k)$: indeed, we may think of the stack $\call P(w_1,\dots,w_k)$ as a weighted projective space $\bb P(w_1,\dots,w_k)$ to which we have attached a $\mu_{w_i}$-stabiliser at the origin in the $i$th chart for all $i\leq k$.

Using \eqref{eq: Rext}, we can give chart equations for the blow-up $Bl_{\cF_\bullet}Y$. The $i$th chart is the $\tilde x_i=1$ ``slice''
\begin{equation}
	\brk{
		\pa{\text{Spec}_Y\frac{\cR}{(\tilde x_i-1)}
		}
		/\mu_{w_i}
	}
	=\brk{
		\pa{
			\text{Spec}_Y\frac{\cO_Y\brk{
					s,\tix_1,\dots,\widehat{\tix_i},\dots \tix_k
				}
			}{\pa{
					x_j-s^{w_j}\tix_j:j\neq i
				}
			}
		}/\mu_{w_i}
	}
	,
	\label{eq: chart DM}
\end{equation}
where $\mu_{w_i}$ acts in parameters (say over the algebraic closure for simplicity) via
\[
\zeta\cdot\pa{
	s,\tix_1,\dots,\widehat{\tix_i},\dots, x_n
}
=\pa{
	\zeta^{-1}s,\zeta^{w_1}\tix_1,\dots,\widehat{\zeta^{w_i}\tix_i},\dots,\zeta^{w_k}\tix_k,x_{k+1},\dots,x_n
}
.
\]
The stabilisers of $Bl_{\cF_\bullet}Y$ are $\mu_{w_1},\dots,\mu_{w_k}$, as exhibited in \eqref{eq: chart DM}. $Bl_{\cF_\bullet}Y$ is therefore a Deligne--Mumford stack; see \cite[Theorem~3.6.4]{alper2024stacks} for instance. We refer to 
\[
\brk{
	\pa{
		BX\setminus \mbf V
	}
	/\bb G_m
}
\]
as the \tbf{proper transform} of $X$ under the blow-up. The blow-up and the proper transform are birational over $Y$ and $X$ respectively, provided $Z$ has everywhere positive codimension in $Y$ and $X$ respectively. 
\begin{remark}[Coarse space]
	With schematic blow-ups, since taking Veronese subalgebras does not affect Proj, the blow-up at $\cF_1$ is the same as the blow-up at $\cF_1^j$ for all $j>0$. The same does not hold for $\call Proj$, the stack-theoretic Proj. Indeed, taking the degree $c$ Veronese subalgebra $\cR^{(c)}\subset\cR$ induces a $\mu_{c}$-gerbe $\call Proj_Y\cR^{(c)}\ra \call Proj_Y\cR^{(c)}\!\!\!\fatslash \mu_c=\call Proj_Y\cR$  (see \cite[Example~3.1.8]{quek2021weighted}, \cite[Exercise~6.4.24.(b)]{alper2024stacks} and \cite[Example 6.4.56]{alper2024stacks}). Nevertheless, the coarse space of $\call Proj$ is the classical relative Proj (see \cite[Proposition 1.6.1]{quek2021weighted}). One can show that any Veronese subalgebra $\cR^{(c)}\subset\cR$ of sufficiently divisible degree $c$ is generated in degree one, i.e. it corresponds to the schematic blow-up at $\cF_c$ (see \cite[Proposition 1.6.3]{quek2021weighted}, and \cite[Theorem 1.18 and Theorem 1.20]{muller2025amplesheavesweightedprojective} for sharp bounds on $c$). Therefore, the coarse moduli space of a weighted blow-up is always a schematic blow-up. 
\end{remark}
\begin{remark}[Base change]\label{rem: base change}
	The notions of weighting and smooth centres have been defined for varieties, but as mentioned in \autoref{prop: pulling back weightings}, they pull back along smooth (surjective) morphisms. Therefore, it is natural to define a weighting on an Artin stack via a weighting on a smooth atlas whose data descends to the stack. Equivalently, one may define it as a filtration of ideal sheaves on the stack that pulls back to an ordinary weighting along any smooth (surjective) morphism from a variety.
	
	Likewise, the relative stack-theoretic Proj over $Y$ satisfies base change, so one can define weighted blow-ups for Artin stacks. Indeed, this is essential for our algorithm: we need to perform multiple blow-ups, and the first blow-up already yields a Deligne--Mumford stack.
\end{remark}
\begin{remark}[Cobordant blow-ups]\label{rem:cob-b-u}
	We have presented weighted blow-ups using the degeneration to the normal bundle, which is a deformational object over the affine line. There exists another useful perspective grounded in the birational aspects of $\bb G_m$-actions. Given a nice enough $\bb G_m$-action on a space $B$, we have a birational equivalence $B_+/\bb G_m\dashrightarrow B_-/\bb G_m$, where
	\begin{align*}
		&B_+:=\{p:\lim_{\lambda\to \infty}\lambda\cdot p\text{ does not exist}\}\text{ and},\\
		&B_-:=\{p:\lim_{\lambda\to 0}\lambda \cdot p\text{ does not exist}\}.
	\end{align*}
	This construction was developed in \cite{Wlodarczyk2000} under the name of \emph{birational cobordism}. It was subsequently used to prove the weak factorisation theorem; see \cite{wlodarczyk2002toroidalvarietiesweakfactorization,AbramovichKaruMatsukiWlodarczyk2002}.
	
	In our situation, the space $B=\text{Spec}_Y\cR$ is a birational cobordism, and the induced birational equivalence (which is actually a birational morphism) $B_+/\bb G_m\to B_-/\bb G_m$ is the blow-down as defined above. In this framework, \cite{wlodarczyk2023functorialresolutiontorusactions, włodarczyk2023coxringsmorphismsresolution, Włodarczyk2023} refers to $B$ as the (full) \emph{cobordant blow-up} of $Y$ (with respect to $\cF_\bullet$).
\end{remark}
\subsection{What is needed for resolution?}\label{subsec: what is needed for res}
In this section, we show that given some properties of the associated centre and invariant, the proof that the invariant decreases on the blow-up can be degenerated to the weighted normal bundle. This proof was first presented in \cite{wlodarczyk2023functorialresolutiontorusactions}, and recently also appeared in \cite{Włodarczyk2023}, \cite{abramovich2025logarithmicresolutionsingularitiescharacteristic} and \cite{abramovich2025resolutionsingularitiesdynamicalmathematician}.

Suppose that we have a rule which for any finite type separated $\Bbbk$-scheme\footnote{We \emph{do not} assume that $X$ is a variety, since we will need to work with the invariant and associated centre on the weighted normal cone, which is not reduced in general.} $X$ embedded in a smooth variety $Y$ in positive and pure codimension,\footnote{Pure codimension eases things for functoriality.} returns a function on closed points 
\[
\inv_{X\subset Y}=\inv_X:|X|\ra \Gamma^{(m)}
\]
such that
\begin{enumerate}[label=(\roman*)]
	\item $\Gamma^{(m)}$ is a well-ordered set that depends only on the codimension $m$ of $X$ in $Y$;\label{eq: (i)}
	\item $\inv_X(p)=\min\:\Gamma^{(m)}$ if and only if  $X$ is smooth at $p$.\label{eq: (ii)}
	\item $\inv_X$ has a maximum $\maxinv_X:=\max_{p\,\in X}\inv_X(p)$, and the locus $Z$ where this maximum is attained is closed, smooth, and has pure codimension;\label{eq: (iii)}
\end{enumerate}

The function $\inv_X$ is called the \tbf{invariant}, and it should stratify the space from the smooth locus to the worst singularities. Suppose furthermore that to each such $X\subset Y$, we assign a smooth centre $\cF(X\subset Y)_\bullet=\cF_\bullet$ called the \tbf{associated centre} to $X\subset Y$, such that $Z_{\cF_\bullet}\subset X$ is the locus where $\inv_X(p)=\maxinv_X$. We now examine which properties the invariant and the associated centre must satisfy to ensure that repeated blow-ups at the associated centre will eventually resolve the singularities (albeit yielding a Deligne--Mumford stack).

The first issue we encounter is that the invariant and associated centre are defined for varieties but not yet for stacks. To be able to define them for Deligne--Mumford pairs (that is, a Deligne--Mumford stack $\call X$ embedded in smooth Deligne--Mumford stack $\call Y$), we first require two properties:
\begin{enumerate}[label=(\Alph*)]
	\item For any smooth morphism $\pi: Y'\ra Y$, we have $\pi^*\inv_X= \inv_{\pi^{-1}(X)}$;\label{eq: A}
	\item For any smooth \emph{and surjective} $\pi: Y'\ra Y$, if $\cF_\bullet$ is the centre associated to $X\subset Y$, then $\cF_\bullet\cdot\cO_{Y'}$ is the associated centre to $\pi^{-1}(X)\subset Y'$.\label{eq: B}
\end{enumerate}
The reason we require surjectivity in condition \ref{eq: B} is straightforward: if $Z$ is the locus of maximal invariant in $X$, consider the open immersion $Y\setminus Z \hookrightarrow Y$. There is no reason for the associated centre of $X\setminus Z \subset Y\setminus Z$ to be related to that of $X \subset Y$.

Using \ref{eq: A}, we define the invariant for a Deligne--Mumford  pair $\call X\subset \call Y$  by considering any smooth morphism $\pi:Y\ra \call Y$ from a scheme mapping $q\mapsto p$ and letting $\inv_\call X(p)=\inv_{\pi^{-1}(\call X)}(q)$.

Using \ref{eq: B}, we can define the associated centre to $\call X\subset \call Y$ to be the ideal filtration on $\cO_\call Y$ that pulls back to the associated centre along any smooth and surjective map from a variety; it is straightforward to see that this descent data is effective.

The final property we require of the invariant and associated centre ensures that we can degenerate to the weighted normal cone in the proof of resolution. Let $B=B_{\cF_\bullet}$ be the degeneration to the normal bundle of $Z$ with respect to $\cF_\bullet$, let $\mbf V\subset B$ be the vertex, and let $BX\subset B$ be the corresponding degeneration to the normal cone to $X$.
\begin{enumerate}
	\myitem{(C)} $\maxinv_X=\maxinv_{BX}$, and the locus with $\inv_{BX}(p)=\maxinv_{BX}$ in $BX$ is $\mbf V$.\label{eq: C}
\end{enumerate}
\begin{theorem}\label{thm: abc res}
	If the invariant and the associated centre satisfy properties \ref{eq: A}, \ref{eq: B} and \ref{eq: C}, then repeated blow-ups at the associated centre will resolve all singularities.
\end{theorem}

Compare \autoref{thm: abc res} with \cite[3.3.33]{wlodarczyk2023functorialresolutiontorusactions}, \cite[4.4.18]{Włodarczyk2023} and \cite[Lemma 4.1.1]{abramovich2025resolutionsingularitiesdynamicalmathematician}.
\begin{proof}
	Let $(X_0\subset Y_0):=(X\subset Y)$. For $j>0$, if $X_{j-1}$ is singular, define $Y_j$ to be the blow-up of $Y_{j-1}$ at the associated centre to $X_{j-1}\subset Y_{j-1}$, and let $X_j$ be the proper transform of $X_{j-1}$. If we are able to show that $\maxinv_{X_{j}}<\maxinv_{X_{j-1}}$ for all $j$, then, since $\Gamma^{(m)}$ is well-ordered and its minimal value indicates smoothness, there is some $m>0$ such that $X_m$ is smooth. Furthermore, note that after an étale base-change to an atlas, it is enough to establish this for schemes; therefore, it suffices to show that $\maxinv_{X_{1}}<\maxinv_X$.
	
	Since $\mbf V\subset BX$ is precisely the locus of maximal invariant and $\maxinv_{BX}=\maxinv_X$, we conclude that $\maxinv_{BX\setminus \mbf V}<\maxinv_X$. However, since the map
	\[
	BX\setminus \mbf V \ra \brk{
		(BX\setminus \mbf V)/\bb G_m
	}
	=X_1
	\]
	is smooth, it preserves the invariant by \ref{eq: A}. Therefore, $\maxinv_{X_1}<\maxinv_{X}$.
\end{proof}

It follows that the resulting resolution functor on resolution pairs is functorial for all smooth, surjective morphisms. In fact, one may show that it is functorial for all smooth morphisms. 

It is natural to ask whether we can treat the case of varieties that are not embeddable in any smooth variety (see \cite[Example 4.9]{roth2004affinestratificationnumbermoduli} for instance). As it turns out, after restricting the associated centre to $X$, the filtration $\cF_\bullet|_X$ is independent of the embedding (see \autoref{cor: embedding invariance} below). Therefore, by locally embedding $X$ in smooth varieties, the proper transform $\text{Spec}_X(\cR|_X)$ is defined globally, and so we can also resolve $X$. This yields a ``non-embedded'' resolution functor which is functorial for smooth morphisms; see \cite[Theorem 8.1.1]{abramovich2024functorial}.

Finally, one can apply a destackification algorithm from \cite{Bergh_2017, bergh2019functorialdestackificationweakfactorization} to obtain a resolution by a \emph{variety}, if desired, and the resulting procedure remains functorial (see \cite[Theorem 8.1.3]{abramovich2024functorial}).

The following properties of the invariant are easier to derive than \ref{eq: C}.
\begin{enumerate}
	\myitem{(D)} For $W\subset Y$ a smooth subvariety such that $X_W=X\cap W\subset W$ is of positive and pure codimension, and $p\in W$, we have $\inv_{X\subset Y}(p)\leq\inv_{X_W\subset W}(p)$;\label{eq: D}
	\myitem{(E)} Degenerating to the weighted normal cone preserves the invariant. That is, considering the embedding $BX_0\subset B_0$, for any $p\in Z$, where $Z$ is embedded via $Z=\mbf V_0$ in $B_0$, we have $\inv_{BX_0}(p)=\inv_X(p)$.\label{eq: E}
\end{enumerate}
\begin{proposition}
	Properties \ref{eq: A}, \ref{eq: D} and \ref{eq: E} imply property \ref{eq: C}.
\end{proposition}\label{prop: reduction conditions}
\begin{proof}
	First note that over $s\neq 0$, $(BX_{s\neq 0}\subset B_{s\neq 0})\simeq (X\times \bb G_m\subset Y\times\bb G_m)$, so that property \ref{eq: A} implies that $\maxinv_{BX_{s\neq 0}\subset B_{s\neq 0}}=\maxinv_X$ and that this invariant is attained exactly on $\mbf V_{s\neq 0}$. By \ref{eq: A} again, $\inv_{BX_{s\neq 0}\subset B_{s\neq 0}}(p)=\inv_{BX\subset B}(p)$ for any $p\in BX_{s\neq 0}$. Let us now move to the normal cone $BX_0$. Since $\maxinv_X=\inv_{BX_0}(p)$ for all $p\in Z$, and any closed $\bb G_m$-invariant subvariety of $BX_0$ intersects $Z$, the fact that the locus of maximal invariant is closed and $\bb G_m$-invariant implies that $\maxinv_{BX_0}=\maxinv_X$.
	
	By \ref{eq: D}, for any $p\in BX_0$, we have the inequality $\inv_{BX}(p)\leq \inv_{BX_0}(p)$. Thus, $\maxinv_X=\maxinv_{BX}$. Now, the locus of maximal invariant is smooth, closed, and $\bb G_m$-invariant. The only such subvariety of $B$ which restricts to $\mbf V_{s\neq 0}$ over $s\neq 0$ is $\mbf V$. Therefore, $\mbf V$ is the locus of maximal invariant with respect to $BX\subset B$.
\end{proof}
In what follows, we construct an invariant and an associated centre that satisfy properties \ref{eq: A}, \ref{eq: B}, \ref{eq: D} and \ref{eq: E}\,---\,properties that are sufficient to ensure the success of the resolution algorithm.

%
% !TEX root = main.tex
\section{The associated centre}\label{sec:associated centre}
Throughout this section, let $\cI\subset\cO_Y$ be an ideal such that $X=V(\cI)\subset Y$ has pure and positive codimension. We develop \hyperref[subsec: method 1]{Method 1} for constructing the centre associated to $X\subset Y$. We show that it has an interpretation using the Newton graph of $\cI$. We also show that the invariant and associated centre satisfy conditions \ref{eq: A} and \ref{eq: D}, and that the construction is independent of the embedding.
\begin{definition}\label{I-admissible}
	We say that a marked centre $\call J=(\cF_\bullet,d)$ is \tbf{$\call I$-admissible} if $v_{\call J}(\call I)\geq 1$, that is, if $\call I\subset \call F_d$. A pre-invariant is $\cI$-\tbf{admissible} if it is the invariant of an admissible marked centre.
\end{definition}
\begin{remark}
	Observe that if a smooth centre $\cF_\bullet$ has $v_{\cF_\bullet}(\cI)>0$, then the marked centre $(\cF_\bullet,v_{\cF_\bullet}(\cI))$ is $\cI$-admissible. Thus, a smooth centre can be endowed with a marking that makes it $\cI$-admissible if and only if $v_{\cF_\bullet}(\cI)>0$.
	
	Just as one only considers the proper transform of $X\subset Y$ under the blow-up of $Y$ at an ideal $\cF_1$ if $V(\cF_1)\subset X = V(\cI)$ (equivalently, $\cI\subset \cF_1$), one should only consider the proper transform of $X$ under a weighted blow-up of $Y$ at a smooth centre $\cF_\bullet$ when $v_{\cF_\bullet}(\cI)>0$.	 $\cI$-admissibility may be regarded as a normalisation of this condition, and the marking adequately measures the degree of tangency of the weighting to the singularities of $X$.
\end{remark}
\begin{remark}[Graphical interpretation of admissibility]\label{rem: graph int adm}
	The concept of $\cI$-admissibility has a very natural graphical interpretation, which we now describe. We work in the local ring at $p\in Y$, where we choose a system of parameters $\mx$. The Newton graph of $\cI$ with respect to $\mx$ is the subset of $\Z_{\geq 0}^n\subset \bb Q^n$ defined by plotting all monomials appearing in the expansion of some element in $\cI$:
	\[
	\mbf{Newt}(\cI,\mx)=\mbf{Newt}(\mx):=\brc{
		\alpha=(\alpha_1,\dots,\alpha_n): \partial_\mx^\alpha\cI=(1)
	}.
	\]
	Let $\cJ$ be an $\mx$-compatible marked centre with invariant $\mbf a=(a_1,\dots,a_k)$. Consider the linear functional $\lambda_{\mbf a}$ on $\bb Q^n$ defined by 
	\begin{equation*}
		\lambda_{\mbf a}(\mbf v):=\frac{v_1}{a_1}+\cdots+\frac{v_k}{a_k},
	\end{equation*}
	where $\mbf v=(v_1,\dots,v_n)$ denotes the standard coordinates on  $\bb Q^n$. 
	The equation $\lambda_{\mbf a}(\mbf v)\geq 1$ is a graphical counterpart of  the equation $v_\cJ(\cI)\geq 1$; indeed, $\cJ$ is $\cI$-admissible if and only if $\mbf{Newt}(\mx)$ is contained in the $\lambda_{\mbf a}(\mbf v)\geq 1$ half-space. Equivalently,  considering the hyperplane
	\begin{equation}
	H(\mbf a):=\{\mbf v: \lambda_{\mbf a}(\mbf v)=1\}\subset \bb Q^n,\label{eq: H(a)}
	\end{equation}
	$\cJ$ is $\cI$-admissible if and only if $\mbf{Newt}(\mx)$ lies in the half-space determined by $H(\mbf a)$ that does not contain the origin.
	We will write this condition as
	\[
	H(\mbf a)\leq \mbf{Newt}(\mx),
	\]
	and say that $H(\mbf a)$ \tbf{lies below} $\mbf{Newt}(\mx)$.
	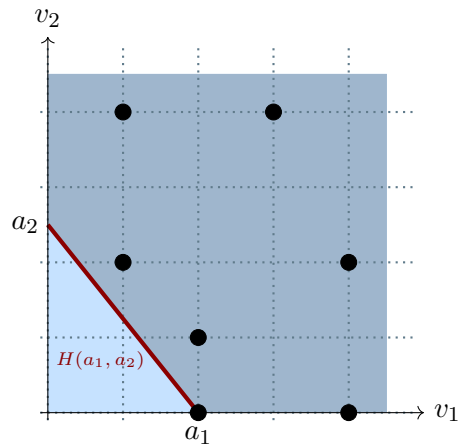
\begin{figure}[htpb]
		\begin{tikzpicture}
			\filldraw[SlateGray1] (0,0)--(0,4.5)--(4.5,4.5)--(4.5,0);
			\filldraw[SlateGray3] (2,0)--(0,2.5)--(0,4.5)--(4.5,4.5)--(4.5,0);
			\draw[thick,color=LightSkyBlue4, dotted] (-0.1,-0.1) grid (4.9,4.9); 
			\draw[->] (-0.1,0) -- (5.0,0) node[right] {$v_1$};
			\draw[ultra thick, -, color=Red4] (2,0)--(0,2.5);
			\draw[->] (0,-0.1) -- (0,5.0) node[above] {$v_2$};
			\node at (0.7,0.7) {\tiny$\color{Red4}{H(a_1,a_2)}$};
			\node at (2,-0.3) {$a_1$};
			\node at (-0.3,2.5) {$a_2$};
			\filldraw (4,0) circle (3pt);
			\filldraw (1,4) circle (3pt);
			\filldraw (2,1) circle (3pt);
			\filldraw (1,2) circle (3pt);
			\filldraw (3,4) circle (3pt);
			\filldraw (4,2) circle (3pt);
			\filldraw (2,0) circle (3pt);
		\end{tikzpicture}
		\caption{Illustration of $H(a_1,a_2)\leq\mbf{Newt}(\mx)$}\label{fig0}
	\end{figure}
	\autoref{fig0} illustrates this condition in the plane.
\end{remark}
\subsection{Uniqueness of a maximal marked centre}

Let $p\in X$ be a point. We will work in local rings at $p$, and so equality of ideals is understood there. Similarly, whenever we talk about admissibility, we mean \hyperref[I-admissible]{$\cI$-admissibility} at $p$. \hyperref[marked centre]{Marked centres}, \hyperref[parameters]{parameters}, and so on are also understood to be at $p$.

In what follows, we use the lexicographic order on \hyperref[def: pre-invariant]{pre-invariants}, with the slight modification that truncated sequences are considered greater. For instance, we have
$$(1,2,2,2)<(1,2,2)<(1,3,3)<(1)<(2,2)<(2,3).$$
This is equivalent to putting the usual lexicographic order on sequences of length $n$ where we set the last entries to be $\infty$. It is useful to retain this perspective.
Note that this is a total order. It induces a preorder on marked centres where $\cJ\leq\cJ'$ if $\text{inv}(\cJ)\leq\text{inv}(\cJ')$. The goal of this section is to prove the following theorem.
\begin{theorem}[{\citeplain[Theorem~5.3.1]{abramovich2024functorial}}]
	\label{thm: unique maximum} 
	The set of admissible marked centres at $p$ has a unique maximal element.\end{theorem}
	We can use \autoref{rem: graph int adm} to obtain a statement equivalent to \autoref{thm: unique maximum} in the language of convex geometry.
	\begin{theorem}[{\citeplain{Bierstone1997}}]\label{thm: unique maximum graph}
		The supremum
		\begin{equation}
			(a_1,\dots,a_k)=\sup_{(\mx,\mbf b)\;:\;H(\mbf b)\leq\mbf{Newt}(\mx)}\mbf b,\label{eq: graph int inv}
		\end{equation}
		where $\mx$ ranges over systems of parameters at $p$ and $\mbf b$ ranges over pre-invariants, has a maximiser. Moreover, any choice of maximising system of parameters $\mx$ gives the same admissible marked centre $\cJ=(x_1^{a_1},\dots,x_k^{a_k})$ (i.e. $\cJ$ is independent of the choice of maximiser).
	\end{theorem}
\begin{remark}\label{rem: graph int thm}

It is shown in \cite[Remark 1.11]{Bierstone1997} that the invariant satisfies the equality \eqref{eq: graph int inv} and that this supremum is achieved by a system of parameters. Although the problem was not framed in terms of marked centres, they already observed that any two maximising systems of parameters $\mx$ and $\my$ have to be related by a parameter transformation that preserves the weighted degree of the weighted leading terms. Thus, in light of \autoref{lem:equal-weightings}, we may extract the full content of \autoref{thm: unique maximum graph} from \cite{Bierstone1997}.
\end{remark}
Before giving a proof of \autoref{thm: unique maximum}, we will need a few definitions.

\begin{definition}
	We say that a pre-invariant of length $l\geq j\geq 0$ \tbf{dominates} a pre-invariant of length $j$ if the first $j$ entries of both sequences agree.	 We say that a marked centre $\cJ'$ \tbf{dominates} a marked centre $\cJ$ if $\inv(\cJ')$ dominates $\inv(\cJ)$ and both centres have a common compatible system of parameters.
\end{definition}

\begin{remark}
	More visually, $\cJ'$ dominates $\cJ$ if and only if there exists a system of parameters $\mx$ at $p$ and integers $0\leq j\leq l\leq n$ such that $\cJ=(x_1^{a_1},\dots,x_j^{a_j})$ and $\cJ'=(x_1^{a_1},\dots,x_j^{a_j},\dots,x_l^{a_l})$ where $(a_1,\dots,a_l)=\inv(\cJ')$.
\end{remark}
\begin{definition}\label{def: semi-associativity}
	A pre-invariant of length $j$ is \tbf{$j$-semi-admissible} if it is dominated by an admissible sequence.
	A marked centre $\cJ$ is \tbf{$j$-semi-associated} if the following two conditions hold:
	\begin{itemize}
		\item The maximal $j$-semi-admissible pre-invariant exists and is equal to $\inv(\cJ)$;
		\item If $\cJ'$ is admissible and $\inv(\cJ')$ dominates $\inv(\cJ)$, then $\cJ'$ dominates $\cJ$.\qedbarhere
	\end{itemize}
\end{definition}
\begin{remark}
	More visually, a marked centre $(x_1^{a_1},\dots,x_j^{a_j})$, where $(a_1,\dots,a_j)$ is the maximal $j$-semi-admissible pre-invariant, is $j$-semi-associated if and only if any admissible marked centre whose invariant starts by $(a_1,\dots,a_j)$ is of the form $(x_1^{a_1},\dots,x_l^{a_l})$ for some $l\geq j$.
\end{remark}
\begin{proposition}\label{k-associated implies associated}
	Let $k$ be an integer. If $\cJ$ is both $k$-semi-associated and admissible, then it is the unique maximal admissible marked centre.
\end{proposition}
\begin{proof}
	Since $k$-semi-admissible sequences bound admissible sequences from above (recall truncated sequences are greater), we know that the maximal admissible pre-invariant exists and equals $\inv(\cJ)$. Suppose $\cJ'$ is admissible and also has this maximal invariant. Then, since $\cJ$ is $k$-semi-associated and $\cJ'$ is admissible,  $\cJ'$ must dominate $\cJ$, so that $\cJ'=\cJ$.
\end{proof}
Another construction shall prove useful. Let $\cJ$ be a marked centre at $p$ and let $\mx$ be a compatible system of parameters. Let $(a_1,\dots,a_k)=\inv(\cJ)$. For $b\geq a_k$, we denote
\[
\cJ[\mx;b]:=(x_1^{a_1},\dots,x_k^{a_k},x_{k+1}^b,\dots,x_n^b).
\]
\begin{proposition}
	Let $\my$ be another $\cJ$-compatible system of parameters at $p$. Then, $\cJ[\mx;b] =\cJ[\my;b]$.	
\end{proposition}
\begin{proof}
	Let $i\leq k$. Note that $\langle dy_l:a_l=a_i\rangle=\langle dx_l:a_l=a_i\rangle$ by looking at the \hyperref[conormal filtration]{conormal filtration} on the underlying weighting of $\cJ$. Thus, in its $\mx$-expansion, $y_i$ has an $x_l$ linear term where $a_l=a_i$, meaning $v_{\cJ[\mx;b]}(y_i)\leq 1/a_i=v_{\cJ}(y_i)$. But since we have $v_{\cJ}\leq v_{\cJ[\mx;b]}$ (on functions, rather than on differential operators), we obtain $v_{\cJ[\mx;b]}(y_i)=1/a_i$. 	
	Let $i>k$. Since $dy_i\notin\langle dy_1,\dots,dy_k\rangle=\langle dx_1,\dots,dx_k\rangle$, $y_i$ must have a linear term in $x_{k+1},\dots,x_n$, so that $v_{\cJ[\mx;b]}(y_i)=1/b$. By \autoref{lem:equal-weightings}, we conclude that $\cJ[\mx;b]=\cJ[\my;b]$.
\end{proof}
\begin{definition}\label{def: b-comp}
	We define the \tbf{$b$-completion} of $\cJ$ at $p$ to be the marked centre 
	\[
	\cJ[b]=\cJ[p,b]:=\cJ[\mx;b]
	\]
	where $\mx$ is any $\cJ$-compatible system of parameters at $p$.	
\end{definition}

\subsection{Constructing $\cJ^{(j+1)}$}
Using \autoref{k-associated implies associated}, our strategy to prove \autoref{thm: unique maximum} is to show inductively the existence of a $j$-semi-associated marked centre $\cJ^{(j)}=(x_1^{a_1},\dots,x_j^{a_j})$ for each $j\leq k$, until we obtain an admissible marked centre $\cJ^{(k)}$. At each step, we will find the next entry $a_{j+1}$ in the maximal admissible pre-invariant and a suitable parameter $\bar x_{j+1}$ so that $\cJ^{(j+1)}=(x_1^{a_1},\dots,x_j^{a_j},\bar x_{j+1}^{a_{j+1}})$ is $(j+1)$-semi-associated. Note that while $a_{j+1}$ is uniquely determined, the choice of $\bar{x}_{j+1}$ is never unique. This provides freedom when constructing the semi-associated marked centres. The following notation will prove useful.

\begin{notation}
	Let $(a_1,\dots,a_j)$ be a pre-invariant, and let $\beta$ be a multi-index of length $l(\beta)\geq j$. We write
	\[
	\Delta(a_1,\dots,a_j;\beta)=\sum_{i=1}^j\frac{\beta_i}{a_i}.
	\]
	If moreover $l(\beta)\geq j+1$, we write
	\[
	\Xi(a_1,\dots,a_j;\beta)
	:=\frac{\sum_{i=j+1}^{l(\beta)}\beta_i}{1-\sum_{i=1}^j\frac{\beta_i}{a_i}}.
	\]
	When the pre-invariant is clear from the context, we will simply write $\Delta(\beta)$ and $\Xi(\beta)$.	 The choice of the symbol ``$\Delta$'' comes from the fact that the equation $\Delta(\beta)<1$ cuts out a weighted simplex in the space of multi-indices.
\end{notation}
\begin{remark}
	The ratio $\Xi(\beta)$ is similar to ones appearing in \cite[Proposition 5.4]{schober2014polyhedralapproachinvariantbierstone} and in \cite[Definition 18.2]{schoberhircharpol}, where  they are used for a polyhedral description of the Bierstone--Milman invariant \cite{Bierstone1997}. 
\end{remark}
We let $\cJ^{(j)}=(x_1^{a_1},\dots,x_j^{a_j})$ denote a $j$-semi-associated marked centre, which we assume not to be admissible.
\begin{lemma}\label{lem: J[b]}
	If $(a_1,\dots,a_j,b_{j+1},\dots, b_l)$ is admissible, then $\cJ^{(j)}[b_{j+1}]$ is admissible.\end{lemma}
\begin{proof}
	We have an admissible centre $\cJ'$ with this invariant, and it dominates $\cJ^{(j)}$ by $j$-semi-associativity, meaning that $\cJ'$ has a common compatible system of parameters with $\cJ^{(j)}[b_{j+1}]$. Therefore, the inequality $v_{\cJ^{(j)}[b_{j+1}]}\geq v_{\cJ'}$ holds on functions, and so  $\cJ^{(j)}[b_{j+1}]$ is admissible.
\end{proof}
\begin{corollary}\label{cor: a_{j+1} in terms of J[b]}
	If either of the two following maxima exists, then both exist and coincide:
	\begin{itemize}
		\item $\max\{b\in\bb Q: (a_1,\dots,a_j,b)\text{ is $(j+1)$-semi-admissible}\}$, 
		\item $\max\{ b\in\bb Q:\cJ^{(j)}[b]\text{ is admissible}\}$.
	\end{itemize} 
\end{corollary}

Let $\mx$ be a system of parameters at $p$ completing $x_1,\dots,x_j$, and consider the set 
\[
S=S_\mx=\brc{
	\beta :\Delta(\beta)<1\text{ and }\partial_{\mx}^\beta\cI=(1)
}
.
\]
\begin{lemma}\label{lem: existence minimum}
	The set $\Xi(S)=\{\Xi(\beta):\beta\in S\}$ has a minimal element.
\end{lemma}
\begin{proof}
	There are finitely many choices of $\beta_1,\dots,\beta_j$ such that $\Delta(\beta)<1$. It follows that
	\[\bigcup_{\Delta(\beta)<1}
	\frac{\Z_{\geq 0}}{1-\Delta(\beta)}\]
	is well-ordered, and so is $\Xi(S)$, \emph{a fortiori}. Since $\cJ^{(j)}$ is not admissible, \autoref{prop: duality} implies that $S$ and $\Xi(S)$ are non-empty.\end{proof}

\begin{proposition}\label{prop: a_{j+1}=c}
	The maxima of \autoref{cor: a_{j+1} in terms of J[b]} exist and are equal to $\min \Xi(S)$. 	
\end{proposition}

\begin{proof}
	Let $c=\min\Xi(S)$. We show that $c=\max\{b: \cJ^{(j)}[b]\text{ is admissible}\}$.
	Let $b>c$ and $b\geq a_j$ so that we may consider $\cJ^{(j)}[b]$.	Let $\beta\in S$ be a minimiser for $\Xi(S)$. We have 
	\[v_{\cJ^{(j)}[b]}(\partial_\mx^\beta)=-\sum_{i=1}^j\frac{\beta_i}{a_i}-\frac{\sum_{i=j+1}^n\beta_i}{b}>-\sum_{i=1}^j\frac{\beta_i}{a_i} -\frac{\sum_{i=j+1}^n\beta_i}{c}=-1\]
	by substituting $c=\Xi(\beta)$, which shows inadmissibility of $\cJ^{(j)}[b]$ by \autoref{prop: duality}. By \autoref{lem: J[b]}, we may conclude that $c\geq a_j$ since $\cJ^{(j)}[a_j]$ is admissible. Let $\gamma\in S$. We have 
	\[v_{\cJ^{(j)}[c]}(\partial_\mx^\gamma)=-\sum_{i=1}^j\frac{\beta_i}{a_i}-\frac{\sum_{i=j+1}^n\gamma_i}{c}\leq -\sum_{i=1}^j\frac{\gamma_i}{a_i}-\frac{\sum_{i=j+1}^n\gamma_i}{\Xi(\gamma)}=-1\]
	by minimality of $c$. By \autoref{prop: duality}, it follows that $\cJ^{(j)}[c]$ is admissible.\end{proof}
We have thus constructed the next entry in the maximal admissible pre-invariant (whose existence we are trying to establish), which we denote by \[a_{j+1}:=\min\Xi(S)=\max\{b\in\bb Q:(a_1,\dots,a_j,b) \text{ is $(j+1)$-semi-admissible}\}.\]
We now proceed to find the next parameter.
Let $\beta\in S$ be a minimiser of $\Xi(S)$ and $\bar\beta=\beta-e_l$ where $e_l$ is the vector with only a $1$ in the $l$th entry, $l>j$ and $\beta_l>0$. By assumption, we may choose $f\in\cI$ such that $\partial_\mx^\beta f$ does not vanish at $p$. Let $\bar x_{j+1}=\partial_\mx^{\bar \beta}f$. 

\begin{proposition}\label{prop: ind semi-ass}
	
	The marked centre $\cJ^{(j+1)}:=(x_1^{a_1},\dots,x_j^{a_j},\bar x_{j+1}^{a_{j+1}})$ is $(j+1)$-semi-associated.
\end{proposition}
\begin{proof}
	Let $\cJ'$ be an admissible marked centre whose invariant dominates $(a_1,\dots,a_{j+1})$. Since $\cJ^{(j)}$ is $j$-semi-associated, $\cJ'$ dominates $\cJ^{(j)}$. By \autoref{lem: fdmtl chain rule lmma}, we have that 
	\[
	v_{\cJ'}(\partial_\mx^{\bar\beta})\geq -\sum_{i=1}^j\frac{\beta_i}{a_i}-\frac{\pa{
			\sum_{i=j+1}^n\beta_i
		}
		-1}{a_{j+1}}.
	\]
	Substituting $a_{j+1}=\Xi(\beta)$ yields $v_{\cJ'}(\partial_\mx^{\bar\beta})\geq -1+1/a_{j+1}$. Thus, 
	\[
	v_{\cJ'}(\bar x_{j+1})\geq v_{\cJ'}(f)-1+1/a_{j+1}\geq 1/a_{j+1}
	\]
	by admissibility of $\cJ'$. However, since $\partial_\mx^\beta f=\partial_{\mx_l}\bar x_{j+1}$ does not vanish at $p$, in its $\mx$-expansion, $\bar x_{j+1}$ has an $x_l$ linear term so that in fact $v_{\cJ'}(\bar x_{j+1})=1/a_{j+1}$. That $\bar x_{j+1}$ has an $x_l$ linear term further ensures that $(x_1,\dots,x_j,\bar x_{j+1})$ is a sequence of parameters. By \autoref{lem:equal-weightings}, we conclude that $\cJ'$ dominates $\cJ^{(j+1)}$.
\end{proof}

We now combine all the results to prove \autoref{thm: unique maximum}.
\begin{proof}[Proof of \autoref{thm: unique maximum}]
	By \autoref{k-associated implies associated}, we need only show that there exists an admissible $k$-semi-associated marked centre $\cJ^{(k)}=\cJ$. Clearly, $\cJ^{(0)}=()$ is $0$-semi-associated, and so by \autoref{prop: ind semi-ass}, there exists a $j$-semi-associated marked centre $\cJ^{(j)}$ for all $j$ until we get to an admissible $k$-semi-associated centre (note that this must happen for $k\leq n$ since $n$-semi-admissibility is equivalent to admissibility).
\end{proof}	
\subsection{Method 1} We may now define the associated centre at a point.
\begin{definition}[Associated centre]\label{def: ass centre at point}
	We say that the unique maximal marked centre $\cJ$ is the \tbf{associated marked centre} to $\cI$ at $p$, and that its underlying smooth centre $\cF_\bullet$ is the \tbf{associated centre} to $\cI$ at $p$. When $\cI$ is clear from the context, we may simply call $\cJ$ (\emph{resp.} $\cF_\bullet$) the associated marked (\emph{resp.} associated) centre at $p$. We call $\inv(\cJ)$ the \tbf{invariant} of $\cI$ at $p$, for which we write $\inv_{X\subset Y}(p)=\inv_X(p)$. If $\inv_X(p)$ has length $k$, 
	for any $j\leq k$, we write $a_{j}^{X\subset Y}(p)=a_{j}^X(p)$ for the $j$th entry in $\inv_X(p)$. If $\cI$ is clear from the context, we may simply write $\inv(p)$ and $a_j(p)$.
\end{definition}
\begin{method}\label{met1}
	The above construction of $\cJ$ may be summarised in the form of an algorithm to compute the associated marked centre at a point. This is \hyperref[subsec: method 1]{Method 1} as presented in the \hyperref[sec: intro]{Introduction}. The input will be a $j$-semi-associated marked centre, and the output will be a $(j+1)$-semi-associated marked centre. The procedure starts with $()$ and terminates when we obtain a $k$-semi-associated marked centre $\cJ^{(k)}=\cJ$ that is admissible.
	
	Let $\cJ^{(j)}=(x_1^{a_1},\dots,x_j^{a_j})$ be $j$-semi-associated and let $\mx$ be a system of parameters completing $x_1,\dots, x_j$. Choose $\beta$ to be a multi-index of length $n$ that is a minimiser of
	\begin{equation}
		a_{j+1}=\min\brc{
			\Xi(\beta): \partial_\mx^\beta\cI=(1)\text{ and }\Delta(\beta)<1
		}
		.
		\label{eq: a_{j+1} is descr}
	\end{equation}
	Choose $l>j$ such that $\beta_l>0$, let $\bar\beta=\beta-e_l$. Choose $f\in\cI$ such that $\partial_\mx^{\beta}f$ does not vanish at $p$ and let $\bar x_{j+1}=\partial_\mx^{\bar \beta}f$. Then, $\cJ^{(j+1)}=(x_1^{a_1},\dots,x_j^{a_j},\bar x_{j+1}^{a_{j+1}})$ is $(j+1)$-semi-associated. 
\end{method}
\begin{remark}
	The invariant detects how bad a singularity at a point is. In particular, $X$ is smooth at $p$ if and only if $\call I=(x_1,\dots,x_m)$ for some parameters $x_1,\dots,x_m$ at $p$, if and only if the invariant at $p$ is the sequence
	\[
	(\underbrace{1,\dots,1}_{\text{length } m}).
	\]
	Any larger invariant indicates a singularity. Conversely, since $X$ has codimension $m$, this is the smallest attainable invariant.
\end{remark}
The procedure implies that the invariant takes its values in some precise set of pre-invariants. Specifically, if $(b_1,\dots,b_j)$ are the first entries of an invariant, the next entry $b_{j+1}$ must be in the set
\[
\Gamma(b_1,\dots,b_j):=[b_j,\infty)\cap \bigcup_{\Delta(b_1,\dots,b_j;\beta)<1}\frac{\Z}{1-\Delta(b_1,\dots,b_j;\beta)},
\]
which is well-ordered. Define now inductively $\Gamma_1=\Z_{>0}$ and 
\[
\Gamma_{j+1}=\brc{
	\pa{
		b_1,\dots,b_{j+1}
	}
	:\pa{
		b_1,\dots,b_{j}
	}
	\in\Gamma_{j},b_{j+1}\in \Gamma\pa{
		b_1,\dots,b_j
	}
}
.
\]
\begin{proposition}\label{prop: well-orderedness}
	$\Gamma=\sqcup_{j\leq n}\Gamma_j$ is well-ordered.
\end{proposition}
\begin{proof}
	Suppose that $\Gamma_j$ is well-ordered and let $A\subset \Gamma_{j+1}$ be a set. Since $\Gamma_j$ is well-ordered, the truncation of $A$ has a minimum $(b_1,\dots,b_j)\in\Gamma_j$, and so to check if $A$ has a minimum, we may assume all its elements dominate $(b_1,\dots,b_j)$. Since $\Gamma(b_1,\dots,b_j)$ is well-ordered, $A$ has a minimum and so $\Gamma_{j+1}$ is well-ordered.
\end{proof}
Let $m$ be the codimension of $X$ in $Y$, and let $\Gamma^{(m)}$ be the subset of $\Gamma$ of those pre-invariants which have value at least
\[
(\underbrace{1,\dots,1}_{\text{length }m}),
\]
or equivalently of those pre-invariants which have length at least $m$.
The invariant function and the set $\Gamma^{(m)}$ satisfy the conditions \ref{eq: (i)} and \ref{eq: (ii)} imposed in the beginning of \autoref{subsec: what is needed for res}, and we shall show in \autoref{sec: global associated centre} that they also satisfy condition \ref{eq: (iii)}.
\begin{remark}
	The set $\Gamma$ has been independently discovered by \cite{temkin2025dreamresolutionprincipalizationii}, where it is denoted $\underline{\text{ORD}}$ and given an equivalent characterisation, which we now detail.\footnote{We thank Dan Abramovich for identifying the connection between $\Gamma$ and $\underline{\text{ORD}}$.} Given $(b_1,\dots,b_k)\in \Gamma$, for all $j\leq k$, we have by definition $b_{j}\in \Gamma(b_1,\dots,b_{j-1})$, which is equivalent to requiring that $b_{j}\geq b_{j-1}$ and that there exist $\alpha_1,\dots,\alpha_j$ such that $\sum_{i=1}^{j}\frac{\alpha_i}{b_i}=1$. Thus, the set $\Gamma$ can be characterised precisely as those pre-invariants $(b_1,\dots,b_k)$ which satisfy an integrality condition, namely that for all $j\leq k$ there is a multi-index $\alpha^{(j)}$ such that 
	\[
	\sum_{i=1}^j\frac{\alpha^{(j)}_i}{b_i}=1.
	\]
	Moreover, it is shown in \emph{op.\ cit.} that any $(b_1,\dots,b_k)\in \Gamma$ is the invariant at $p$ of some ideal: choosing any sequence of parameters $x_1,\dots,x_k$ at $p$, the marked centre $(x_1^{b_1},\dots,x_k^{b_k})=\cJ=(\cF_\bullet,d)$ is associated to the ideal $\cF_d=\pa{f\in\cO_Y:v_{\cJ}(f)\geq 1}$, which, in the terminology of \emph{op.\ cit.}, is called the \emph{rounding} of $\cJ$; see \cite[Theorem 2.1.9]{temkin2025dreamresolutionprincipalizationii}. This refines previous estimates of the set of possible invariants; compare \cite[Remark 5.1.2]{abramovich2024functorial} and \cite[Remark 2.1.10.ii]{temkin2025dreamresolutionprincipalizationii}. The study of $\underline{\text{ORD}}$ was a cornerstone for the extension of these resolution techniques to the case of quasi-excellent schemes of characteristic zero.
\end{remark}
\subsection{The Newton graph picture}\label{sec: Newton graph and numerical inequalities}

Using \autoref{rem: graph int adm} and \autoref{rem: graph int thm}, we give a graphical interpretation of the construction of the associated marked centre based on \autoref{met1}. This is in line with \cite{Schober2013, schober2014polyhedralapproachinvariantbierstone}, where the polyhedral interpretation of the Bierstone--Milman invariant was first proved.

Consider $\bb Q^n$ with standard coordinates $\mbf v=(v_1,\dots,v_n)$. 

At the $(j+1)$st step, we are given the following data:
\begin{itemize}
	\item A system of parameters $(x_1,\dots,x_j,y_{j+1},\dots,y_n)$, where $x_1,\dots,x_j$ are the parameters to be compatible with the associated centre and the rest of the parameters need to be refined;
	\item A hyperplane
	\begin{equation*}
		H_j:=H(a_1,\dots,a_j,\dots,a_j)\subset \bb Q^n
	\end{equation*}
	where we use the notation from \eqref{eq: H(a)} and $(a_1,\dots,a_j,\dots,a_j)$ has length $n$. $H_j$ should be thought of as the hyperplane corresponding to $\cJ^{(j)}[a_j]$.
\end{itemize}
 Note importantly that the $v_i$-intercept of $H_j$, that is the value at which $H_j$ meets the $v_i$-axis, is $a_i$ for $i< j$ and $a_j$ otherwise .
 
By construction, the following condition is satisfied.
\begin{equation}\label{eq: H_j lies under}
	H_j\leq\mbf{Newt}:= \mbf{Newt}(x_1,\dots,x_j,y_{j+1},\dots,y_n).
\end{equation}
We now want to modify our hyperplane $H_j$ by simultaneously increasing the $v_{j+1}$-,\dots, $v_n$-intercepts by the same amount, and we want to do so maximally while preserving condition \eqref{eq: H_j lies under}. That is, we want to find the minimal increase of these intercepts that touches $\mbf{Newt}$. Playing with the equations of these hyperplanes readily gives that we need to increase the intercepts to \begin{equation*}
	a_{j+1}=\min\brc{
		\Xi(\beta):\beta\in\mbf{Newt}\text{ and } \Delta(\beta)<1
	},
\end{equation*}
and we define 
\[
H_{j+1}:=H(a_1,\dots,a_{j+1},\dots,a_{j+1}),
\]
where $(a_1,\dots,a_{j+1},\dots,a_{j+1})$ has length $n$. To play this optimising game in the next steps, we need to find the next compatible parameter $x_{j+1}$. By construction, there exists $\beta\in \mbf{Newt}\cap H_{j+1}$ such that $\Delta(\beta)<1$. We choose $f\in\cI$ such that $\beta$ is a monomial of $f$ with respect to our system of parameters. We let $\bar\beta=\beta-e_i$ where $\beta_i>0$ and define
\[
x_{j+1}:=\partial_{(x_1,\dots,x_j,y_{j+1},\dots,y_n)}^{\bar\beta} f.
\]
We verify whether  $\cJ^{(j+1)}=(x_1^{a_1},\dots,x_{j+1}^{a_{j+1}})$ is admissible, i.e. whether\footnote{Here, we may need to reorder $y_{j+2},\dots,y_n$ to ensure that $(x_1,\dots,x_{j+1},y_{j+2},\dots,y_n)$ is a system of parameters.}
\[
H(a_1,\dots,a_{j+1})\leq \mbf{Newt}(x_1,\dots,x_{j+1},y_{j+2},\dots,y_n),
\]
and if it is not the case, we proceed to the next step. 

For plane curves, our graphical interpretation of the construction of the associated centre was found by Lapointe in a private draft from 2022 and recently appeared in \cite{abramovich2024resolvingplanecurvesusing,abramovich2025torusactionsweightedblowups}. We illustrate this process with an example of such a curve, as we can plot things in the plane.
\begin{example}\label{ex: newt curve}
	Let $f=x^4+xy^4+y^6$, and consider the corresponding curve $C\subset\bb A^2$. We compute the associated centre at the origin. Running \autoref{met1}, we see $a_1=4$ and that $x$ may be taken as the first parameter. This gives
	\[
	a_2=\min\brc{\frac{4}{1-\frac{1}{4}},\; 6}=16/3,
	\]
	and we can take $y$ as second parameter.
	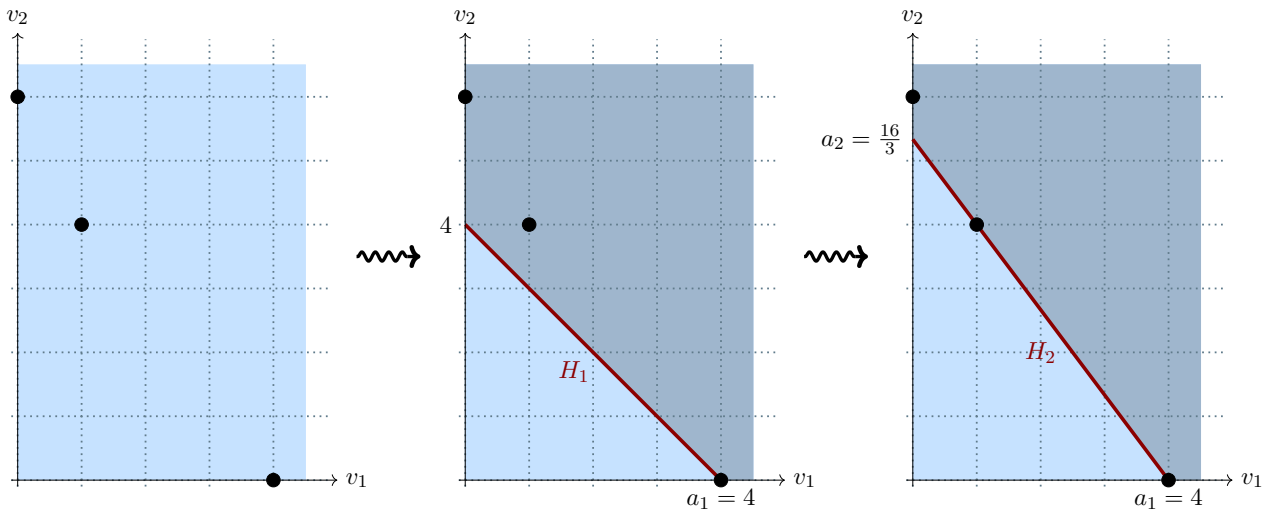
\begin{figure}[htbp]
		\centering
		\scalebox{0.85}{ 
			\begin{tikzpicture}
				\begin{scope}[xshift=0cm]
					\filldraw[SlateGray1] (0,0)--(0,6.5)--(4.5,6.5)--(4.5,0);
					\draw[thick,color=LightSkyBlue4, dotted] (-0.1,-0.1) grid (4.9,6.9); 
					\draw[->] (-0.1,0) -- (5.0,0) node[right] {$v_1$};
					\draw[->] (0,-0.1) -- (0,7.0) node[above] {$v_2$};
					\filldraw (4,0) circle (3pt);
					\filldraw (1,4) circle (3pt);
					\filldraw (0,6) circle (3pt);
				\end{scope}
				
				\node at (5.8, 3.5) {\thicksquigarrow};
				
				\begin{scope}[xshift=7cm]
					\filldraw[SlateGray1] (0,0)--(0,6.5)--(4.5,6.5)--(4.5,0);
					\filldraw[SlateGray3] (4,0)--(0,4)--(0,6.5)--(4.5,6.5)--(4.5,0);
					\draw[thick,color=LightSkyBlue4, dotted] (-0.1,-0.1) grid (4.9,6.9); 
					\draw[->] (-0.1,0) -- (5.0,0) node[right] {$v_1$};
					\draw[ultra thick, -, color=Red4] (4,0)--(0,4);
					\draw[->] (0,-0.1) -- (0,7.0) node[above] {$v_2$};
					\node at (1.7,1.7) {$\color{Red4} H_1$};
					\node at (-0.3,4) {$4$};
					\node at (4,-0.3) {$a_1=4$};
					\filldraw (4,0) circle (3pt);
					\filldraw (1,4) circle (3pt);
					\filldraw (0,6) circle (3pt);
				\end{scope}
				
				\node at (12.8, 3.5) {\thicksquigarrow};
				
				\begin{scope}[xshift=14cm]
					\filldraw[SlateGray1] (0,0)--(0,6.5)--(4.5,6.5)--(4.5,0);
					\filldraw[SlateGray3] (4,0)--(0,16/3)--(0,6.5)--(4.5,6.5)--(4.5,0);
					\draw[thick,color=LightSkyBlue4, dotted] (-0.1,-0.1) grid (4.9,6.9); 
					\draw[->] (-0.1,0) -- (5.0,0) node[right] {$v_1$};
					\draw[->] (0,-0.1) -- (0,7.0) node[above] {$v_2$};
					\draw[ultra thick, -, color=Red4] (4,0)--(0,16/3);
					\node at (2,2) {$\color{Red4} H_2$};
					\node at (-0.8,16/3) {$a_2=\frac{16}{3}$};
					\filldraw (4,0) circle (3pt);
					\filldraw (1,4) circle (3pt);
					\filldraw (0,6) circle (3pt);
					\node at (4,-0.3) {$a_1=4$};
				\end{scope}
				
			\end{tikzpicture}
		}
		
		\caption{Graphical construction of the associated centre for \autoref{ex: newt curve}}\label{fig1}
	\end{figure}
	\autoref{fig1} illustrates the construction in the plane.
\end{example}

\subsection{Smooth invariance, embedding invariance, restriction}
\begin{proposition}[Local functoriality;\;{\citeplain[Theorem 5.1.1]{abramovich2024functorial}}]\label{prop: local functoriality} If $\pi: Y'\ra Y$ 
	is a smooth morphism mapping $ p'\mapsto p$ and $\cF_\bullet$ and $\cF'_\bullet$ are the associated centres to $\cI$ at $p$ and to $\cI\cdot\cO_{Y'}$ at $p'$ respectively, then, $\cF'_\bullet=\cF_\bullet\cdot\cO_{Y'}$.
	
	In particular, the invariant satisfies condition \ref{eq: A}, that is $\pi^*\inv_X=\inv_{\pi^{-1}(X)}$.
\end{proposition}
\begin{proof}
	It suffices to observe that $(\partial_{\pi^*\mx}^\beta)(\pi^*f)=\pi^*(\partial_\mx^\beta f)$ to be convinced that every step in the procedure is preserved by smooth pullback. \end{proof}
Once we find the correct definition of the \emph{global} associated centre, \autoref{prop: local functoriality} will imply property \ref{eq: B} with ease.		

If we consider instead a closed immersion of a smooth variety $W\hookrightarrow Y$, the associated centre cannot be pulled back in the same manner as the morphism might ``destroy parameters''. Nevertheless, as with the order of vanishing, our invariant increases upon restriction to a smooth subvariety.
\begin{proposition}\label{prop: invariant increases on closed subvarieties}
	The invariant satisfies \ref{eq: D}. That is, if $W\subset Y$ is a smooth closed subvariety such that $X_W:=X\cap W\subset W$ is of positive and pure codimension, then, for all $p\in W$, 
	\[
	\inv_{X\subset Y}(p)\leq \inv_{X_W\subset W}(p).
	\]
\end{proposition}

\begin{proof}
	We prove the statement by induction. Suppose that for $j\geq 0$ there is a \hyperref[def: semi-associativity]{$j$-semi-associated} marked centre to $\cI|_W$ at $p$ which is the restriction $\cJ^{(j)}|_W$ of a $j$-semi-associated marked centre $\cJ^{(j)}$ to $\cI$ at $p$. Explicitly, we may write $\cJ^{(j)}=(x_1^{a_1},\dots,x_j^{a_j})$ for parameters $x_1,\dots,x_j$ which restrict to parameters on $W$. We may also complete $x_1,\dots,x_j$ into a system of parameters $\mx$ at $p$ such that $W=V(x_{n-w+1},\cdots,x_n)$ where $w=\codim_Y W$. We then have the inequality
	\begin{align*}
		a_{j+1}^{X}(p)&=\min_{\beta}\brc{
			\Xi(\beta):\Delta(\beta)<1\text{ and } \partial^\beta_\mx\cI=(1)
		}\\
		&\leq \min_{\beta}\brc{
			\Xi(\beta):\beta=(\beta_1,\dots,\beta_{n-w},0\dots,0),\;\Delta(\beta)<1 \text{ and }\partial_\mx^\beta\cI=(1)
		}
		=a_{j+1}^{X_W}(p).
	\end{align*}
	In particular, they are equal only if there is a common minimiser $\beta=(\beta_1,\dots,\beta_{n-w},0,\dots, 0)$ since the second minimum ranges over a subset of the range of the first minimum. In this case, for any $f\in\cI$, we have $(\partial_\mx^\beta f)|_W=\partial_{(x_1,\dots,x_{n-w})}^{(\beta_1,\dots,\beta_{n-w})}(f|_W)$, so that we can take the same next parameter for both associated marked centres.
	
	We have thus proved that if our inductive hypothesis holds for $j\geq 0$, then it either holds for $j+1$ or $\inv_{X}(p)<\inv_{X_W}(p)$. Since it holds for $j=0$, we can conclude that $\inv_{X}(p)\leq\inv_{X_W}(p)$.\end{proof}
As stated so far, the associated centre depends on a choice of embedding $X\hookrightarrow Y$. We now show that the filtration we get by restricting to $X$ is actually independent of the embedding. 
We first need a proposition about étale local equivalence of closed immersions of the same codimension.

\begin{proposition}
	Let $\iota_i:X\hookrightarrow Y_i$ for $i=1,2$ be two closed immersions with the same codimension at $p$. Then, there exists a neighbourhood $X_0\subset X$ of $p$ and a closed immersion $\iota:X_0\hookrightarrow W$ with étale morphisms $f_i:W\ra Y_i$ such that $\iota_i=f_i\circ \iota$. 
\end{proposition}
\begin{proof}
	Let $m=\dim_pY_i$ and consider the diagonal embedding $X\ra Y_1\times Y_2$. We have that
	\[
	\dim\ker\pa{
		T_p^*\pa{
			Y_1\times Y_2
		}
		\ra T_p^*X
	}
	-\dim \ker\pa{
		T_p^*Y_i\ra T_p^*X
	}
	=m,
	\]
	so we may choose linearly independent $dt_1,\dots, dt_m\in \ker(T_p^*(Y_1\times Y_2)\ra T_p^*X)=:K$ that are linearly independent of $\ker(T_p^*Y_i\ra T_p^*X)$ for $i=1,2$. Let $\cA$ be the ideal at $p$ corresponding to the closed immersion $X\to Y_1\times Y_2$. We have the equality
	\[
	K=\frac{\cA+\mf m_p^2}{\mf m_p^2}=\frac{\cA}{\cA\cap\mf m_p^2},
	\]
	so that we may lift $dt_1,\dots,dt_m$ to elements $t_1,\dots,t_m\in\cA$, and resulting zero scheme 
	\[
	W:=V(t_1,\dots,t_m)\supset X
	\]
	is smooth near $p$. By our assumptions, the map 
	\[
	T_p^*Y_i\ra T_p^*(Y_1\times Y_2)/\langle dt_1,\dots,dt_m\rangle=T^*_pW
	\]
	is injective, meaning $W\ra Y_i$ is smooth near $p$ and of relative dimension $0$, i.e. it is étale near $p$.
\end{proof}

\begin{corollary}\label{cor: same dim same centre}
	Let $\cF_\bullet^{Y_i}$ be the associated centre to $X\hookrightarrow Y_i$ at $p$. Then, $\cF_\bullet^{Y_1}|_X=\cF_\bullet^{Y_2}|_X$.
\end{corollary}
\begin{proof}
	The weightings are local, so we may assume without loss of generality that $X_0=X$. From \autoref{prop: local functoriality}, it follows that $\cF_\bullet^{W}=f_i^*\cF_\bullet^{Y_i}$ is the associated centre to $\iota: X\hookrightarrow W$ at $p$. Since $\iota_i=f_i\circ \iota$, we have that $\cF_\bullet^{Y_i}|_X=\iota_i^*\cF_\bullet^{Y_i}=\iota^*f_i^*\cF_\bullet^{Y_i}=\iota^*\cF_\bullet^{W}=\cF^{W}_{\bullet}|_X$.
\end{proof}

The following lemma is known as the \emph{re-embedding principle} in the literature.

\begin{lemma}[{\citeplain[Proposition 2.9.3]{abramovich2017principalizationidealstoroidalorbifolds}}]
	\label{lem: reembed}
	Let $\cF_\bullet$ be the associated centre to $X\hookrightarrow Y$ at $p$ and $\cG_\bullet$ be the associated centre to $X\simeq X\times\{0\}\hookrightarrow Y\times\Spec\Bbbk[x]$ at $p$. Then, $\cG_\bullet|_X=\cF_\bullet|_X$.
\end{lemma}
\begin{proof}
	Running the algorithm, if $(\cF_\bullet,d)=(x_1^{a_1},\dots,x_k^{a_k})$ is the associated marked centre to $X\hookrightarrow Y$ at $p$, then $(x,x_1^{a_1},\dots,x_k^{a_k})$ is the associated centre to $X\simeq X\times \{0\}\hookrightarrow Y\times\Spec\Bbbk[x]$ at $p$. If $w_1,\dots,w_k$ are the weights of $\cF_\bullet$, $(d,w_1,\dots,w_k)$ will be the weight sequence of the smooth centre underlying $(x,x_1^{a_1},\dots,x_k^{a_k})$. Since $x$ vanishes on $X$, restricting both smooth centres yields the same filtration on $\cO_X$.   \end{proof}
Combining \autoref{cor: same dim same centre} and \autoref{lem: reembed}, we deduce:
\begin{corollary}[Embedding invariance]\label{cor: embedding invariance}
	Let $\cF_\bullet$ be the associated centre to $X\hookrightarrow Y$ at $p$. Then, $\cF_\bullet|_X$ is independent of the embedding.
\end{corollary}
As explained in \autoref{subsec: what is needed for res}, this ensures the existence of a ``non-embedded'' resolution functor, once the existence of an embedded resolution functor has been proved. It also ensures that we have a canonical resolution for those varieties that do not embed in smooth varieties.

% 
% !TEX root = main.tex
\section{The global associated centre}\label{sec: global associated centre}
In this section, we first show that we can globalise the definition of the associated centre. The underlying subvariety of the resulting global associated centre is the locus where the invariant attains its maximum, and this locus is smooth and closed. We also develop \hyperref[met2]{Method 2} for computing the associated centre by using derivatives of ideals, which is best suited for a global algorithm.
Let $\cI\subset\cO_Y$ be an ideal such that $X=V(\cI)$ is of pure and positive codimension.

\subsection{Globalising}
For any $p\in X$, we write $\cJ(p)$ for the associated marked centre at $p$. $\cJ(p)$ being a marked centre at a point, $Z_{\cJ(p)}$ is the germ of a variety which is smooth at $p$. In particular, it extends to a unique irreducible subvariety of $Y$, which we also denote $Z_{\cJ(p)}$. Similarly, note that $\cJ(p)$ extends to a marked centre on some neighbourhood of $p$ (at least that where the parameters are defined), which shall also be denoted by $\cJ(p)$.

\begin{proposition}\label{prop: ass centre glblises}
	Near $p$, the following hold: 
	\begin{itemize}
		\item The maximal invariant is $(a_1(p),\dots,a_k(p))=(a_1,\dots,a_k)$;
		\item The locus where this invariant is attained is exactly $Z_{\cJ(p)}$;
		\item For all $q$ in this locus, $\cJ(q)=\cJ(p)$.
	\end{itemize}
\end{proposition}
\begin{proof}
	We show this by induction. Let $(x_1^{a_1},\dots,x_j^{a_j})$ be $j$-semi-associated at $p$. Suppose by induction that near $p$, if
	\[
	\inv(q)\geq (\underbrace{a_1,\dots,a_j,\dots,a_j}_{\text{length $n$}}),
	\]
	then $\cJ(q)$ dominates $(x_1^{a_1},\dots,x_j^{a_j})$. Then, for any such $q$, we have that
	\begin{equation}
		a_{j+1}(q)=\min_\beta \brc{
			\Xi(\beta): \Delta(\beta)<1\text{ and } \partial_\mx^\beta\cI=\cO_{Y,q}
		}
		\label{eq: a_{j+1}(q)},
	\end{equation}
	where $\mx$ is a system of parameters completing $x_1,\dots,x_j$. Let $\beta$ minimise \eqref{eq: a_{j+1}(q)} at $q=p$. Since the equality $\partial_\mx^\beta\cI=\cO_Y$ holds in a neighbourhood of $p$, we may conclude that near $p$, if 
	\[\inv(q)\geq (\underbrace{a_1,\dots,a_{j+1},\dots,a_{j+1}}_{\text{length $n$}}),\]
	then $a_{j+1}(q)=a_{j+1}$. Choose $f\in\cI$ such that $\partial_\mx^\beta f$ does not vanish at $p$ (and hence in a neighbourhood of $p$). Let $\overline{x}_{j+1}=\partial_\mx^{\overline\beta}f$, $\overline\beta=\beta-e_l$, where $l\geq j+1$ and $\beta_l>0$. By construction, $(x_1^{a_1},\dots,x_j^{a_j},\overline x_{j+1}^{a_{j+1}})$ is $(j+1)$-semi-associated at all $q$ near $p$ such that $\inv(q)\geq (a_1,\dots,a_{j+1},\dots,a_{j+1})$.
	
	Since the inductive hypothesis holds trivially for $j=0$, we conclude that near $p$, if
	\[
	\inv(q)\geq (\underbrace{a_1,\dots,a_k,\dots,a_k}_{\text{length } n})
	\]
	then $\cJ(q)$ dominates $\cJ(p)$, implying that $\inv(p)\geq \inv(q)$. Since $(a_1,\dots,a_k,\dots,a_k)\leq\inv(p)$, this means that $\inv(p)$ is maximal near $p$ and that the locus where this maximal invariant is attained is contained in $Z_{\cJ(p)}$. By \autoref{prop: val from pt to global}, $\cJ(p)$ is admissible at every point in $Z_{\cJ(p)}$ in a neighbourhood of $p$ (where $\cJ(p)$ is defined), which concludes the proof.
\end{proof}
\begin{proposition}\label{prop: locus of maximal invariant exists}
	The invariant function satisfies condition \ref{eq: (iii)}, that is, the maximal invariant $(a_1,\dots,a_k)=\maxinv_{X}:=\max_{p\,\in X}\inv(p)$ exists, and the locus where it is attained is a smooth and closed subvariety of $Y$ of pure codimension $k$.
\end{proposition}
\begin{proof}
	Cover $Y=\bigcup_{p\in X}U_p$ where each $U_p$ is chosen to be a neighbourhood of $p$ on which $\inv(p)$ is maximal and attained on $Z_{\cJ(p)}\cap U_p$, which is smooth. Because $Y$ is quasi-compact, we can refine this to a finite subcover, ensuring that the maximal invariant exists. Now, on each $U_p$ in our finite cover, the locus of maximal invariant is either a smooth and closed subvariety of codimension $k$, or empty, which implies the statement.\end{proof}

\begin{corollary}
	The function $\inv:|X|\ra \Gamma$ is upper semicontinuous (in the Zariski topology). In particular, the level sets are locally closed. Moreover, there are finitely many level sets.
\end{corollary}	
\begin{proof}
	This is shown inductively by discarding successive loci of maximal invariant and applying \autoref{prop: locus of maximal invariant exists}. Since the invariant is upper semicontinuous and valued in a well-ordered set, we conclude that it has finitely many level sets.
\end{proof}	

\begin{definition}[Associated centre]\label{def: global associated marked centre}
	The (\emph{global}) \tbf{associated marked centre} to $\cI$ is the unique admissible marked centre $\cJ=(\cF_\bullet,d)$ satisfying $\inv(\cJ)=\maxinv_X$ and whose underlying subvariety $Z_\cJ$ is the locus of maximal invariant. We call $\cF_\bullet$ the \tbf{associated centre} to $\cI$.
\end{definition}

\autoref{prop: local functoriality} directly implies property \ref{eq: B} for the associated centre.
\begin{corollary}\label{prop: global functoriality} The associated centre satisfies \ref{eq: B}. That is, if $\pi: Y'\ra Y$ is a smooth surjective morphism and $\cF_\bullet$ is the associated centre to $X\subset Y$, then $\cF_\bullet\cdot\cO_{Y'}$ is the associated centre to $\pi^{-1}(X)\subset Y'$.
	
\end{corollary}
It remains only to verify that the associated centre satisfies condition \ref{eq: E}. We first recall the content of \ref{eq: E}. Consider the normal cone $BX_0$ to $X$ at the associated centre $\cF_\bullet$. Let $Z\subset BX_0$ be embedded via its zero section, i.e. $Z=\mbf V_0$, where $\mbf V\subset B$ is the vertex inside the degeneration to the normal bundle. Condition \ref{eq: E} says that for any $p\in Z$, $\inv_X(p)=\inv_{BX_0\subset B_0}(p)$. Recall also that $BX_0\subset B_0$ corresponds to the ideal $\tilde \cI=(\tilde f:f\in\cI)\subset \cO_{B_0}$. If $\mx$ is an $\cF_\bullet$-compatible system of parameters at $p$, we obtain parameters $\tilde \mx$ at $p$ on $B_0$, and for
\begin{equation}
	f=\sum_{\beta}c_\beta^{(f)}\mx^\beta\in \cO_Y,\label{eq: f}
\end{equation}

we have
\begin{equation}
	\tilde f=\sum_{\substack{\beta\\v_{\cF_\bullet}(\mx^\beta)=v_{\cF_\bullet}(f)}}c_\beta^{(f)}\tilde\mx^\beta\in \cO_{B_0}.\label{eq: f twidle}
\end{equation}
\begin{proposition}\label{prop: invariant to normal cone}
	Let $\cJ=(\cF_\bullet,d)$ be the associated marked centre to $\cI$. Then, the associated marked centre to $\tilde \cI$ is the initialisation $\tilde\cJ=(\tilde\cF_\bullet,d)$ of $\cJ$, where $\tilde\cF_i:=(\tilde f : f\in\cF_i)$. In particular, the invariant satisfies condition \ref{eq: E}. 
\end{proposition}
\begin{proof} Recall that the locus of maximal invariant of $BX_0\subset B_0$ is contained in $Z$ since it is a closed and smooth subvariety of $B_0$ of positive codimension which is $\bb G_m$-invariant. Thus, it suffices to establish the statement at all points $p\in Z\subset B_0$.
	
	Let $\inv(\cJ)=(a_1,\dots,a_k)$. Consider the following inductive hypothesis. Suppose that for some $0\le j<k$, there exists a marked centre $\cJ^{(j)}$ that is $j$-semi-associated to $\cI$ at a point $p\in Z\subset Y$, and that $\tilde\cJ^{(j)}:=(\tilde x_1^{a_1},\dots,\tilde x_j^{a_j})$ is $j$-semi-associated to $\tilde \cI$ at the same point $p\in Z\subset B_0$.
	Let $\mx$ be $\cJ$-compatible parameters completing $x_1,\dots,x_j$. Using \autoref{met1}, \eqref{eq: f} and \eqref{eq: f twidle}, we obtain the inequality
	\begin{align}
		a_{j+1}^{BX_0}(p)&=\min_{(\beta,f)}\brc{
			\Xi(\beta): f\in\cI,\;c_\beta^{(f)}\neq0,\; \Delta(\beta)<1 \text{ and }v_\cJ(\mx^\beta)=v_\cJ(f)
		}\label{eq: a_{j+1} BX0}
		\\
		&\geq \min_{(\beta,f)}\brc{
			\Xi\pa{
				\beta
			}
			:f\in\cI,\; c_\beta^{(f)}\neq 0\text{ and } \Delta(\beta)<1
		}
		=a_{j+1}
		\label{eq: a_{j+1} hom}
	\end{align}
	because the index set of the minimum \eqref{eq: a_{j+1} BX0} is a subset of that of the minimum \eqref{eq: a_{j+1} hom}. Let $(\beta,f)$ be a minimiser of the minimum \eqref{eq: a_{j+1} hom}. We have
	\[v_\cJ(\mx^\beta)=\sum_{i=1}^k\frac{\beta_i}{a_i}\leq \sum_{i=1}^j\frac{\beta_i}{a_i}+\frac{\sum_{i=j+1}^n\beta_i}{a_{j+1}}= 1\]
	by substituting $a_{j+1}=\Xi(\beta)$. Thus, by the admissibility of $\cJ$, we must have $v_\cJ(\mx^\beta)=v_{\cJ}(f)=1$. Therefore, $(\beta,f)$ is also a minimiser of the minimum \eqref{eq: a_{j+1} BX0}, implying that $a_{j+1}^{BX_0}(p)=a_{j+1}$.
	
	Choose $l>j$ such that $\beta_l>0$ and let $\bar\beta=\beta-e_l$. We may choose $ y=\partial_\mx^{\bar\beta} f$ 
	as the next parameter of the associated centre to $\cI$ at $p$, and $\partial_{\tilde\mx}^{\bar\beta}\tilde f$ as the next parameter for the associated centre to $\tilde\cI$ at $p$. But by \autoref{lem: diff and hom}, $\tilde y=\partial_{\tilde\mx}^{\bar\beta}\tilde f$, and so the inductive hypothesis holds for $j+1$. 
	
	The inductive hypothesis holding for $j=0$ with $\cJ^{(0)}=()$, we may conclude that $\tilde\cJ=(\tilde x_1^{a_1},\dots,\tilde x_k^{a_k})=(\tilde\cF_\bullet,d)$ is $k$-semi-associated to $\tilde \cI$ at $p$. But by definition, this marked centre is $\tilde \cI$-admissible at $p$, and therefore associated to $\tilde \cI$ at $p$ by \autoref{k-associated implies associated}. 	
	In particular, $\inv_{BX_0}(p)=\inv_X(p)$.
\end{proof}
\subsection{Revisiting the associated centre at a point}

We now present a method for calculating the associated centre that uses derivatives of ideals.
\begin{notation}
	Let $x_1,\dots,x_j$ be local parameters and let $\beta$ be a multi-index of length $l\leq j+1$. We construct the ideals $\cD[\beta;x_1,\dots,x_j,\cI]=\cD[\beta]$ recursively. If $l=1$, we define
	\[
	\cD[\beta]:=\cD^{\leq\beta_1}\cI.
	\]
	Otherwise, we define
	\begin{flalign*}
		&&&&\cD[\beta]:=\cD^{\leq\beta_l}\pa{
			\cD\brk{
				\pa{
					\beta_1,\dots,\beta_{l-1}
				}
			}
			|_{V(x_{l-1})}
		}
		.&&\qedbarhere
	\end{flalign*}
\end{notation}

\begin{lemma}
	Let $\mx$ be a system of parameters at $p$ completing $x_1,\dots,x_j$. Consider any length $n$ multi-indices $\gamma^{(1)}+\cdots+\gamma^{(j+1)}=\gamma$ satisfying $\gamma^{(i)}_l=0$ for $l<i$. We have the equality
	\[
	\partial_\mx^{\gamma^{(j+1)}}\pa{
		\pa{
			\partial_\mx^{\gamma^{(j)}}
			\pa{
				\cdots\pa{
					\partial_{\mx}^{\gamma^{(2)}}\pa{
						\pa{
							\partial_\mx^{\gamma^{(1)}}\cI
						}
						|_{V(x_1)}
					} 
				}
				|_{V(x_2)}\cdots
			}
		}
		|_{V(x_j)}
	}
	=(\partial_\mx^\gamma\cI)|_{V(x_1,\dots,x_j)}.
	\]
\end{lemma}
\begin{proof}
	By linearity of differential operators and restriction, this equality may be checked on monomials where it is trivial.
\end{proof}

\begin{corollary}\label{cor: gen of coeff}
 The following equality of ideals holds:
	\[
	\cD[\beta]=\sum_\gamma\left(\partial_\mx^\gamma\cI\right)|_{V(x_1,\dots,x_j)},
	\]
	where this sum ranges over multi-indices $\gamma$ of length $n$ satisfying
	\begin{itemize}
		\item $\sum_{i=1}^l\gamma_i\leq\sum_{i=1}^l\beta_i$ for all $l\leq j$ and
		\item $\sum_{i=1}^n\gamma_i\leq\sum_{i=1}^{j+1}\beta_i$. 
	\end{itemize}
\end{corollary}
We intend to use the ideals $\cD[\beta]$ similarly to the ideals $\partial_\mx^{\beta}\cI$ in \autoref{met1}. In fact, the following proposition holds analogously to \eqref{eq: a_{j+1} is descr}.
\begin{proposition}\label{prop: a_{j+1} coeff}
	Let $(x_1^{a_1},\dots,x_j^{a_j})$ be $j$-semi-associated at $p$ but not admissible. Then, the next entry in $\inv(p)$ is equal to
	\begin{equation}
		\min_{l(\beta)=j+1}\brc{
			\Xi(\beta) :\Delta(\beta)<1\text{ and }\cD[\beta]=(1)
		}
		.
		\label{eq: min global}	
	\end{equation}
\end{proposition}
It is clear that \eqref{eq: min global} is less than or equal to
\begin{equation}
	a_{j+1}=\min_{l(\alpha)=n}\brc{
		\Xi(\alpha):\Delta(\alpha)<1\text{ and }\partial_\mx^\alpha\cI=(1)
	}
	,
	\label{eq: a_{j+1} v1}
\end{equation}
which we know to be the next entry in the invariant by using \autoref{met1}. Indeed, if $\alpha$ is a minimiser of \eqref{eq: a_{j+1} v1}, then we have 
\[
(1)=(\partial_\mx^\alpha\cI)|_{V(x_1,\dots,x_j)}\subset\cD[\beta]
\]
where $\beta=(\alpha_1,\dots,\alpha_j,\sum_{i=j+1}^n\alpha_i)$. Therefore, $\cD[\beta]=(1)$ and $\Xi(\beta)=\Xi(\alpha)=a_{j+1}$ so that this side of the inequality holds. 

The other inequality\,---\,namely, that \eqref{eq: min global} $\geq$ \eqref{eq: a_{j+1} v1} — is more challenging to prove. The following numerical inequality is crucial for establishing it.
\begin{theorem}\label{thm: num}
	Let $\beta$ and $\gamma$ be multi-indices of length $j+1$ satisfying
	\begin{flalign*}
		&&&&\sum_{i=1}^l\gamma_i\leq\sum_{i=1}^l\beta_i&&\text{for $l=1,\dots,j+1$.}
	\end{flalign*}
	Let $(b_1,\dots,b_j)$ be a pre-invariant. Then, either
	\begin{itemize}
		\item $\Xi(b_1,\dots,b_j;\gamma)\leq\Xi(b_1,\dots,b_j;\beta)$, or
		\item $\Xi(b_1,\dots,b_j;\gamma)<b_j$.
	\end{itemize}
\end{theorem}
The proof of \autoref{thm: num} is elementary but rather technical and is deferred to \autoref{subsec: numerical inequalities}. 
Applying \autoref{thm: num} to our situation, the generators $(\partial_\mx^\gamma\cI)|_{V(x_1,\dots,x_j)}$ of $\cD[\beta]$ from \autoref{cor: gen of coeff} must satisfy either $\Xi(\gamma)\leq \Xi(\beta)$ or $\Xi(\gamma)<a_j$. In particular, if $\Xi(\beta)<a_{j+1}$, then $\Xi(\gamma)<a_{j+1}$.
\begin{proof}[Proof of \autoref{prop: a_{j+1} coeff}]
	It remains to be shown that \eqref{eq: min global} is greater than or equal to \eqref{eq: a_{j+1} v1}. Suppose that $\Xi(\beta)<a_{j+1}$. For any generator $(\partial_\mx^\gamma\cI)|_{V(x_1,\dots,x_j)}$ of $\cD[\beta]$ as above, \autoref{thm: num} tells us that $\Xi(\gamma)< a_{j+1}$. Therefore, by the minimality of $a_{j+1}$, $\partial_\mx^\gamma\cI\neq (1)$ and, equivalently, $(\partial_\mx^\gamma\cI)|_{V(x_1,\dots,x_j)}\neq (1)$. Since we work at a point, this means $\cD[\beta]\neq(1)$.
\end{proof}
We next establish that compatible parameters for the associated centre may likewise be constructed in this framework.
\begin{proposition}
	Let $\beta$ be a minimiser of \eqref{eq: min global}, and let $x'_{j+1}$ be any maximal contact element of $\cD[(\beta_1,\dots,\beta_j,\beta_{j+1}-1)]$.
	Any lift $x_{j+1}$ of $x_{j+1}'$ to $\cO_Y$ can be taken as the next compatible parameter in the associated centre at $p$.	
\end{proposition}
\begin{proof}
	Let $\cJ=(x_1^{a_1},\dots,x_j^{a_j},y_{j+1}^{a_{j+1}},\dots,y_k^{a_k})$ be the associated centre at $p$. The element $x_{j+1}$ has a linear term in $y_{j+1},\dots, y_k$ so that $v_\cJ(x_{j+1})\leq 1/a_{j+1}$.
	
	Conversely, since any lift of $x'_{j+1}$ has the same valuation (we are playing with terms containing $x_1,\dots,x_j$ and $v_\cJ(x_{j+1})\leq 1/a_{j+1}$), we may assume
	\[x_{j+1}\in\sum_\gamma\partial_\mx^\gamma\cI,\]
	where the sum ranges over multi-indices $\gamma$ of length $n$ satisfying $\sum_{i=1}^l\gamma_i\leq\sum_{i=1}^l\beta_i$ for $l\leq j$ and $\sum_{i=1}^n\gamma_i\leq(\sum_{i=1}^{j+1}\beta_i)-1$. These conditions on $\gamma$ and admissibility of $\cJ$ imply that
	\[
	v_{\cJ}(\partial_\mx^\gamma\cI)\geq 1-\sum_{i=1}^k\frac{\gamma_i}{a_i}\geq 1-\sum_{i=1}^j\frac{\beta_i}{a_i}-\frac{\beta_{j+1}-1}{a_{j+1}}=\frac{1}{a_{j+1}},
	\]
	by substituting $a_{j+1}=\Xi(\beta)$. Thus, we may conclude that $v_\cJ(x_{j+1})= 1/a_{j+1}$. By \autoref{lem:equal-weightings}, $x_{j+1}$ may be taken as the next parameter in the associated centre.
\end{proof}

\begin{method}\label{met2}
	The above proofs may be summarised in a new procedure for constructing the associated centre at a point. This is \hyperref[subsec: method 2]{Method 2} from the \hyperref[sec: intro]{Introduction}.
	
	Let $\cJ^{(j)}=(x_1^{a_1},\dots,x_j^{a_j})$ be $j$-semi-associated. Then, the next entry in the invariant is
	\begin{equation}
		a_{j+1}=\min_{l(\beta)=j+1}\brc{
			\Xi(\beta) :\Delta(\beta)<1\text{ and }\cD[\beta]=(1)
		}
		.
		\label{eq: a_{j+1} global pt 2}
	\end{equation}
	Let $\beta$ be a minimiser and let $x_{j+1}$ be the lift of a maximal contact element of $\cD[(\beta_1,\dots,\beta_j,\beta_{j+1}-1)]$ at $p$. The marked centre $\cJ^{(j+1)}=(x_1^{a_1},\dots,x_{j+1}^{a_{j+1}})$ is $(j+1)$-semi-associated.
	The procedure stops when $\cJ^{(k)}=\cJ$ is admissible.
	
	Note that one could instead choose as the next parameter the lift of a maximal contact element of the larger ideal
	\[
	\pa{
		\cD\brk{
			\pa{
				\beta_1,\dots,\beta_j,\beta_{j+1}-1
			}
		}
		: \Xi(\beta)=a_{j+1}
	}
	\]
	since the same proof shows that any such lift has the correct valuation. For technical reasons, it is necessary to mention this subtlety to outline a global algorithm, as we do in \autoref{const: global algo}.
\end{method}

There is a straightforward criterion to check whether $\cJ^{(j)}$ is admissible in this framework. Let $\mx$ be a system of parameters completing $x_1,\dots,x_j$. Recall that admissibility can be phrased as the following implication: if $\Delta(\beta)<1$, then $\partial_\mx^\beta\cI\neq(1)$. 

Suppose that $\cJ^{(j)}$ is admissible. Let $\beta$ be a multi-index of length $j$ with $\Delta(\beta)<1$. Suppose that $\cD[\beta]|_{V(x_j)}\neq (0)$. Then, there exists an integer $\beta_{j+1}$ such that $\cD[(\beta_1,\dots,\beta_{j+1})]=(1)$. By \autoref{cor: gen of coeff}, $\cD[(\beta_1,\dots,\beta_{j+1})]$ is generated by ideals $\partial_\mx^\gamma\cI|_{V(x_1,\dots,x_j)}$ for some multi-indices $\gamma$ of length $n$ satisfying inequalities which in particular imply the inequality $\Delta(\gamma)\leq\Delta(\beta)<1$. Since we are at a point, there must be one of these ideals such that $\partial_\mx^\gamma\cI|_{V(x_1,\dots,x_j)}=(1)$, or equivalently $\partial_\mx^\gamma\cI=(1)$. But this contradicts admissibility, so we actually have $\cD[\beta]|_{V(x_j)}=(0)$.

Conversely, suppose that for all $\beta$ of length $j$ and $\Delta(\beta)<1$ we have $\cD[\beta]|_{V(x_j)}=(0)$. Let $\gamma$ be a multi-index of length $n$ with $\Delta(\gamma)<1$. Let $b=\sum_{i=j+1}^n\gamma_i$. Since 
\[
(\partial_\mx^{\gamma}\cI)|_{V(x_1,\dots,x_j)}\subset\cD\brk{
	\pa{
		\gamma_1,\dots,\gamma_j,b
	}
}
\subset \cD\brk{
	\pa{
		\gamma_1,\dots,\gamma_j
	}
}
|_{V(x_j)}=(0),
\]
$\partial_\mx^\gamma\cI$ vanishes at $p$, so that $\cJ^{(j)}$ is admissible.
Therefore, we obtain the following proposition:
\begin{proposition}
	$\cJ^{(j)}$ is admissible if and only if for all multi-indices $\beta$ of length $j$ such that $\Delta(\beta)<1$, we have $\cD[\beta]|_{V(x_j)}=(0)$.
\end{proposition}
\begin{construction}[Sketch of a global algorithm]\label{const: global algo}
	The main advantage of \autoref{met2} is that it simplifies the description of the invariant’s level sets, thereby making it easier to formulate a global algorithm. Indeed, consider the following inductive hypothesis: suppose that we have a chart $U\subset Y$ and a marked centre $(x_1^{a_1},\dots,x_j^{a_j})$ on $U$ such that if $q\in U$ has invariant 
	\[
	\inv(q)\geq (\underbrace{a_1,\dots,a_j,\dots, a_j}_{\text{length }n}),
	\]
	then $\cJ(q)$ dominates $(x_1^{a_1},\dots,x_j^{a_j})$. It can be shown by induction that for $b\geq a_j$, the locus where
	\[
	\inv(q)\geq (\underbrace{a_1,\dots,a_j,b,\dots, b}_{\text{length }n}),
	\]
	is cut out by the ideal 
	\[
	\cA\!\pa{
		a_1,\dots,a_j;b
	}
	:=\pa{
		\cD\brk{
			\beta
		}
		:l(\beta)=j+1,\; 0\leq \Xi(\beta)<b
	}
	.
	\]
	Therefore, the next entry in $\maxinv_X$ is
	\[
	a_{j+1}=\max\brc{
		b: \cA\!\pa{
			a_1,\dots,a_j;b
		}
		\neq (1)
	}
	,
	\]
	where inequality of ideals is on $V(x_1,\dots,x_k)$. Let $\cA=\cA(a_1,\dots,a_j;a_{j+1})$. The locus $V(\cA)\subset V(x_1,\dots,x_j)\subset U$ is that where $\inv(q)$ dominates $(a_1,\dots,a_{j+1})$. Choose charts $U_l$ and lifts $x_{j+1}^{(l)}\in\cO_{U_l}$ of maximal contact elements of the ideal
	\[
	\pa{
		\cD\brk{
			\pa{
				\beta_1,\dots,\beta_j,\beta_{j+1}-1
			}
		}
		:\Xi(\beta)=a_{j+1}
	}
	\]
	such that each $x_{j+1}^{(l)}$ vanishes on $V(\cA)\cap U_l$ and $V(x_{j+1}^{(l)})\subset U_l$ is smooth. Then, each chart $U_l$, together with its marked centre 
	\[
	\pa{x_1^{a_1},\dots,x_j^{a_j},\pa{x_{j+1}^{(l)}}^{a_{j+1}}},
	\]
	satisfies the inductive hypothesis for $j+1$. In particular, if we start at step $0$ with the single chart $Y$ and the marked centre $()$, and stop once an admissible marked centre is achieved (discarding any components of the subvariety on which the resulting marked centre is not admissible), the outcome is the associated centre.
	
	\cite{lee2020algorithmicresolutionweightedblowings} has provided a pseudo-code for the construction of the associated centre following \cite{abramovich2024functorial} for which there exists a running implementation in the computational algebra program {\sc singular} \cite{DGPS}; see \cite{resweightedlib}. Replacing their algorithm \texttt{prepareCenter} with the construction we have sketched above would yield a functional implementation of the algorithm based on \autoref{met2}.
\end{construction}
\subsection{Examples}
We now compute the associated centre in two examples: one local and one global. We use \autoref{met1} for the first and \autoref{met2} for the second.
\begin{example}[$A_m$-singularity]
	Let $m\geq 1$ and let $f=x^2-y^{m+1}$. Consider $\cI=(f)\subset\cO_{\bb A^2}$. It cuts out a curve in $Y=\bb A^2$ with a singularity of type $A_m$ at the origin. We wish to compute the associated centre at the origin. As the first parameter, we may take $\frac{1}{2}\partial_{x}f=x$, and so $(x^{2})$ is $1$-semi-associated. We carry out the second step using the same system of parameters $(x,y)$. The only multi-index $\beta=(\beta_1,\beta_2)$ with $\Xi(2;\beta)<1$ such that $\partial^\beta f$ does not vanish at the origin is $\beta=(0,m+1)$. Therefore, we may take $y=\frac{1}{(m+1)!}\partial_{y}^m f$ as the next parameter. The associated marked centre is thus $(x^2,y^{m+1})$. If $m$ is even, the weights of the associated centre are $(m+1,2)$; if $m$ is odd, the weights are $(\frac{m+1}{2},1)$. Note that when $m=1,2$, we recover Examples \ref{ex: gen in deg 1} and \ref{ex: weighting associated to cusp} respectively.
\end{example}
\begin{example}[Whitney umbrella]\label{eg: whitney umbrella} Let $f=x^2-y^2z$, and consider the ideal $\cI=(f)\subset\cO_{\bb A^3}$. This cuts out the Whitney umbrella $X\subset\bb A^3$. We compute the global associated centre. We find $\cD^{\leq 2}\cI=(1)$ and $\cD^{\leq 1}\cI=(x,y^2,yz)$. Thus, we have $a_1=2$ and the maximal contact element $x$ can be chosen as the first parameter. 
	
	For the second step, only $\beta_1=0,1$ are considered. We have $\cD[(1,0)]=(\cD^{\leq 1}\cI)|_{x=0}=(y^2,yz)$ and $\cD[(0,0)]=\cI|_{x=0}=(y^2z)$, with $\maxord\;\cD[(1,0)]=2$ and $\maxord\;\cD[(0,0)]=3$. Since
	\[\frac{2}{1-\frac{1}{2}}=4>\frac{3}{1-0}=3,\]
	we conclude that $a_2=3$. Both $y$ and $z$ can therefore serve as the next parameters. Hence, the associated marked centre is $(x^2,y^3,z^3)$. The associated centre has compatible parameters $(x,y,z)$ with respective weights $(3,2,2)$.
	
	We can also compute the associated centre to $X\setminus\{\mbf 0\}\subset \bb A^3\setminus\{\mbf 0\}$ to see how the singularities behave away from the origin. We again have maximal order $2$, but we now have the equality $\cD^{\leq 1}\cI|_{\bb A^3\setminus\{0\}}=(x,y^2,yz)|_{\bb A^3\setminus\{0\}}=(x,y)$. This means that both $x$ and $y$ may be taken as the first parameter, and $(x^2,y^2)$ is admissible, and therefore associated. Away from the origin, the maximal invariant is thus $(2,2)$, attained on the line $V(x,y)$. $X$ is smooth away from this line so that the invariant there is $(1)$. 
\end{example}
%

%
% !TEX root = main.tex
\section{Efficiency and comparison with previous constructions}\label{sec: eff}

In this section, we wish to compare the efficiency of our techniques with those of \cite{abramovich2024functorial}.

In \cite{abramovich2024functorial}, the associated centre is defined using coefficient ideals, a construction which we now recall.

\begin{definition}
	Let $\cI\subset\cO_Y$ be an ideal and let $b$ be a positive integer. We define the \tbf{$b$th coefficient ideal of $\cI$} to be the ideal
	\begin{flalign}
		&&&&\call C(\cI,b):=\sum_{i=0}^{b-1}\pa{
			\call D^{\leq i}\cI
		}
		^{\frac{b!}{b-i}}.&&\qedbarhere\label{eq:coeff-id}
	\end{flalign}
\end{definition}
\begin{remark}
	In \cite{abramovich2024functorial}, a larger ``saturated'' version of the coefficient ideal is used; see also \cite{kollar2007resolutionsingularitiesseattle, abramovich2020relativedesingularizationprincipalizationideals}. However, as mentioned in \cite[5.1]{lee2020algorithmicresolutionweightedblowings}, both versions are equivalent for computational purposes.
\end{remark}
\begin{remark}
	Various notions of coefficient ideals have been explored in several contexts throughout the years; see \cite{villamayor1989constructiveness, Bierstone1991, Bierstone1997, wlodarzcyk2005simple, kollar2007resolutionsingularitiesseattle}. Earlier definitions of coefficient ideals were expressed within the framework of $\bb Q$-ideals, which avoids the use of factorials; indeed, the exponents in \eqref{eq:coeff-id} serve to tune the weight of each summand appropriately, and the factorials arise when avoiding rational powers. In that regard, the original approach of \cite{Bierstone1997} is absent of factorial growth.\footnote{We thank Bernd Schober and Dan Abramovich for explaining this.} Recently, \cite{wlodarczyk2023functorialresolutiontorusactions,Włodarczyk2023} presented the results of \cite{abramovich2024functorial} in a similar framework, thereby achieving computational gains before the present paper. One should also look at \cite{Schober2013, schober2014polyhedralapproachinvariantbierstone, schoberhircharpol} for treatments of the Bierstone--Milman invariant avoiding factorial growth.
\end{remark}
\begin{construction}[Associated marked centre following \cite{abramovich2024functorial}]\label{con: ass centre coeff}
	Let $\cI=\cI_X\subset\cO_Y$ be a proper and nowhere zero ideal, and let $p\in X$ be a point. We define $\cI[1]:=\cI$, let $b_1=\ord_p\cI$, and choose $x_1$ to be a maximal contact element of $\cI$ at $p$. For the second step, we define \[\cI[2]:=\call C(\cI,b_1)|_{V(x_1)},
	\]
	let $b_2=\ord_p\cI[2]$ and choose $x_2$ to be the lift of a maximal contact element of $\cI[2]$ at $p$. We continue this process inductively, defining 
	\[
	\cI[j+1]:=\call C(\cI[j],b_j)|_{V(x_j)}
	\]
	and choosing maximal contact elements. We stop when $\cI[k+1]=(0)$. For all $j\leq k$, we define 
	\[
	a_{j}:=\frac{b_{j}}{(b_{j-1}-1)!\cdots(b_1-1)!}.
	\]
	The associated marked centre at $p$ is then $(x_1^{a_1},\dots,x_k^{a_k})$.
\end{construction}

\begin{remark}The numbers $b_j$ explode quickly. For example, the seemingly innocuous invariant $(a_1,\dots,a_5)=(2,3,3,\frac{20}{3},9)$ becomes 
	$(b_1,\dots,b_5)=(2,3,6,1600,2160\cdot1599!)$.
	
	In particular, the number of derivatives required to calculate the coefficient ideals increases factorially as the algorithm progresses.  In \autoref{met2}, the role of the coefficient ideals is played by the simpler ideals $\cD[\beta]$, for which the number of derivatives grows much more slowly.
\end{remark}

\begin{example}\label{ex: comparison efficiency}
	Consider the equation $x^4+y^5+z^6=0$ in $\bb A^3$.
	Let us first compute the associated marked centre at the origin $\mbf 0$ with \autoref{met2}.
	
	At the first step, we differentiate $\cI=(x^4+y^5+z^6)$:
	\begin{itemize}
		\item $\cD^{\leq 0}\cI=\cI=(x^4+y^5+z^6)$;
		\item $\cD^{\leq 1}\cI=(x^3,y^4,z^5)$;
		\item $\cD^{\leq 2}\cI=(x^2,y^3,z^4)$;
		\item $\cD^{\leq 3}\cI=(x,y^2,z^3)$;
		\item $\cD^{\leq 4}\cI=(1)$.
	\end{itemize}
	Thus, $\cI$ has order $4$ at $\mbf 0$, and $x$ may be taken as the first parameter. Now consider $(\cD^{\leq\beta_1}\cI)|_{V(x)}$ for $\beta_1<4$. We can observe that it has order $5-\beta_1$ at $\mbf 0$. Thus,
	\[
	a_2=\min_{\beta_1<4}\brc{
		\frac{5-\beta_1}{1-\frac{\beta_1}{4}}
	}
	=5,
	\]
	with minimiser $\beta_1=0$. Since $y$ is a (lift of a) maximal contact element of $(\cD^{\leq 0}\cI)|_{V(x)}=(y^5+z^6)$, we take it as the next parameter. We now look at the non-negative integral solutions of $\frac{\beta_1}{4}+\frac{\beta_2}{5}<1$. There are fourteen of these:
	\begin{itemize}
		\item $\beta_1=0$ and $\beta_2=0,1,2,3,4$;
		\item $\beta_1=1$ and $\beta_2=0,1,2,3$;
		\item $\beta_1=2$ and $\beta_2=0,1,2$;
		\item $\beta_1=3$ and $\beta_2=0,1$.
	\end{itemize}
	We see that for each of these, the order of $\cD[(\beta_1,\beta_2,0)]$ at $\mbf 0$ is $6-\beta_1-\beta_2$. We thus have
	\[
	a_3=\min_{\frac{\beta_1}{4}+\frac{\beta_2}{5}<1}\brc{
		\frac{6-\beta_1-\beta_2}{1-\frac{\beta_1}{a_1}-\frac{\beta_2}{a_2}}
	}
	=6,
	\]
	with minimiser $\beta_1=\beta_2=0$. We choose $z$ as the last parameter, as it is (the lift of) an element of maximal contact of $\cI|_{V(x,y)}=(z^6)$. The associated marked centre is $(x^4,y^5,z^6)$.
	
	Let us now compute the same associated marked centre using the procedure of \cite{abramovich2024functorial}. The first steps are the same; we can take $x$ as the first parameter and $a_1=4$. Next, we compute 
	\begin{align*}
		\cI[2]&=\call C(\cI,3)|_{x=0}=\sum_{0\leq i\leq 3}\pa{
			\call D^{\leq i}\cI
		}
		|_{x=0}^{\frac{4!}{4-i}}\\
		&=\pa{
			y^5+z^6
		}
		^{6}+\pa{
			y^4,z^5
		}
		^{8}+\pa{
			y^3,z^4
		}
		^{12}+\pa{
			y^2,z^3
		}
		^{24}=\pa{
			\pa{
				y^5+z^6
			}
			^6,y^{32},z^{40}
		}
		.
	\end{align*}
	This ideal has order $b_2=30$ at $\mbf 0$ and has $y$ as a maximal contact element, which we take as the next parameter.
	We compute
	\[
	\cI[3]=\call C\pa{
		\cI[2],30
	}
	|_{x=y=0}=\sum_{0\leq i\leq 29}\pa{
		\call D^{\leq i}\cI[2]
	}
	|_{x=y=0}^{\frac{30!}{30-i}}.
	\]	 	
	With a little thought, we see that $(\call D^{\leq i}\cI[2])|_{x=y=0}=(z^{36-i})$. Thus, $\cI[3]=(z^{36})^{29!}$, implying that $b_3=36\cdot 29!$ and $z$ may be taken as the last parameter. However, many computational algebra programs output an error when calculating this ideal because the exponents are too large. 
\end{example}

At each step of \autoref{met2}, we need to compute the derivatives $\cD[\beta]|_{V(x_j)}$ of $\cI$ for each multi-index $\beta$ of length $j$ with $\Delta(a_1,\dots,a_j;\beta)<1$. These multi-indices form a simplex, which we denote by $\Sigma$. We then have to compute the order at $p$ and possibly maximal contacts for each $\beta\in \Sigma$. Therefore, understanding the size $\sigma(a_1,\dots,a_j)$ of $\Sigma$ informs us directly about the algorithm's complexity.
\begin{proposition}
	The following inequality holds: 
	\[
	\sigma(a_1,\dots,a_j)\leq \sum_{A\subset\{1,\dots,j\}}\frac{1}{|A|!}\prod_{l\in A}a_l.
	\]
\end{proposition}

\begin{proof}
	\cite[(1.5)]{yau2006upper} states that the interior of $\Sigma$, i.e. the set $\{\beta_1,\dots,\beta_j\neq 0\}\subset \Sigma$, has size bounded by $\frac{1}{j!}\cdot a_1\cdots a_j$. Applying the same inequality to the interior of the simplex 
	\[
	\{\beta_l=0\text{ for }l\notin A\}\subset \Sigma
	\]
	for each proper subset $A\subset\{1,\dots,j\}$  covers all other cases.
\end{proof}
The function $\sigma(a_1,\dots,a_j)$ has a very controllable growth, being bounded by a multi-linear function in $a_1,\dots,a_j$ whose coefficients $\frac{1}{|A|!}$ decay rapidly as $|A|$ grows. This gives a bound on the complexity of the algorithm in terms of the invariant of the singularity.  Note, however, that since the algorithm is used to calculate the invariant, this bound is of limited use absent some additional a priori bound on the invariant in a given class of examples.
%

%
% !TEX root = main.tex
\appendix
\section{Numerical inequalities}
\label{subsec: numerical inequalities}
The purpose of this appendix is to prove the following numerical theorem, which was crucial in the development of \autoref{met2}.

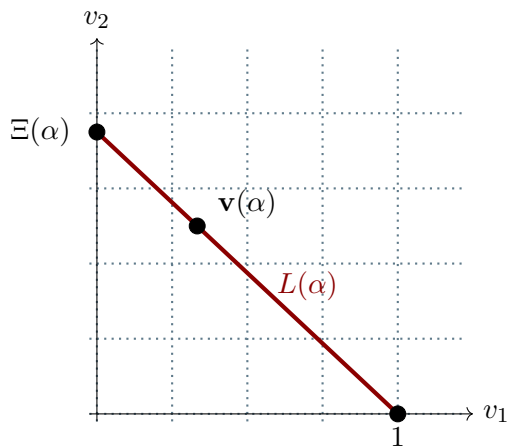
\begin{figure}[htbp]
	\scalebox{1}{
		\begin{tikzpicture}
			\draw[thick,color=LightSkyBlue4, dotted] (-0.1,-0.1) grid (4.9,4.9);
			\draw[->] (-0.1,0) -- (5.0,0) node[right] {$v_1$};
			\draw[ultra thick, -, color=Red4] (4,0)--(0,3.75);
			\draw[->] (0,-0.1) -- (0,5.0) node[above] {$v_2$};
			\node at (4,-0.3) {$1$};
			\filldraw (4,0) circle (3pt);
			\filldraw (4/3,5/2) circle (3pt);
			\filldraw (0,3.75) circle (3pt);
			\node at (-.75,3.75) {$\Xi(\alpha)$};
			\node at (2.8,1.7) {$\color{Red4} L(\alpha)$};
			\node at (2,2.8) {$\mbf v(\alpha)$};
		\end{tikzpicture}
	}
	\caption{$\Xi(\alpha)$ is the $v_2$-intercept of $L(\alpha)$}\label{fig: int v(alpha)}
\end{figure}

\begin{reptheorem}{thm: num}
	Let $\beta$ and $\gamma$ be multi-indices of length $j+1$ satisfying
	\begin{flalign}
		&&&&\sum_{i=1}^l\gamma_i\leq\sum_{i=1}^l\beta_i&&\text{for $l=1,\dots,j+1$.}\label{eq: Cond1}
	\end{flalign}
	Let $(b_1,\dots,b_j)$ be a pre-invariant. If \[\Xi(b_1,\dots,b_j;\gamma)>\Xi(b_1,\dots,b_j;\beta),\] then \[\Xi(b_1,\dots,b_j;\gamma)<b_j.\]
\end{reptheorem}
In order to prove this theorem, we first arrange things so that the situation may be drawn in the plane by collapsing some dimensions. Let $\alpha$ be a pre-invariant of length $j+1$. Let
\[
\mbf v(\alpha):=\pa{\sum_{i=1}^j\frac{\alpha_i}{b_i},\alpha_{j+1}}
\]
in the $\mbf v=(v_1,v_2)$-plane.
This way, the $v_2$-intercept of the line $L(\alpha)$ passing through $(1,0)$ and $\mbf v(\alpha)$ is $\Xi(b_1,\dots,b_j;\alpha)$, as sketched in \autoref{fig: int v(alpha)}.

We now interpret \autoref{thm: num} in the plane. The conditions \eqref{eq: Cond1} imply that $\sum_{i=1}^j\frac{\gamma_i}{b_i}<\sum_{i=1}^j\frac{\beta_i}{b_i}$ since $b_1\leq\cdots\leq b_j$, so that $\mbf v(\gamma)$ lies to the left of $\mbf v(\beta)$. \autoref{thm: num} says that if \eqref{eq: Cond1} holds and if the $v_2$-intercept of the line $L(\gamma)$ is greater than that of $L(\beta)$, then it is less than $b_j$. \autoref{fig3} illustrates this.
\begin{figure}[htbp]
	\scalebox{1}{
		\begin{tikzpicture}
			\draw[thick,color=LightSkyBlue4, dotted] (-0.1,-0.1) grid (4.9,5.9);
			\draw[->] (-0.1,0) -- (5.0,0) node[right] {$v_1$};
			\draw[->] (0,-0.1) -- (0,6.0) node[above] {$v_2$};
			\node at (-0.3,5) {$b_j$};
			\node at (4,-0.3) {$1$};
			\draw[ultra thick, -, color=Red4] (4,0)--(0,8/5);
			\draw[ultra thick, -, color=Blue4] (4,0)--(0,2+2/3);
			\node at (0.5,2.7) {$\color{Blue4}L(\gamma)$};
			\node at (0.5,1.7) {$\color{Red4}L(\beta)$};
			\filldraw (3,2/5) circle (3pt);
			\node at (2.6,0.2) {$\mbf v(\beta)$};
			\filldraw (1,2) circle (3pt);
			\node at (1.3,2.3) {$\mbf v(\gamma)$};
		\end{tikzpicture}
	}
	\caption{Illustration of \autoref{thm: num}} \label{fig3}
\end{figure}
The following lemma, an extremal case of \autoref{thm: num},  shall be established first, as it can be leveraged for a proof.
\begin{lemma}\label{lem: num}
	Let $\beta$ be a multi-index of length $j+1$ and let $(b_1,\dots,b_j)$ be a pre-invariant. If $\;\sum_{i=1}^{j+1}\beta_i>\Xi(b_1,\dots,b_j;\beta)$, then $\;\sum_{i=1}^{j+1}\beta_i<b_j$.
\end{lemma}
\begin{proof}
	We have the following equivalence of inequalities:
	\begin{align}
		\sum_{i=1}^{j+1}\beta_i&>\Xi(b_1,\dots,b_j;\beta)=\frac{\beta_{j+1}}{1-\sum_{i=1}^j\frac{\beta_i}{b_i}}\notag\\
		\pa{
			\sum_{i=1}^{j+1}\beta_i
		}
		\pa{
			1-\sum_{i=1}^j\frac{\beta_i}{b_i}
		}
		&>\beta_{j+1}\notag\\
		\sum_{i=1}^j\beta_i-\pa{
			\sum_{i=1}^{j+1}\beta_i
		}
		\pa{
			\sum_{i=1}^j\frac{\beta_i}{b_i}
		}
		&>0\notag\\
		\sum_{i=1}^j\beta_i&>\pa{
			\sum_{i=1}^{j+1}\beta_i
		}
		\pa{
			\sum_{i=1}^j\frac{\beta_i}{b_i}
		}
		\notag\\
		\frac{\sum_{i=1}^j\beta_i}{\sum_{i=1}^{j+1}\beta_i}&>\sum_{i=1}^j\frac{\beta_i}{b_i}.\label{eq: ineq eq}
	\end{align}
	Now suppose by contradiction that $\sum_{i=1}^{j+1}\beta_i\geq b_j$. Then \eqref{eq: ineq eq} implies that
	\begin{align}
		\frac{\sum_{i=1}^j\beta_i}{b_j}>\sum_{i=1}^j\frac{\beta_i}{b_i},\notag
	\end{align} which is impossible because $b_1\leq \cdots\leq b_j$. Thus, $\sum_{i=1}^{j+1}\beta_i<b_j$.
\end{proof}
Note that if $\mbf v(\beta)=(0,\sum_{i=1}^{j+1}\beta_i)$, then $\beta_1,\dots,\beta_j=0$. Therefore, if $\gamma$ satisfies \eqref{eq: Cond1}, we also have $\gamma_1,\dots,\gamma_j=0$, implying that $\Xi(b_1,\dots,b_j,\gamma)=\gamma_{j+1}\leq \beta_{j+1}=\Xi(b_1,\dots,b_j,\beta)$, and so \autoref{thm: num} holds vacuously. In order to prove \autoref{thm: num}, we may therefore assume that $\mbf v(\beta)\neq(0,\sum_{i=1}^{j+1}\beta_i)$, so that the line $M$ from $(0,\sum_{i=1}^{j+1}\beta_i)$ to $\mbf v(\beta)$ may be considered. Recall that $\mbf v(\gamma)$ lies to the left of $\mbf v(\beta)$. Therefore, if we show that the conditions \eqref{eq: Cond1} imply that $\mbf v(\gamma)$ lies under $M$, it will imply
\begin{equation}
	\Xi(b_1,\dots,b_j;\gamma)\leq\max\brc{
		\Xi(b_1,\dots,b_j;\beta),\sum_{i=1}^{j+1}\beta_i
	}
	.
	\label{eq: leq max}
\end{equation}
\autoref{fig4} illustrates this. The green line is $M$, and the red and blue lines are, respectively, $L(\beta)$ and $L(\gamma)$ as before. It explains better than words why $\mbf v(\gamma)$ lying under $M$ implies that the $v_2$-intercept of $L(\gamma)$ is either bounded from above by the $v_2$-intercept of $L(\beta)$, or by $\sum_{i=1}^{j+1}\beta_i$.
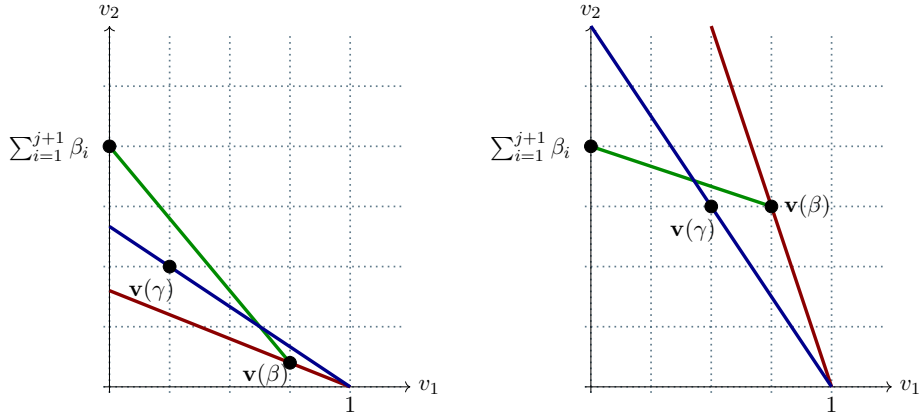
\begin{figure}[htbp]
	\scalebox{.8}{
		\begin{tikzpicture}
			
			\draw[thick,color=LightSkyBlue4, dotted] (-0.1,-0.1) grid (4.9,5.9);
			\draw[->] (-0.1,0) -- (5.0,0) node[right] {$v_1$};
			\draw[->] (0,-0.1) -- (0,6.0) node[above] {$v_2$};
			\draw[ultra thick, -, color=Green4] (0,4)--(3,2/5);
			\draw[ultra thick, -, color=Red4] (4,0)--(0,8/5);
			\draw[ultra thick, -, color=Blue4] (4,0)--(0,2+2/3);
			\node at (-1,4) {$\sum_{i=1}^{j+1}\beta_i$};
			\filldraw (0,4) circle (3pt);
			\node at (4,-0.3) {$1$};
			\filldraw (3,2/5) circle (3pt);
			\node at (2.6,0.2) {$\mathbf{v}(\beta)$};
			\filldraw (1,2) circle (3pt);
			\node at (0.7,1.6) {$\mathbf{v}(\gamma)$};

			\begin{scope}[xshift=8cm] 
				\draw[thick,color=LightSkyBlue4, dotted] (-0.1,-0.1) grid (4.9,5.9);
				\draw[->] (-0.1,0) -- (5.0,0) node[right] {$v_1$};
				\draw[->] (0,-0.1) -- (0,6.0) node[above] {$v_2$};
				\draw[ultra thick, -, color=Green4] (0,4)--(3,3);
				\draw[ultra thick, -, color=Red4] (4,0)--(2,6);
				\draw[ultra thick, -, color=Blue4] (4,0)--(0,6);
				\node at (-1,4) {$\sum_{i=1}^{j+1}\beta_i$};
				\filldraw (0,4) circle (3pt);
				\node at (4,-0.3) {$1$};
				\filldraw (3,3) circle (3pt);
				\node at (3.6,3) {$\mathbf{v}(\beta)$};
				\filldraw (2,3) circle (3pt);
				\node at (1.7,2.6) {$\mathbf{v}(\gamma)$};
			\end{scope}
		\end{tikzpicture}
	}
	\caption{Two illustrations of the inequality \eqref{eq: leq max} when $\mathbf{v}(\gamma)$ lies under $M$} \label{fig4}
\end{figure}

Using \autoref{lem: num}, inequality \eqref{eq: leq max} immediately implies \autoref{thm: num}. To show \eqref{eq: leq max}, note first that $M$ is given by the equation 
\[
v_2=\sum_{i=1}^{j+1}\beta_i - \frac{\sum_{i=1}^j\beta_i}{\sum_{i=1}^j\frac{\beta_i}{b_i}}v_1.
\]
On the other hand, condition \eqref{eq: Cond1} for $l=j+1$ states exactly that
\[
\gamma_{j+1}\leq \sum_{i=1}^{j+1}\beta_i-\sum_{i=1}^j\gamma_i
\]
In the $(v_1,v_2)$-plane, this means that $\mbf v(\gamma)$ satisfies
\begin{equation}
	v_2\leq \sum_{i=1}^{j+1}\beta_i-\frac{\sum_{i=1}^j\gamma_i}{\sum_{i=1}^j\frac{\gamma_i}{b_i}}v_1\label{eq: Cond2}
\end{equation}
Thus, to establish that $\mbf v(\gamma)$ lies under $M$, it suffices to show that for any $\gamma_1,\dots,\gamma_j$ such that 
\begin{flalign}
	&&&&\sum_{i=1}^l\gamma_i\leq\sum_{i=1}^l\beta_i&&\text{for $l=1,\dots,j$,}\label{eq: Cond3}
\end{flalign}
we have the following inequality:
\begin{equation}
	\frac{\sum_{i=1}^j\gamma_i}{\sum_{i=1}^j\frac{\gamma_i}{b_i}}\geq\frac{\sum_{i=1}^j\beta_i}{\sum_{i=1}^j\frac{\beta_i}{b_i}}.\label{eq: fdt num ineq}
\end{equation}
\begin{proof}[Proof of \autoref{thm: num}]
	We need to show that the conditions \eqref{eq: Cond3} imply the inequality \eqref{eq: fdt num ineq}. Suppose first that $\sum_{i=1}^j\gamma_i=\sum_{i=1}^j\beta_i$. As observed before, the conditions \eqref{eq: Cond3} imply that $\sum_{i=1}^j\frac{\gamma_i}{b_i}\leq\sum_{i=1}^j\frac{\beta_i}{b_i}$ since $b_1\leq\cdots\leq b_j$. Therefore \eqref{eq: fdt num ineq} holds in this case.
	
	Suppose otherwise that $\sum_{i=1}^j\gamma_i<\sum_{i=1}^j\beta_i$. Let $c=\sum_{i=1}^j\beta_i-\sum_{i=1}^j\gamma_i$. Consider the following equivalences of inequalities:
	\begin{align}
		\frac{\sum_{i=1}^j\gamma_i+c}{\sum_{i=1}^j\frac{\gamma_i}{b_i}+\frac{c}{b_j}}&\geq\frac{\sum_{i=1}^j\gamma_i}{\sum_{i=1}^j\frac{\gamma_i}{b_i}}\label{eq: deb}\\
		\pa{
			\sum_{i=1}^j\gamma_i+c
		}
		\frac{\pa{
				\sum_{i=1}^j\frac{\gamma_i}{b_i}
			}
		}{\sum_{i=1}^j\frac{\gamma_i}{b_i}+\frac{c}{b_j}}&\geq \sum_{i=1}^j\gamma_i \notag\\
		\pa{
			\sum_{i=1}^j\gamma_i+c-\frac{c\pa{
					\sum_{i=1}^j\gamma_i+c
				}
			}{b_j\sum_{i=1}^j\frac{\gamma_i}{b_i}+c}
		}
		&\geq \sum_{i=1}^j\gamma_i\notag\\
		\frac{\sum_{i=1}^j\gamma_i+c}{b_j\sum_{i=1}^j\frac{\gamma_i}{b_i}+c}&\leq 1\notag\\
		\sum_{i=1}^j\gamma_i&\leq b_j\sum_{i=1}^j\frac{\gamma_i}{b_i}\notag.
	\end{align}
	The last inequality holds since $b_1\leq\cdots\leq b_j$, and therefore \eqref{eq: deb} also holds. Since $\gamma_1,\dots,\gamma_{j-1},\gamma_j+c$ satisfy \eqref{eq: Cond3} and $(\sum_{i=1}^j\gamma_i)+c=\sum_{i=1}^j\beta_i$, we may conclude by \eqref{eq: deb} and the earlier case considered that
	\[
	\frac{\sum_{i=1}^j\gamma_i}{\sum_{i=1}^j\frac{\gamma_i}{b_i}}\leq \frac{\sum_{i=1}^j\gamma_i+c}{\sum_{i=1}^j\frac{\gamma_i}{b_i}+\frac{c}{b_j}}\leq\frac{\sum_{i=1}^j\beta_i}{\sum_{i=1}^j\frac{\beta_i}{b_i}}.\qedhere
	\]
\end{proof}
\printbibliography
\end{document}